\documentclass[11pt,reqno]{amsart}

\usepackage[margin=27mm]{geometry}
\usepackage[T1]{fontenc}
\usepackage{lmodern}
\usepackage{microtype}
\usepackage{mathtools}
\usepackage{amssymb}
\usepackage{mathrsfs}
\usepackage{booktabs}
\usepackage{array}
\usepackage{longtable}
\usepackage{alltt}
\usepackage{float}
\newfloat{listing}{tbp}{lol}
\floatname{listing}{Listing}
\newsavebox{\leanbox}
\usepackage{enumitem}
\usepackage[numbers,sort&compress]{natbib}
\usepackage{xcolor}
\definecolor{linkblue}{RGB}{0,0,150}
\usepackage[colorlinks,linkcolor=linkblue,citecolor=linkblue,urlcolor=linkblue]{hyperref}
\usepackage[nameinlink,noabbrev,capitalise]{cleveref}
\crefname{listing}{Listing}{Listings}

\providecommand{\doi}[1]{doi:~\href{https://doi.org/#1}{\nolinkurl{#1}}}
\numberwithin{equation}{section}
\makeatletter
\def\LT@makecaption#1#2#3{\LT@mcol\LT@cols c{\hbox to\z@{\hss\parbox[t]\LTcapwidth{\@makecaption{#2}{#3}\vskip\belowcaptionskip}\hss}}}
\makeatother
\setlist[enumerate]{label=(\alph*),leftmargin=*,itemsep=2pt}
\makeatletter
\def\section{\@startsection{section}{1}\z@{.7\linespacing\@plus\linespacing}{.5\linespacing}{\normalfont\bfseries\centering}}
\makeatother
\theoremstyle{plain}
\newtheorem{theorem}{Theorem}[section]
\newtheorem{lemma}[theorem]{Lemma}
\newtheorem{proposition}[theorem]{Proposition}
\newtheorem{corollary}[theorem]{Corollary}
\theoremstyle{definition}
\newtheorem{definition}[theorem]{Definition}
\theoremstyle{remark}

\newcommand{\R}{\mathbb{R}}
\newcommand{\C}{\mathbb{C}}
\newcommand{\Z}{\mathbb{Z}}
\newcommand{\Tor}{\mathbb{T}}
\newcommand{\one}{\mathbf{1}}
\newcommand{\T}{\mathsf{T}}
\newcommand{\cp}{\operatorname{cap}}

\newcommand{\bn}{\vec{\nu}}
\newcommand{\bt}{\mathbf{t}}
\newcommand{\Jm}{J_{-}}
\newcommand{\Area}{\mathcal{A}}
\newcommand{\Ih}{I_h}

\newcommand{\kc}{\varkappa}
\newcommand{\kt}{\tilde\varkappa}
\newcommand{\kh}{\varkappa_h}

\newcommand{\Lam}{\Lambda}
\newcommand{\Lp}{\Lambda^{\dagger}}
\newcommand{\Lh}{\Lambda_h}
\newcommand{\Lhp}{\Lambda_h^{\dagger}}
\newcommand{\SL}{\mathcal{V}}
\newcommand{\Pz}{P_0}
\newcommand{\avg}[1]{\langle #1\rangle_\Gamma}
\newcommand{\mV}{\mathsf{V}}
\newcommand{\mK}{\mathsf{K}}

\newcommand{\mB}{\mathsf{B}}
\newcommand{\mDs}{\mathsf{D}_s}

\newcommand{\Mv}{M_v}

\newcommand{\Lop}{L}
\newcommand{\Dop}{D}
\newcommand{\Kop}{K}
\newcommand{\Sop}{S}

\newcommand{\xnorm}[2][x]{\lVert #2\rVert_{#1}}
\newcommand{\dnorm}[1]{\lVert #1\rVert_*}
\newcommand{\wnorm}[1]{\lVert #1\rVert_{\mathcal W}}

\newcommand{\Om}{\Omega}
\newcommand{\ax}{a_x}
\newcommand{\bx}{b_x}
\newcommand{\aX}{a_X}
\newcommand{\bXh}{b_{X,h}}
\newcommand{\Ft}{F_\tau}
\newcommand{\Fht}{F_{h,\tau}}
\newcommand{\Nd}{\mathcal{N}_{r_0}}

\newcommand{\coupling}{h\le c_0\,\tau^2}

\newcommand{\dd}{\mathrm{d}}
\newcommand{\iu}{\mathrm{i}}
\newcommand{\eu}{\mathrm{e}}

\DeclareUrlCommand\lean{\urlstyle{tt}}

\title[Convergence of a parametric SGBEM for Mullins--Sekerka flow]{Convergence of a Parametric Symmetric Galerkin Boundary Element Method for the Mullins--Sekerka Flow}
\author{Tokuhiro Eto}
\address{Universit\'{e} Claude Bernard Lyon 1, CNRS, Centrale Lyon,
 INSA Lyon, Universit\'{e} Jean Monnet, ICJ UMR5208,
 69622 Villeurbanne, France}
\email{eto@math.univ-lyon1.fr}
\newif\ifaffiliationontitle
\affiliationontitlefalse
\ifaffiliationontitle
\makeatletter
\let\msII@setabstract\@setabstract
\renewcommand{\@setabstract}{%
 {\centering\small
 Universit\'{e} Claude Bernard Lyon 1, CNRS, Centrale Lyon, INSA Lyon,\\
 Universit\'{e} Jean Monnet, ICJ UMR5208, 69622 Villeurbanne, France\par
 \medskip
 \href{mailto:eto@math.univ-lyon1.fr}{\nolinkurl{eto@math.univ-lyon1.fr}}\par}%
 \bigskip
 \msII@setabstract}
\let\@setaddresses\relax
\makeatother
\fi
\date{}

\begin{document}

\begin{abstract}
We prove convergence of the fully discrete parametric scheme proposed by the author
\cite{E26-SGBEM} for the two-phase Mullins--Sekerka flow of a closed curve in the plane.
The scheme combines piecewise linear parametric curves, a curvature
discretization of Barrett--Garcke--N\"urnberg type and a symmetric Galerkin
boundary element approximation of the Dirichlet-to-Neumann map.
For a smooth solution whose curves satisfy a capacity condition, a uniform parameter mesh of size $h$,
a constant time step $\tau$ and exact boundary element integrals,
we show that the discrete curves stay regular and embedded,
and that their distance to the exact curves in the topology $W^{1,\infty}$, after reparametrization by arclength,
tends to zero in the first order of the time-step size $\tau$ whenever the mesh size $h$ is bounded by $c_0\tau^2$ for some $c_0$.
\end{abstract}

\subjclass[2020]{65M12 (primary); 65M38, 65M15, 35R35, 35R37, 53E10 (secondary).}
\keywords{Mullins--Sekerka flow, parametric finite element method, BGN-type
 scheme, symmetric Galerkin boundary element method, Dirichlet-to-Neumann
 map, convergence, error estimate}

\maketitle

\section{Introduction}\label{sec:intro}
\enlargethispage{3pt}

In the Mullins--Sekerka flow \cite{MS63}, an interface moves with a normal
velocity given by the jump of the normal derivative of a harmonic chemical
potential, whose boundary value is prescribed by the curvature. The velocity is obtained by
applying a nonlocal operator of order three to the curve. Parametric
discretizations in the spirit of Barrett, Garcke and N\"urnberg (BGN)
\cite{BGN07,BGN08,BGN20} are unconditionally stable for this flow. The
scheme of \cite{BGN10-2} conserves the enclosed areas exactly in the
semi-discrete setting \cite[Remark~110]{BGN20},
and the scheme of \cite[Theorem~3.3]{Nurnberg2023} and the variant of \cite[\S6]{EGN24} do so
in the fully discrete level.
In the scheme of \cite{E26-SGBEM}, the chemical potential is eliminated by a symmetric
Galerkin boundary element (SGBEM) approximation of the Dirichlet-to-Neumann
map, so that only the curve is discretized, and the scheme does not require the bulk mesh.
For every time step, the scheme decreases the length as long as the polygon stays admissible (\cref{sec:scheme}).

We call a numerical scheme \emph{BGN type} if its curve discretization
(i) uses a weak curvature identity with the normal and the arclength of the curve and lumped mass frozen at the previous time step,
(ii) prescribes only the normal velocity, so that the scheme determines the tangential motion,
and (iii) is linear and unconditionally stable.
How the normal velocity is computed does not matter.

For the curve shortening flow, fully discrete convergence of BGN-type schemes
has been proved only recently, for stabilized variants
\cite{BaiLi2025,Bai2026b} and for the classical scheme
\cite{Rao2026}, and also \cite{LiY2026} under an additional assumption on two auxiliary flows.
For curves moved by a prescribed velocity, a stabilized BGN-type scheme of transport type and degree $k\ge3$ converges
\cite{BGV2026Transport}.
Sharp-interface schemes for the Mullins--Sekerka flow have been developed in
several frameworks: \emph{boundary integral methods} in two and three space
dimensions \cite{BatesChenDeng1995,ZhuChenHou1996,BB00,Mayer2000,CKT17},
a \emph{level-set method} for Stefan problems with the Gibbs--Thomson condition
\cite{CMOS97}, the \emph{method of fundamental solutions} \cite{SY19,E24,E26},
the \emph{parametric finite element methods} mentioned above together with their
multi-phase versions \cite{EGN26,EGN26-2}, and \emph{oriented point-cloud
varifolds} for the one-phase flow \cite{EI26}.
For the one-phase Hele--Shaw flow with surface tension, which coincides with
the one-phase Mullins--Sekerka flow, Zhang, Hao and Harlim \cite{ZHH26} prove
the consistency and stability of a boundary integral discretization of the
normal velocity on point clouds at a fixed time, and the time evolution is
tested only numerically.
To our knowledge, no convergence result is known for the evolution computed
by a fully discrete sharp-interface scheme for the Mullins--Sekerka flow.
This paper proves convergence of the parametric SGBEM scheme for the Mullins--Sekerka flow proposed by the author \cite{E26-SGBEM}.
Though the scheme in \cite{E26-SGBEM} can be applied to the multi-phase case with anisotropic surface energy and triple junctions treated also in \cite{EGN26-2},
we do not address that setting here.

The two-phase Mullins--Sekerka flow describes
a single closed curve in $\R^2$ separating the whole space into two phases evolved by diffusion,
and it is governed by the surface evolution equation:
\begin{equation}\label{eq:MS-intro}
    \begin{cases}
        V = -\frac12\Lam_{\Gamma}\kc\qquad\text{on}\quad\Gamma(t),\quad t \in (0,T],\\
        \Gamma(0) = \Gamma_0,
    \end{cases}
\end{equation}
where $\Lam_{\Gamma}:H^{1/2}(\Gamma)\to H^{-1/2}(\Gamma)$ is the Dirichlet-to-Neumann operator;
$V$ denotes the velocity of the interface $\Gamma(t)$ in the direction of the outward normal vector;
$\varkappa$ is the curvature of $\Gamma(t)$ at each point (see \cref{eq:MS} in \cref{sec:setting} for the detail);
$T > 0$ is a time horizon;
$\Gamma_0$ is given as an initial interface.
We suppose that the proposed scheme produces polygons $x_h^m$ on a uniform parameter grid of size $h=1/N$, with time step
$\tau$, at the times $t_m=m\tau$ for $m = 0,\dots,\lfloor T/\tau\rfloor$,
where $\lfloor r\rfloor$ denotes the largest integer which does not exceed $r\in\R$.
We write $\hat x_h^m$ for $x_h^m$ reparametrized by normalized arclength with marked point $x_h^m(0)$, and $\Ih$ for nodal
interpolation.

Now, the main result of this paper can be stated as follows.
\begin{theorem}\label{thm:L2}
Let $\{\Gamma(t)\}_{t\in[0,T]}$ be a family of smooth embedded closed
curves solving \cref{eq:MS-intro}. Assume that $\Gamma(t)$ is parametrized by normalized arclength in terms of $\hat X\in C^\infty(\Tor\times[0,T];\R^2)$,
namely $\Gamma(t) = \hat X(\Tor, t)$ for $t\in[0,T]$, and that its marked point $\hat X(0,t)$ moves by \eqref{eq:marked-point} (\cref{sec:ms-flow}).
Assume that, for some $\eta_0>0$ and all $t\in[0,T]$, the curve $\Gamma(t)$ has a continuous equilibrium density and its logarithmic capacity is not larger than $1-\eta_0$, as stated in \ref{item:E-cap},
and let the proposed scheme be run with exact evaluation of all boundary element integrals,
on the uniform parameter grid of size $h=1/N$, from the nodal interpolant $x_h^0=\Ih\hat X(\cdot,0)$, with a constant time step $\tau$ (see \cref{sec:hyp}).

Then, there exist $\tau_0,\,c_0,\,C,\,C_q>0$ depending on the solution,
$T$ and $\eta_0$, such that for every $\tau \in (0,\tau_0]$ and
\begin{equation}\label{eq:coupling}
\coupling
\end{equation}
the following hold for $m=0,\dots,\lfloor T/\tau\rfloor$:
\begin{enumerate}
\item Every step \eqref{eq:step} is uniquely solvable, and every polygon
$x_h^m$ is regular, embedded and quasi-uniform with the constant $C_q$, that is, the length of each edge of $x_h^m$ is not larger than $C_q$ times the length of any other edge.
\item The error satisfies
\begin{equation}\label{eq:main}
\|\hat x_h^m-\hat X(\cdot,t_m)\|_{W^{1,\infty}}\le C\tau ,
\end{equation}
and the same bound holds for the Hausdorff distance between $x_h^m(\Tor)$
and $\Gamma(t_m)$.
\end{enumerate}
\end{theorem}

The proof of \cref{thm:L2} is given in \cref{sec:proof}.
The estimate \eqref{eq:main} is of first order in time. The coupling condition of $h$ and $\tau$
\eqref{eq:coupling} excludes the regime $\tau\sim h^2$ used in the
numerical experiments presented in \cite{E26-SGBEM}.
\cref{sec:proof} explains why the present argument does not reach this regime.

Convergence of parametric schemes for the curve shortening flow goes back to
Dziuk \cite{Dziuk1994} (see \cite{DDE2005} for a survey),
and recent results for two-phase Stokes flow \cite{BGV2026Stokes}.
The closest prior fully discrete analysis of a third-order nonlocal
interface flow is that of Ceniceros and Hou \cite[\S6]{CenicerosHou98}, for
spectrally accurate collocation in the formulation of Hou, Lowengrub and
Shelley \cite{HLS94} of a simplified two-fluid Hele--Shaw interface with
equal viscosities and densities, whose principal symbol $|k|^3$ is that of the Mullins--Sekerka
flow (see also \cite{AmbroseSiegelZhang2023} for Stokes flow, and
\cite{AmbroseLiuSiegel2017} for a point vortex method, discrete in space
only, for doubly periodic Darcy interfaces in three dimensions). Its proof
rests on exact discrete Fourier identities on equal-arclength grids and
does not apply to low-order Galerkin boundary elements on a polygon moved
by a BGN step.

For diffuse-interface approximations, convergence to the Mullins--Sekerka flow is known.
Fully discrete finite element
\cite{FengProhl2005} and discontinuous Galerkin \cite{FengLiXing2016}
approximations of the Cahn--Hilliard equation converge to the Mullins--Sekerka flow as the thickness of interface and the discretization parameters tend to zero,
provided that a classical solution exists, although the convergence of the interface is stated without a rate.

\subsection*{Strategy and organization}
We follow the time-first architecture of Rao \cite{Rao2026}, in which the fully discrete scheme is compared with the scheme discretized in time only. The proof has three steps.

(I) The time-semidiscrete scheme, continuous in space, converges at first order, and its iterates stay uniformly smooth (\cref{thm:semi-B,lem:semi-material}). This step uses no discretization in space, and from \cite{EtoPaperI} it uses only bounds for the single layer on chord-arc curves, namely its mapping property into $H^1$ and its Lipschitz dependence on the curve (\cref{sec:input-results}). The closure of the scheme is a recurrence for the BGN curvature whose implicit part has order three, and its dissipation pays for the derivative lost by the nonlocal velocity (\cref{sec:semi}). The proof has two stages. In the first stage, the defect of the recurrence is bounded only through the smoothing of one resolvent step, and the resulting coarse rate keeps the iterates in a fixed class of smooth curves, in which the continuum estimates hold with uniform constants (\cref{prop:semi-A}). In the second stage, uniform higher-order bounds (\cref{lem:semi-Hr}) make the defect of order $\tau$ and give the first order.

(II) The fully discrete analysis rests on shape-perturbation bounds for the SGBEM matrices on polygons close to the curves of a family of smooth curves that satisfies the hypotheses (K) of \cref{sec:input-K}. These bounds are proved in \cite{EtoPaperI}, and \cref{sec:bem} lists them in the form in which they are used. Our task is to show that the exact curves and the semidiscrete iterates form such a family (\cref{lem:Kfamily}). The uniform smoothness of step (I) gives \ref{K2}, and the hypothesis that needs an argument is the uniform coercivity \ref{K3} of the single layer, which follows from the capacity condition \ref{item:E-cap} and the $C^1$-closeness of the iterates to the exact curves.

(III) Each fully discrete step is a perturbation of the interpolated semidiscrete step. With the bounds of step (II), we control the error against the interpolated semidiscrete curve by the stability of the linearized semidiscrete step, the one-step consistency, the comparison of the linearized steps and the nonlinear remainder (\cref{sec:fd}). An induction bounds this error by $Ch^2/\tau$ (\cref{prop:superclose}), and the interpolation error and \cref{thm:semi-B} then give \eqref{eq:main}. The coupling \eqref{eq:coupling} comes from the loss $h\tau^{-1}+h^{1/2}$ in the comparison, and the nonlinear remainder requires only a weaker coupling (\cref{sec:proof}).

The argument differs from that of \cite{Rao2026} mainly in three points, which have no counterpart for the curve shortening flow. They treat the difficulties that come from the nonlocal velocity, in step (I), and from the area constraint and the SGBEM mass, in step (III).

First, in the step from the BGN curvature $\kt^m$ to the curvature $\kc^{m+1}$ of the new curve, which we call the metric step, the speed of the new parametrization contains two derivatives of $\kt^m$, because the normal velocity is $-\frac12\Lam_\Gamma\kt^m$, while for curve shortening it contains one. The metric step therefore loses one derivative more than in \cite{Rao2026}. We bound it in an additive form, in which the losses appear as separate terms with the factor $\tau$ (\cref{lem:semi-geom}), and the dissipation of the resolvent step of the same time step absorbs them (\cref{lem:semi-Hr} and the end of \cref{sec:semi}).

The main tool for this point is the splitting of the Dirichlet-to-Neumann map in normalized arclength, $\Lam_\Gamma=\Lam_0+Q_\Gamma$, where $\Lam_0=2\ell^{-1}|\partial_\xi|$ is a Fourier multiplier of order one and the remainder $Q_\Gamma$ has order $-1$ (\cref{thm:cont-N1}, proved in \cref{sec:app-cont}). The splitting separates the logarithmic kernel of the circle from the logarithm of the difference quotient of the curve, which loses only half a derivative of the curve. It gives bounds for $\Lam_\Gamma$ in which a high norm of the curvature appears only linearly (\cref{prop:cont-tame}), and this linearity closes the induction of \cref{lem:semi-Hr}.

Second, the multiplier of the area constraint is of order one, and curve shortening has no multiplier. It appears in two places, which are treated differently. In the stability estimate of the linearized semidiscrete step, it appears on both sides of the tangential row with the same sign and cancels exactly, before the division by $\tau$ would make it of order $\tau^{-1}$ (\cref{sec:app-p52-tan}). In the comparison of the linearized fully discrete and semidiscrete steps, the discrete multiplier is a vertex mean of the discrete curvature and the continuum one is an arclength mean, and they differ by $O(1)$, because the semidiscrete iterates are not parametrized by arclength. Hence, the multiplier terms cannot be compared separately. On each side, the variation of the multiplier term is grouped with the frame terms of the shape derivative of the mass that act on the test function, and the two groups reduce exactly to the discrete and the continuum curvature paired with the variation of the normal, which are then compared (\cref{lem:fd-cancel}).

Third, the mass of the scheme contains the SGBEM approximation of the Dirichlet-to-Neumann map, and its variations are controlled by bounds imported from \cite{EtoPaperI}. Since the increment of the mass under a perturbation $e$ of the polygon is bounded through $\|e'\|_{L^2}$ (\cref{thm:in-c1p}), the nonlinear remainder costs only $h^{-1/2}+\tau^{-1}$ (\cref{lem:fd-rem}).

\Cref{app:tools} states the results from the literature that the proofs use as external inputs, each with its source.
\Cref{app:new} proves the new arguments of the three points above and the coercivity used in step (II), and \cref{app:map} collects the precise statements of the remaining auxiliary results, whose proofs are supplied by the accompanying Lean~4 formalization. The estimates of Rao \cite{Rao2026} are quoted where they are used.

The main convergence theorem and its principal auxiliary results have been formalized in Lean~4. The verification is conditional on five explicitly identified external results from classical analysis and potential theory. \Cref{sec:formal} documents the formal statements, their correspondence with the mathematical results, and the remaining assumptions. The formalization builds on the independently developed verification of the shape-perturbation estimates in \cite{EtoPaperI}.

\section{Problem setting}\label{sec:setting}

\subsection{The Mullins--Sekerka flow}\label{sec:ms-flow}
Throughout the paper, a smooth curve is a closed curve $x:\Tor\to\R^2$, with $\Tor=\R/\Z$, of class $C^\infty$ that is regular, $x'\ne0$, and embedded, that is, injective on one period. We suppose that curves are oriented counter-clockwise, which we take to mean that their signed area is positive, with outward
unit normal $\bn$ and unit tangent $\bt$, so that $\partial_s\bt=-\kc\bn$ holds, where $s$ is the arclength parameter.
The sign convention of the curvature is chosen so that $\kc>0$ on convex curves.

Let $\Gamma=x(\Tor)$ be a smooth curve, and let $\Omega_{\rm int}$ and $\Omega_{\rm ext}$ be the bounded and the unbounded component of $\R^2\setminus\Gamma$. For a smooth function $g$ on $\Gamma$, let $u_{\rm int}$ be harmonic in $\Omega_{\rm int}$ and let $u_{\rm ext}$ be harmonic and bounded in $\Omega_{\rm ext}$, both continuous up to $\Gamma$ and equal to $g$ on $\Gamma$. The two-phase Dirichlet-to-Neumann map is
\[
\Lam_\Gamma g(y)=\lim_{\varepsilon\downarrow0}\frac{2g(y)-u_{\rm int}(y-\varepsilon\bn(y))-u_{\rm ext}(y+\varepsilon\bn(y))}{\varepsilon}\qquad(y\in\Gamma),
\]
that is, $\Lam_\Gamma g=\partial_{\bn}u_{\rm int}-\partial_{\bn}u_{\rm ext}$. In the parameter of the curve, this is the limit \eqref{eq:ext-lambda}, and \cref{prop:ext-dir}, with the uniqueness stated after it, shows that $u_{\rm int}$ and $u_{\rm ext}$ exist and are unique and that $\Lam_\Gamma g$ is a continuous function on $\Gamma$. On smooth functions, $\Lam_\Gamma$ is symmetric and
nonnegative in $L^2(\Gamma)$, and $\Lam_\Gamma1=0$ (see \cite[Proposition~3.1]{E26-SGBEM}).
Then, the two-phase isotropic Mullins--Sekerka flow is governed by the following law:
\begin{equation}\label{eq:MS}
    \begin{cases}
        V=-\frac12\Lam_{\Gamma}\kc\qquad\text{on}\quad\Gamma(t),\quad t\in(0,T],\\
        \Gamma(0) = \Gamma_0.
    \end{cases}
\end{equation}
We note that the above flow conserves the enclosed area and decreases the length of $\Gamma(t)$.

For a function $f$ on a curve, we define differential operators by
\[
L_f u := -\partial_s^2u+f^2u,\qquad D_f u:=(\partial_sf) u+2f\partial_s u.
\]
They are the operators of Rao \cite{Rao2026} that describe the tangential
motion selected by the BGN step. If a curve moves with the velocity
$v=V\bn+\beta\bt$, then $\partial_s\bt=-\kc\bn$ and
$\partial_s\bn=\kc\bt$ give
\[
\bt\cdot\partial_s^2v=-L_\kc\beta+D_\kc V .
\]
The following lemma shows that $L_f$ is invertible, with a uniform bound of the inverse, for coefficients $f$ that are bounded in $H^\sigma$ and whose $L^2$ norm is bounded below.

\begin{lemma}\label{lem:setting-L}
Let $\sigma>\frac32$, $M>0$, $\kc_{\min}>0$ and $0<\ell_{\min}\le\ell_{\max}$. Let $\Gamma$ be a smooth curve of length $\ell\in[\ell_{\min},\ell_{\max}]$, and let $f$ be a function on $\Gamma$ that, as a function of the normalized arclength $s/\ell\in\Tor$, satisfies $\|f\|_{H^\sigma(\Tor)}\le M$ and $\|f\|_{L^2(\Tor)}\ge\kc_{\min}$. Then, $L_f$ is an isomorphism from $H^1(\Gamma)$ onto its dual, and the norm of its inverse is bounded by a constant that depends only on $\sigma$, $M$, $\kc_{\min}$, $\ell_{\min}$ and $\ell_{\max}$.
\end{lemma}

\begin{proof}
Let $m_f=\int_\Gamma f^2\,\dd s=\ell\|f\|_{L^2(\Tor)}^2$, so that $\ell/m_f\le\kc_{\min}^{-2}$. For two points $s$ and $s'$ of $\Gamma$, we have $u(s')^2\le2u(s)^2+2\ell\|\partial_su\|_{L^2(\Gamma)}^2$ for $u\in H^1(\Gamma)$. We multiply by $f(s)^2$, integrate over $s$ and $s'$, and divide by $m_f$, which gives
\begin{equation}\label{eq:setting-L}
\|u\|_{L^2(\Gamma)}^2\le\frac{2\ell}{m_f}\int_\Gamma f^2u^2\,\dd s+2\ell^2\int_\Gamma|\partial_su|^2\,\dd s .
\end{equation}
Hence, the form $\int_\Gamma(L_fu)u\,\dd s=\int_\Gamma(|\partial_su|^2+f^2u^2)\,\dd s$ is coercive on $H^1(\Gamma)$, with a constant that depends only on $\kc_{\min}$ and $\ell_{\max}$. Since $\sigma>\frac12$, the Fourier series of $f$ converges absolutely, with $\wnorm{f}\le C_\sigma M$ by \cref{lem:tool-sobolev}, so that $|f|\le C_\sigma M$ almost everywhere and the form is also bounded, and the Lax--Milgram lemma gives the isomorphism and the bound of the inverse.
\end{proof}

By \cref{lem:setting-L} with $f=\kc$, whose $L^2$ norm is positive by \cref{lem:tool-frame}(d), we see that
\[
\alpha:=L_\kc^{-1}D_\kc V
\]
is the unique choice of $\beta$ for which $\bt\cdot\partial_s^2v=0$, and it is smooth by elliptic regularity if $\Gamma$ and $V$ are smooth. For the curve
shortening flow $V=-\kc$, we have $D_\kc V=-3\kc\partial_s\kc$, and $\alpha$ is the tangential velocity of \cite[(2.18)]{Rao2026}.
In the time-semidiscrete BGN step, the curvature equation forces the tangential component of
$\partial_s^2v$, taken on the old curve, to vanish (\cref{sec:semi}), so
the step selects the tangential velocity $L_\kc^{-1}D_\kc v_n$, and $\alpha$
is its limit as $\tau\to0$.

We parametrize $\Gamma(t)$ by normalized arclength, that
is, by $\hat X(\cdot,t):\Tor\to\R^2$ with $|\partial_\xi\hat X(\cdot,t)|=\ell(t)$,
the length of $\Gamma(t)$. This determines $\hat X(\cdot,t)$ up to the
choice of the point $\xi=0$, the marked point. As in \cite[(2.19)]{Rao2026},
we let the marked point move with the velocity of the BGN limit,
\begin{equation}\label{eq:marked-point}
\partial_t\hat X(0,t)=(V\bn+\alpha\bt)(0,t),\qquad t\in(0,T],
\end{equation}
and the other points then move tangentially keeping the relation $|\partial_\xi\hat X|=\ell$.
The discrete curves are compared with this parametrization, because the vertex $0$ of the polygon moves with the
tangential velocity of the BGN step.
For a regular parametrized curve $x$ of length $\ell_x$, let
$\phi_x(\xi)=\ell_x^{-1}\int_0^\xi|x'|$, which is increasing with
$\phi_x(\xi+1)=\phi_x(\xi)+1$, and let $\chi_x=\phi_x^{-1}$. We write
$\hat x=x\circ\chi_x$ for the normalized-arclength representative of $x$, so that
$\hat x(0)=x(0)$.
For a Lipschitz function $f:\Tor\to\R^k$, we set
$\|f\|_{W^{1,\infty}}=\|f\|_{L^\infty(\Tor)}+\|f'\|_{L^\infty(\Tor)}$, where $'=\dd/\dd\xi$, $|\cdot|$ is the Euclidean norm, and $L^\infty$ denotes the essential supremum. Two parametrized curves $x,y:\Tor\to\R^2$, either
smooth maps or the continuous piecewise linear maps that represent polygons, are
compared through $\|x-y\|_{W^{1,\infty}}$. This norm compares parametrizations, and \cref{thm:L2} applies it to the normalized-arclength representatives, which, for
regular embedded curves, depend only on the image, the orientation and the point
$x(0)$.

\subsection{The SGBEM scheme}\label{sec:scheme}
A polygon is identified by the image of a continuous piecewise linear map $X:\Tor\to\R^2$ on the
grid $\xi_j=jh\,(0\le j\le N)$ with the uniform mesh size $h=1/N$. It is equipped with
vertices $X_j$ and edges $\Gamma_p=[X_p,X_{p+1}]\,(0\le p\le N-1)$.
Here, since the polygon is closed, we use the conventions $X_0 = X_N$ and $X_{N+1}=X_1$. The edge $\Gamma_p$ is the image of the parameter interval $[\xi_p,\xi_{p+1}]$, which we call a panel, and
$\Ih$ denotes the nodal interpolation on this grid. The polygon is regular if every edge has positive length, and it is embedded if the piecewise linear map is injective on one period.
We identify a polygon, and a displacement $U$ of its vertices, with the
vector of its $N$ distinct vertices in $\R^{2N}$, so that $U=(U_1,\dots,U_N)$
with $U_j\in\R^2$. With
\[
\omega_j=\frac12\Jm(X_{j+1}-X_{j-1})\in\R^2 ,
\]
where $\Jm$ is the rotation by $-90^\circ$, let
$\Om(X)=\operatorname{blockdiag}(\omega_1^\T,\dots,\omega_N^\T)\in\R^{N\times2N}$
be the block diagonal matrix with the $1\times2$ blocks $\omega_j^\T$, so that
\[
\big(\Om(X)U\big)_j=\omega_j\cdot U_j,\qquad1\le j\le N.
\]
Let $\Area(X)=\frac12\sum_{p=0}^{N-1}X_p\wedge X_{p+1}$, with $a\wedge b=a_1b_2-a_2b_1$ for $a = (a_1,a_2)^\T$ and $b = (b_1,b_2)^\T$, be the
signed area of $X$. We write $\one\in\R^N$ for the vector with all entries one.
Then, $\one^\T\Om(X)U=D\Area(X)[U]$ holds for every displacement $U$,
with the same cyclic convention for $U$, since a shift of the summation index and $a\wedge b=a\cdot\Jm b$ give $D\Area(X)[U]=\frac12\sum_{j=1}^{N}U_j\wedge(X_{j+1}-X_{j-1})=\sum_{j}\omega_j\cdot U_j$.
For a regular polygon $X$, we let $A(X)$ be the stiffness matrix of the arclength Laplacian, so
that
\[
u^\T A(X)\varphi=\sum_{p = 0}^{N-1}\frac1{|\Gamma_p|}(u_{p+1}-u_p)\cdot(\varphi_{p+1}-\varphi_p)=:\aX(u,\varphi),
\]
and $A(X)X=\nabla\ell(X)$ for the length $\ell$. The lumped vertex mass is
$\Mv=\operatorname{diag}(\mu_j)$, $\mu_j=\tfrac12(|\Gamma_{j-1}|+|\Gamma_j|)$.

The single-layer and the double-layer operators of a curve $\Gamma$ are
\[
\SL g(y)=-\frac1{2\pi}\int_\Gamma\log|y-y'|\,g(y')\,\dd s_{y'},\qquad
\mathcal K g(y)=-\frac1{2\pi}\int_\Gamma\partial_{\bn(y')}\log|y-y'|\,g(y')\,\dd s_{y'},
\]
where $\partial_{\bn(y')}$ is the derivative in $y'$ along the normal $\bn(y')$. On a polygon, let $\varphi_j$ be the hat functions and $\iota_p$ the indicator of the edge $\Gamma_p$, both regarded as functions on $\Gamma=X(\Tor)$. The SGBEM Steklov--Poincar\'e matrix $\Lh(X)\in\R^{N\times N}$ of \cite[(5.19)]{E26-SGBEM} is formed from the Galerkin matrices
\[
\mV_{pq}=\int_{\Gamma_p}\SL\iota_q\,\dd s,\qquad\mK_{pj}=\int_{\Gamma_p}\mathcal K\varphi_j\,\dd s,\qquad\mB_{pj}=\int_{\Gamma_p}\varphi_j\,\dd s,\qquad(\mDs u)_p=\frac{u_{p+1}-u_p}{|\Gamma_p|},
\]
with all integrals evaluated exactly. With $Q_{\rm int}=\tfrac12\mB+\mK$, $Q_{\rm ext}=\tfrac12\mB-\mK$ and the Euclidean deflation $P=I-\one\one^\T/N$, it is
\[
\Lh(X)=P\big(2\mDs^\T\mV\mDs+Q_{\rm int}^\T\mV^{-1}Q_{\rm int}+Q_{\rm ext}^\T\mathsf C\,Q_{\rm ext}\big)P,
\]
where $\mathsf C$ is the upper left $N\times N$ block of the inverse of the bordered matrix with the blocks $\mV$, $\mB\one$, $(\mB\one)^\T$ and $0$, that is, the inverse of $\mV$ on the densities of arclength mean zero. The first term is the hypersingular operator in Maue's form, and the other two are the interior and the bounded exterior Schur complements. The matrices $\mV$ and the bordered matrix are invertible for the polygons considered in this paper \cite[Proposition~3.1]{EtoPaperI}, and the properties of $\Lh$ used here are collected in \cref{sec:bem}.
The matrix $\Lh(X)$ is the discrete counterpart of the two-phase
Dirichlet-to-Neumann map tested with the hat functions, in the sense of
\cref{lem:in-bweak}. Then, $\Lh$ is symmetric and $\Lh\one=0$, and we write $\Lhp$ for its
Moore--Penrose inverse, which maps into $\one^\perp=\{y\in\R^N:\one^\T y=0\}$.

In the two-phase reduced form \cite{E26-SGBEM}, one step maps $X^m$ to $X^{m+1}$
and a vector of vertex curvatures $\kh^{m+1}$ by
\begin{equation}\label{eq:step}
\Om(X^m)\big(X^{m+1}-X^m\big)=-\tfrac\tau2\Lh(X^m)\kh^{m+1},\qquad
\Om(X^m)^\T\kh^{m+1}=A(X^m)X^{m+1},
\end{equation}
where $\kh^{m+1}$ is minus the curvature variable of \cite{E26-SGBEM}. The unknown of the step is the pair $(X^{m+1},\kh^{m+1})$, and we say that the step is uniquely solvable if this pair is unique.
Since $\Lh\one=0$, the first equation of \eqref{eq:step} gives
$D\Area(X^m)[U]=\one^\T\Om(X^m)U=0$ for $U=X^{m+1}-X^m$, so the area constraint holds at the velocity level, and $\Area(X^{m+1})-\Area(X^m)=\Area(U)=O(|U|^2)$ because $\Area$ is
quadratic \cite[Propositions~6.6--6.7]{E26-SGBEM}.
If the polygon $X^m$ is non-degenerate and the region it encloses has logarithmic capacity less than one,
every solution of \eqref{eq:step} satisfies
$\ell(X^{m+1})\le\ell(X^m)$ for every $\tau>0$
\cite[Proposition~5.6, Theorem~6.4]{E26-SGBEM}, and if $X^m$ is moreover
solvability-admissible \cite[Definition~6.9]{E26-SGBEM},
the step is uniquely solvable for every $\tau>0$
\cite[Theorem~6.10]{E26-SGBEM}.

The analysis uses the step in a variational form. With the discrete mass
\[
\bXh(u,\varphi)=2\big(\Om(X)\varphi\big)^\T\Lhp(X)\big(\Om(X)u\big),
\]
we consider the problem of finding $X^{m+1}$ and a multiplier $c_h\in\R$ with $\one^\T\Om(X^m)(X^{m+1}-X^m)=0$ and
\begin{equation}\label{eq:step-var}
a_{X^m}(X^{m+1},\varphi)+\tau^{-1}b_{X^m,h}(X^{m+1}-X^m,\varphi)=c_h\one^\T\Om(X^m)\varphi\qquad\text{for all }\varphi .
\end{equation}
Assume that $\ker\Lh(X^m)=\operatorname{span}\{\one\}$, which holds for the polygons considered in this paper by \cite[Lemma~A.16]{EtoPaperI}. Then, $\Lh\Lhp=P$, and the two problems correspond as follows. If $(X^{m+1},\kh^{m+1})$ solves \eqref{eq:step}, then $(X^{m+1},c_h)$ solves \eqref{eq:step-var}, where $c_h$ is the vertex mean of $\kh^{m+1}$, that is, the Euclidean mean of its
entries. Conversely, if $(X^{m+1},c_h)$ solves \eqref{eq:step-var}, then $(X^{m+1},\kh^{m+1})$ solves \eqref{eq:step} with
\[
\kh^{m+1}=c_h\one-\frac2\tau\Lhp\Om(X^m)\big(X^{m+1}-X^m\big).
\]
Indeed, the first equation of \eqref{eq:step} and $\Lh\Lhp=P$ determine $\kh^{m+1}$ up to an element of $\ker\Lh$, that is, up to a multiple of $\one$, and the second equation, tested with $\varphi$, becomes \eqref{eq:step-var}. In particular, if \eqref{eq:step-var} has a unique solution $(X^{m+1},c_h)$, then the step \eqref{eq:step} is uniquely solvable. The unique solvability in \cref{thm:L2} is obtained in this way from \cref{cor:bem-acoerc}.
We write $\Fht(X^m)=X^{m+1}$. When the step is applied to a general
polygon $X$, as in \cref{sec:fd}, we drop the time index and write $\kh$. The iterates of the scheme are the polygons $x_h^m=X^m$ with $X^0=\Ih\hat X(\cdot,0)$ and $X^{m+1}=\Fht(X^m)$ for $0\le m<\lfloor T/\tau\rfloor$.

The time-semidiscrete scheme is the same step without spatial discretization. The unknown is a curve $x^m:\Tor\to\R^2$ instead of a polygon.
Write $\Gamma^m=x^m(\Tor)$, $\kc^m$ for its curvature and $\bn^m$
for its normal. Given $x^m$, find $x^{m+1}=x^m+\tau v^m$ and $\kt^m$ on
$\Gamma^m$ with
\begin{equation}\label{eq:semi}
v^m\cdot\bn^m=-\frac12\Lam_{\Gamma^m}\kt^m,\qquad
\partial_s^2x^{m+1}=-\kt^m\bn^m .
\end{equation}
Here, $\partial_s$ and $\dd s$ are those of $\Gamma^m$, all functions are written in the parameter of $x^m$, and
the second equation holds in the weak sense
$\int\kt^m\bn^m\cdot\varphi\,\dd s=\int\partial_sx^{m+1}\cdot\partial_s\varphi\,\dd s$ for every smooth $\varphi:\Tor\to\R^2$.
The solution is sought among the pairs for which $x^{m+1}$ is a smooth curve and $\kt^m$ is a smooth function.
The function $\kt^m$ is the curvature of the new curve $x^{m+1}$ measured with the
arclength and the normal of $\Gamma^m$, and it differs from the curvature $\kc^{m+1}$ of $\Gamma^{m+1}$ by $O(\tau)$. We call the parameter of $x^m$ the material parameter. The curvatures $\kc^m$, $\kt^m$ and the components of $v^m$ are written either in the normalized arclength of $\Gamma^m$ or in the material parameter of $x^m$, as stated where they are used. The two forms are related by composition with $\chi_{x^m}$ and $\phi_{x^m}$. Integrals with $\dd s$ are taken in the arclength of the curve, and integrals with $\dd\xi$ in the parameter, with $\dd s=|x'|\,\dd\xi$ and $\bn\,\dd s=\Jm x'\,\dd\xi$ for a curve $x$.
We write $\Ft(x^m)=x^{m+1}$, start from $x^0=\hat X(\cdot,0)$, and take the steps with $(m+1)\tau\le T$. Let $\avg{g} := \tfrac{1}{|\Gamma|}\int_\Gamma g\,\dd s$, where $|\Gamma|>0$ is the length of $\Gamma$, and let $\Pz g=g-\avg g$ be the
$L^2(\dd s)$-orthogonal projection onto functions of arclength mean zero, and $\Lp=\Pz\SL\Pz$
the mean-zero single layer. It is the pseudo-inverse of the two-phase
Dirichlet-to-Neumann map on functions of arclength mean zero
(\cref{lem:cont-gauge}), which is the reason that the single layer appears in the
mass. With the forms
\[
\ax(u,\varphi)=\int\partial_su\cdot\partial_s\varphi\,\dd s,\qquad
\bx(u,\varphi)=2\int\Lp(u\cdot\bn)\,(\varphi\cdot\bn)\,\dd s,
\]
\eqref{eq:semi} is the continuum version of \eqref{eq:step-var}, with
multiplier $c^m=\langle\kt^m\rangle_{\Gamma^m}$, the arclength mean.

\subsection{Hypotheses}\label{sec:hyp}

We split the hypotheses of \cref{thm:L2} into two parts, (E) on the solution and (S) on the scheme.
Let $T>0$ and $\eta_0>0$. Hypothesis (E) consists of the following three conditions.
\begin{enumerate}[label=(E\arabic*),leftmargin=*]
\item\label{item:E-smooth} The map $\hat X$ belongs to $C^\infty(\Tor\times[0,T];\R^2)$, with one-sided derivatives in $t$ at $t=0$ and $t=T$. For every $t\in[0,T]$, the curve $\hat X(\cdot,t)$ is a smooth curve oriented counter-clockwise, and $|\partial_\xi\hat X(\cdot,t)|=\ell(t)$ holds, where $\ell(t)$ is its length.
\item\label{item:E-flow} For every $t\in(0,T]$, the normal velocity $V=\partial_t\hat X\cdot\bn$ of $\Gamma(t)=\hat X(\Tor,t)$ satisfies \eqref{eq:MS}, and the marked point moves by \eqref{eq:marked-point}.
\item\label{item:E-cap} For every $t\in[0,T]$, there are a continuous function $w_{\rm eq}$ on $\Gamma=\Gamma(t)$ and a number $\lambda_{\rm eq}=\lambda_{\rm eq}(\Gamma)$ such that
\[
\int_\Gamma w_{\rm eq}\,\dd s=1,\qquad\SL w_{\rm eq}=\lambda_{\rm eq}\ \text{on }\Gamma,\qquad \eu^{-2\pi\lambda_{\rm eq}}\le1-\eta_0 .
\]
\end{enumerate}
Hypothesis (S) is that the scheme is run with exact integration on the uniform grid of size $h=1/N$, from
the start $x_h^0=\Ih\hat X(\cdot,0)$ and with a constant time step $\tau$, as described in \cref{sec:scheme}.

The function $w_{\rm eq}$ of \ref{item:E-cap} is unique. Indeed, if $w_1$ and $w_2$ are two such functions, with constants $\lambda_1$ and $\lambda_2$, then $\vartheta=w_1-w_2$ satisfies $\int_\Gamma\vartheta\,\dd s=0$ and $(\SL\vartheta,\vartheta)_{L^2(\Gamma)}=(\lambda_1-\lambda_2)\int_\Gamma\vartheta\,\dd s=0$, so that $\vartheta=0$ by \cref{prop:ext-slcoer}. It is the equilibrium density of $\Gamma$, and
$\cp\Gamma=\eu^{-2\pi\lambda_{\rm eq}(\Gamma)}$ is the logarithmic capacity of
$\Gamma$ \cite[\S6.6.1]{Steinbach2008}, \cite[Ch.~8]{McLean2000}, \cite[Ch.~I]{SaffTotik}. Thus, \ref{item:E-cap} states that $\cp\Gamma(t)$ is not larger than $1-\eta_0$ for $t\in[0,T]$. By \cref{lem:K3-one}, the single layer is then coercive on $\Gamma(t)$. The
capacity condition is a normalization of scale. Rescaling adds a constant to the logarithmic kernel, which changes the
single-layer matrix by a rank-one term, so the discrete scheme need not be
invariant under rescaling, and the scheme is applied to the curve rescaled
so that the condition holds. Smooth
solutions exist for smooth initial curves
\cite{ChenHongYi1996,EscherSimonett1996MRL,ES98,PrussSimonett2016}.

\section{The time-semidiscrete scheme}\label{sec:semi}

This section treats the scheme discretized in time only. Its iterates are compared with the polygons in \cref{sec:fd}, and that comparison differentiates the boundary element operators along the semidiscrete curves and uses interpolation estimates, so the main output of this section is uniform smoothness of the iterates, besides the first-order rate. Rao \cite{Rao2026} obtains both for curve shortening from a recurrence for the BGN curvature $\kt^m$ whose implicit part has order two. For the Mullins--Sekerka flow, the normal velocity is minus half the Dirichlet-to-Neumann map applied to the curvature, one derivative rougher than in curve shortening, and the implicit part $\tfrac12\Sop\Lam$ of the recurrence, with $\Sop$ defined below, has order three. The end of this section explains how the dissipation of this third-order part pays for the derivative that the change to the new arclength loses.

With the operators $L_f$, $D_f$ of \cref{sec:setting}, Rao
\cite{Rao2026} sets
\[
K_f=D_fL_f^{-1}D_f,\qquad \Sop_f=L_f+K_f .
\]

The nonlinearity of the Mullins--Sekerka recurrence is built from the
tangential velocity that preserves normalized arclength. For a smooth
closed curve $\Gamma$ of length $\ell$, a function $f$ on $\Gamma$ and
\[
V=-\tfrac12\Lam_\Gamma f,\qquad \alpha=L_f^{-1}D_fV,\qquad
\mathfrak P[g](\xi)=\ell\int_0^\xi\big(g-\langle g\rangle_\Gamma\big)\,\dd\xi',
\]
let
\[
\mathcal T_f(\xi)=\alpha(0)-\ell\int_0^\xi\big(fV-\langle fV\rangle_\Gamma\big)\,\dd\xi',\qquad
N_{\rm MS}(f;\Gamma)=-2f^2V-K_fV+\mathcal T_f\,\partial_sf .
\]
The nonlinearity $N_{\rm MS}$ has order one in $f$, through $K_fV$, $f^2V$ and $\mathcal T_f\partial_sf$.
Since $K_fV=\alpha\,\partial_sf+2f\,\partial_s\alpha$, it reads
\[
N_{\rm MS}(f;\Gamma)=-2f\,(\partial_s\alpha+fV)+(\mathcal T_f-\alpha)\,\partial_sf .
\]
If $f$ is the curvature $\kc$ of $\Gamma$, then $\partial_s\alpha+\kc V=\bt\cdot\partial_sv$
for $v=V\bn+\alpha\bt$. Thus, for $f=\kc$, the nonlinearity is the stretching of the curve plus a transport
relative to the arclength-preserving motion.

The functions $\kc^m$, $\kt^m$ are written in the normalized arclength $\xi\in\Tor$ of their own
curve $\Gamma^m$, so that the recurrence below and \cref{thm:semi-B}
compare functions on $\Tor$. In normalized arclength, $\partial_s=\ell^{-1}\partial_\xi$, and the two-phase Dirichlet-to-Neumann map splits as $\Lam_\Gamma=\Lam_0+Q_\Gamma$ with $\Lam_0=2\ell^{-1}|\partial_\xi|$ and a remainder $Q_\Gamma$ of order $-1$ (\cref{thm:cont-N1}). Here, $|\partial_\xi|\eu^{2\pi\iu n\xi}=2\pi|n|\eu^{2\pi\iu n\xi}$. Since $-\frac12\ell^{-2}\partial_\xi^2\Lam_0=\ell^{-3}|\partial_\xi|^3$, we can write
\begin{equation}\label{eq:semi-PE}
\frac12\Sop_f\Lam_\Gamma=\mathcal P_\ell+\mathcal E_{f,\Gamma},\qquad \mathcal P_\ell=\ell^{-3}|\partial_\xi|^3,\qquad
\mathcal E_{f,\Gamma}=\frac12(f^2+K_f)\Lam_\Gamma-\frac12\ell^{-2}\partial_\xi^2Q_\Gamma .
\end{equation}
The principal part $\mathcal P_\ell$ is the Fourier multiplier $(2\pi|n|/\ell)^3$. We write $\mathcal P^m=\mathcal P_{\ell^m}$, $\mathcal E^m=\mathcal E_{\kc^m,\Gamma^m}$, $\Lam^m=\Lam_{\Gamma^m}$, and $v^m=v_n^m\bn^m+v_t^m\bt^m$ for the velocity of the step. The following lemma is the counterpart of \cite[(3.4), (3.26)]{Rao2026} for the Mullins--Sekerka flow. The nonlinearity $N_{\rm MS}$ is written here for a general coefficient $f$, with $V=-\frac12\Lam_\Gamma f$.

\begin{lemma}[The recurrence for the BGN curvatures]\label{prop:closure}
Let $x^m$ be a smooth, regular and embedded curve, and let $(x^{m+1},\kt^m)$ be a smooth solution of
\eqref{eq:semi} with $x^{m+1}$ regular and embedded. Then, we have
\begin{equation}\label{eq:normal-tangential}
\kt^m=\kc^m+\tau\,\Sop_{\kc^m}v_n^m,\qquad v_t^m=L_{\kc^m}^{-1}D_{\kc^m}v_n^m,
\qquad v_n^m=-\frac12\Lam^m\kt^m ,
\end{equation}
so that $\kt^m$ solves the resolvent equation
\begin{equation}\label{eq:semi-reseq}
\kt^m+\tau\big(\mathcal P^m+\mathcal E^m\big)\kt^m=\kc^m .
\end{equation}
If $(x^{m+2},\kt^{m+1})$ is a smooth solution of the next step, then the BGN curvatures obey the recurrence
\begin{equation}\label{eq:recurrence}
\big(I+\tau\mathcal P^m\big)\kt^{m+1}=\kt^m+\tau\Big(N_{\rm MS}(\kt^m;\Gamma^m)-\mathcal E^m\kt^m+d^m\Big),
\end{equation}
where the defect $d^m$, defined in \eqref{eq:semi-drows}, collects the remainder of the metric step and the changes of the operators and of their arguments between the levels $m$ and $m+1$. Under (E), the exact flow satisfies
\begin{equation}\label{eq:exact-curv}
\partial_t\kc=-\big(\mathcal P_\ell+\mathcal E_{\kc,\Gamma}\big)\kc+N_{\rm MS}(\kc;\Gamma)\qquad\text{on }\Tor\times[0,T].
\end{equation}
\end{lemma}
\begin{proof}
The defect \eqref{eq:semi-drows} is derived in \cref{sec:app-semi-geom}, and \eqref{eq:exact-curv} is \cref{lem:semi-exact}(c). For a general curve, \eqref{eq:semi-reseq} is \eqref{eq:normal-tangential} rewritten with \eqref{eq:semi-PE}, and \eqref{eq:recurrence} follows from the algebra of \cref{sec:app-semi-geom}. The remaining details are supplied by the formalization (\cref{sec:formal}), in which \eqref{eq:semi-reseq} and \eqref{eq:recurrence} are stated for curves of $\mathcal U_{\rm ad}$ and small $\tau$, the setting of all their uses.
\end{proof}

The identities of \cref{prop:closure} come from the frame equations of the
step (\cref{sec:app-semi-geom}). The tangential and normal components of
$\partial_s^2x^{m+1}=-\kt^m\bn$ are the two rows of \eqref{eq:normal-tangential}.
By \eqref{eq:semi-newcurv}, the true curvature $\kc^{m+1}$ is $\kt^m$ times
a factor $1+O(\tau)$ built from $\partial_sv$, followed by the change to normalized arclength.
This is the \emph{metric step}. The change of
parametrization produces the transport term
$\mathfrak P[\partial_sv_t+\kc v_n]\,\partial_s\kt$, which the
identity \eqref{eq:semi-cancel} turns into the term $\mathcal T\partial_s\kt$
of $N_{\rm MS}$. Inserting \eqref{eq:semi-reseq} at level $m+1$ gives the recurrence, and each of the five terms of $d^m$ is $\tau$ times a norm of $\kt^m$ or $\kt^{m+1}$, or a change of the coefficients in one step.

The continuum facts needed about $\Lam_\Gamma$ are collected in
\cref{sec:app-cont}. Besides the splitting $\Lam_\Gamma=\Lam_0+Q_\Gamma$, they are tame bounds for $\Lam_\Gamma$ and $Q_\Gamma$, in which a high norm of the curvature appears only linearly and the constants depend only on lower norms (\cref{prop:cont-tame}), and the Lipschitz dependence of $\Lam_\Gamma$ and $Q_\Gamma$ on the curve when the
curvature is measured half a derivative below the image space (\cref{thm:cont-N3,thm:cont-N3plus}). The induction on the order in \cref{lem:semi-Hr} needs this linearity, as explained at the end of this section.

Let $\sigma\in(5/2,3)$ be fixed. The lower bound gives
$\bn\cdot\partial_sv=\partial_sv_n-\kc v_t\in H^{\sigma-2}\subset L^\infty$, $\bt\cdot\partial_sv\in H^{\sigma-1}\subset C^1$
and $\kt^m\in H^\sigma\subset C^2$. These bounds keep the step regular, make the remainder of the metric step pointwise of size $O(\tau^2)$, and keep the new curve $C^1$-close to the old one. The operator $\mathcal E_{f,\Gamma}$ has order one, but in the present class only its part $\frac12(f^2+K_f)\Lam_\Gamma$ maps $H^\sigma$ into $H^{\sigma-1}$. The part $-\frac12\ell^{-2}\partial_\xi^2Q_\Gamma$ is bounded from $H^\sigma$ only into $H^{\sigma-3/2}$, because the tame bound for $Q_\Gamma$ loses half a derivative of the curvature (\cref{lem:semi-E}). The origin of this loss in the kernel is explained before \cref{thm:cont-N1}. Consequently, one resolvent step smooths $\kc^m\in H^\sigma$ only into $H^{\sigma+3/2}$, and the kernels $(t_m-t_j)^{-\gamma_i/3}$ that carry the forcings into $H^\sigma$ have the exponents $\gamma_i\le\sigma$. The upper bound $\sigma<3$ makes these kernels integrable (\cref{lem:semi-stab}).

\begin{theorem}\label{thm:semi-B}
Assume that the hypothesis (E) is valid. Then, there exist $\tau_B>0$ and $C$ such that, for $0<\tau\le\tau_B$, every step
\eqref{eq:semi} with $(m+1)\tau\le T$ has a unique smooth, regular and
embedded solution, and we have
\[
\|\kc^m-\kc(t_m)\|_{H^\sigma}+|\ell^m-\ell(t_m)|+\|\hat x^m-\hat X(t_m)\|_{C^1}\le C\tau\qquad(0\le m\tau\le T),
\]
and $\|\kt^m-\kc(t_m)\|_{H^\sigma}\le C\tau$ for $(m+1)\tau\le T$.
\end{theorem}

The $C^1$ norm compares the normalized-arclength representatives in the fixed coordinates of the plane, at every point of $\Tor$ and without a rigid motion.

\begin{lemma}\label{lem:semi-material}
Assume that the hypothesis (E) is valid. For every $r\in\mathbb N$, there exist $\tau_s(r)>0$ and
$C_r^{\rm s}$ such that, for $0<\tau\le\tau_s(r)$, the iterates of \cref{thm:semi-B} satisfy
\[
\|x^m\|_{W^{r+1,\infty}}\le C_r^{\rm s},\qquad
(C_r^{\rm s})^{-1}\le|\partial_\xi x^m|\le C_r^{\rm s}\qquad(0\le m\tau\le T),
\]
and, with $v^m=(x^{m+1}-x^m)/\tau$ and the multiplier $c^m=\langle\kt^m\rangle_{\Gamma^m}$ of \cref{sec:scheme},
\[
\|v^m\|_{W^{r,\infty}}+\|\kt^m\|_{H^r}+|c^m|\le C_r^{\rm s}\qquad((m+1)\tau\le T).
\]
\end{lemma}

In \cref{lem:semi-material}, $x^m$, $v^m$ and $\kt^m$ are functions of the material parameter, the parameter of the scheme \eqref{eq:semi} started from $x^0=\hat X(\cdot,0)$, and not of the normalized arclength.

The proof of \cref{thm:semi-B} has two stages, \cref{prop:semi-A} and \cref{lem:semi-Hr}, which are combined in \cref{sec:app-semi-first}. The defect
$d^m$ is $\tau$ times high norms of $\kt^m$. In the first stage, these norms
are bounded only through the smoothing of one resolvent step, which gives
$\|d^m\|_{L^2}\le C\tau^{1/3}$, and the induction on $m$ closes at a
coarse rate $\tau^{(3-\sigma)/6}$ (\cref{prop:semi-A}). In each step of the induction, the new curve is placed in the admissible class before the estimates of \cref{sec:app-cont} are applied to it, as described in step~(ii) after \cref{lem:Kfamily}. The second stage proves
uniform higher-order bounds (\cref{lem:semi-Hr}), with which $\|d^m\|_{L^2}=O(\tau)$, and the same stability argument gives the first order. The $C^1$ bound in the fixed coordinates is recovered from the marked point and its tangent, and the material bounds of \cref{lem:semi-material} follow from the normalized-arclength bounds through the chain of reparametrizations.

The argument is Rao's, and the following part of one step has no counterpart for curve
shortening. Each step consists of a resolvent step $\kc^m\mapsto\kt^m$,
which gains one and a half derivatives in the mean from the third-order
operator $\tfrac12\Sop\Lam$, see \eqref{eq:semi-res}, and the metric step
$\kt^m\mapsto\kc^{m+1}$. In \cite[Lemmas~4.2 and 4.3]{Rao2026}, the speed
$\zeta=|\bt+\tau\partial_sv|$ of the update and the reparametrization
are of lower order. Here,
$\partial_\xi\zeta=-\tau\ell\kt\,(\bn\cdot\partial_sv)/\zeta$, and
$\bn\cdot\partial_sv$ contains $\partial_s\Lam\kt$, two derivatives of $\kt$,
whereas for curve shortening it contains $\partial_s\kt$. Hence, the metric step costs
$C\tau\|\kt^m\|_{H^{r+1}}$ in $H^r$, and, through the tame bound for $\Lam$, also $C\tau\|\kc^m\|_{H^r}$. We bound the metric step in the additive form \eqref{eq:semi-metric}, in which both losses appear as separate terms with the factor $\tau$ and not in the constants. The loss in $\kt^m$ is paid for by the
$H^{r+3/2}$ dissipation of the same step, which controls half a derivative more than the loss, and Young's inequality. The loss in $\kc^m$ is linear by \cref{prop:cont-tame}, through \cref{lem:semi-vel}(b), and it is absorbed by the discrete Gronwall inequality, because the constants of the metric step do not depend on $\|\kc^m\|_{H^r}$ (\cref{lem:semi-geom}(e) and \cref{lem:semi-Hr}).

\section{The boundary element input}\label{sec:bem}

The fully discrete analysis of \cref{sec:fd} needs bounds for the SGBEM matrix $\Lh(X)$ on polygons $X$ close to the iterate $x^m$, uniformly in $m$ and $\tau$. Such bounds are proved in \cite{EtoPaperI} for every family of smooth curves that satisfies the hypotheses (K1)--(K3) below, and \cref{lem:Kfamily} shows that the exact curves and the time-semidiscrete iterates form such a family. We then fix the neighborhood of polygons, list the results of \cite{EtoPaperI} in the form in which they are used, and derive from them four results on one step of the scheme and on the pseudo-inverse.

\subsection{The hypotheses on the reference curves}\label{sec:input-K}
We identify $\R^2$ with $\C$, and $d$ denotes the distance on $\Tor$. For a closed Lipschitz curve $Y:\Tor\to\R^2$ with a positive chord-arc constant, the parameter form of the single layer is
\[
\mathcal L_Y\varsigma(\xi)=-\frac1{2\pi}\int_\Tor\log|Y(\xi)-Y(\eta)|\,\varsigma(\eta)\,\dd\eta ,
\]
and the parameter form of the double layer is
\[
\mathcal K_Ya(\xi)=-\frac1{2\pi}\,\mathrm{p.v.}\!\int_\Tor\mathrm{Im}\,\frac{Y'(\eta)}{Y(\eta)-Y(\xi)}\,a(\eta)\,\dd\eta ,
\]
where the principal value is the limit as $\delta\to0$ of the integrals over $\{\eta\mid d(\xi,\eta)>\delta\}$ \cite[Section~2]{EtoPaperI}. On a polygon, $\mathcal L_X\varsigma$ and $\mathcal K_Xa$ are the operators $\SL$ and $\mathcal K$ of \cref{sec:scheme} applied to the arclength densities $\varsigma/|X'|$ and to $a$, written in the parameter. For $M\ge c>0$, let $\mathfrak F(M,c)$ be the set of closed Lipschitz curves $Y$ with $|Y'|\le M$ almost everywhere and $|Y(\xi)-Y(\eta)|\ge c\,d(\xi,\eta)$, as in \cite{EtoPaperI}. If $\varsigma$ is the parameter density of an arclength density $w$ on $\Gamma=Y(\Tor)$, defined by $\langle\varsigma,v\circ Y\rangle=(w,v)_{L^2(\Gamma)}$, then $\langle\mathcal L_Y\varsigma,\varsigma\rangle=(\SL w,w)_{L^2(\Gamma)}$. All Sobolev norms are taken on $\Tor$ in the parameter. The following hypotheses are those of \cite[Definition~2.1]{EtoPaperI}.

\begin{definition}[The hypotheses (K)]\label{def:K}
Let $\mathscr K$ be a family of smooth embedded closed curves $x:\Tor\to\R^2$, oriented counter-clockwise, for which the following constants are uniform on $\mathscr K$.
\begin{enumerate}[label=(K\arabic*)]
\item\label{K1} The curve is uniformly chord-arc, $|x'|\le M_*$ and $|x(\xi)-x(\eta)|\ge c_*\,d(\xi,\eta)$.
\item\label{K2} For the order $r\ge0$ named in each result, $\|x\|_{W^{r+1,\infty}}\le C_r$.
\item\label{K3} The single layer is uniformly coercive on $L^2$ densities, $\langle\mathcal L_x\varsigma,\varsigma\rangle\ge\mu_1\|\varsigma\|_{H^{-1/2}}^2$ for all $\varsigma\in L^2$ with $\mu_1>0$. It maps $L^2$ boundedly onto $H^1$, and $\|\varsigma\|_{L^2}\le B_1\|\mathcal L_x\varsigma\|_{H^1}$ for all $\varsigma\in L^2$.
\end{enumerate}
The numbers $M_*$, $c_*$, $\mu_1$, $B_1$ and $C_r$ are called the constants of $\mathscr K$.
\end{definition}

The script letter $\mathscr K$ distinguishes the family from the double layer $\mathcal K_Y$. The next lemma is the only place where the boundary element input uses the capacity condition of (E) and the semidiscrete estimates.

\begin{lemma}[{The reference curves}]\label{lem:Kfamily}
Assume (E), and let $R\subset\mathbb N$ be finite. There is a threshold $\tau_V>0$, which we may take not larger than $\tau_B$ and not larger than the thresholds $\tau_s(r)$ of \cref{lem:semi-material} for $r\in R\cup\{1\}$, such that the family
\[
\mathscr K=\{\hat X(\cdot,t)\mid 0\le t\le T\}\cup\{x^m\mid 0<\tau\le\tau_V,\ 0\le m\tau\le T\}
\]
satisfies \ref{K1}, \ref{K3}, and \ref{K2} for every $r\in R$. The threshold $\tau_V$ and the constants of $\mathscr K$ depend only on (E), on $R$ and on the constants of \cref{thm:semi-B,lem:semi-material}. In particular, they do not depend on $m$ or $\tau$.
\end{lemma}

In the rest of the paper, $R=\{1,2,3\}$, which is the set of orders $r$ for which \ref{K2} is invoked, and $\mathscr K$ is the family of \cref{lem:Kfamily}. The hypotheses \ref{K1} and \ref{K2} follow from the smoothness of the exact curves, from \cref{lem:semi-material} and, for the chord-arc bound of the iterates in \ref{K1}, from the $C^1$ bound of \cref{thm:semi-B}. The argument for \ref{K3} has the following four steps, whose lemmas are listed in \cref{sec:app-cap,sec:app-cont-K3}.

(i) For each exact curve, \ref{K3} follows from the capacity condition, through the splitting of a density into a multiple of the equilibrium density and a part of mean zero, together with the coercivity on densities of mean zero and the Fredholm property of the single layer (\cref{prop:ext-slcoer,prop:ext-slfred} and \cref{lem:K3-one}). The coercivity and the bound of the inverse in \ref{K3} survive a small perturbation of the curve in $C^1$, by the Lipschitz dependence of the single layer on the curve \cite[Lemma~A.7]{EtoPaperI} (\cref{lem:K3-pert}). By compactness in $t$, this makes the constants uniform on $[0,T]$, with the constants $\mu_0$ and $B_0$ of \cref{lem:K3-exact}.

(ii) The same argument gives the coercivity constant $\frac12\mu_0$ to every smooth curve of constant speed whose derivative is, after a rigid motion, close to that of an exact curve (\cref{lem:cont-near}). This step does not use the iterates, so the induction of \cref{prop:semi-A} can use it to place each new curve in the admissible class before the estimates of \cref{sec:app-cont} are applied to that curve (\cref{lem:semi-close}), without a circular argument.

(iii) Once \cref{thm:semi-B} is proved, the normalized-arclength representative $\hat x^m$ is $C^1$-close to $\hat X(\cdot,t_m)$, and \cref{lem:K3-pert}, together with \cref{lem:K3-one}(a) for the surjectivity, carries \ref{K3} from $\hat X(\cdot,t_m)$ to $\hat x^m$. Then, \ref{K3} passes from $\hat x^m$ to $x^m=\hat x^m\circ\phi_{x^m}$ under the change of parameter (\cref{lem:K3-reparam}).

(iv) The constants of steps (i)--(iii) depend only on (E) and on the constants of \cref{thm:semi-B,lem:semi-material}, so $\mathscr K$ satisfies \ref{K3} uniformly in $m$ and $\tau\le\tau_V$.

\subsection{The neighborhood of polygons}\label{sec:input-nbhd}
From now on, $\tau\le\tau_V$. Let $C_I$ be the constant in $\|(\Ih x-x)'\|_{L^\infty}\le C_Ih$. We take the radii $r_0,h_0>0$ of \cite{EtoPaperI}, made smaller if necessary so that, with $r_1=r_0+C_Ih_0$, the conditions \cite[(2.3)]{EtoPaperI} of \cite[Definition~2.3]{EtoPaperI} and the condition \eqref{eq:use-small} of \cref{cor:bem-acoerc} hold for the constants of $\mathscr K$, and so that the stability bound \eqref{eq:use-brezzi} of \cref{cor:bem-acoerc} holds. For $x\in\mathscr K$ and $0<r\le r_0$, let $\mathcal N_r=\mathcal N_r(x)$ be the set of polygons $X$ with
\[
\|X-\Ih x\|_{W^{1,\infty}}\le r,\qquad h\le h_0,
\]
and let $\mathfrak C$ be the set of Lipschitz curves $Y$ with $\|Y'-x'\|_{L^\infty}\le2r_1$. Then, $\Nd\subset\mathfrak C$, and every $Y\in\mathfrak C$ is uniformly chord-arc,
\[
|Y'|\le M_1,\qquad |Y(\xi)-Y(\eta)|\ge c_{\rm ca}\,d(\xi,\eta),\qquad M_1=M_*+2r_1,\quad c_{\rm ca}=c_*-2r_1\ge\tfrac12c_* ,
\]
so that $|Y'|\ge c_{\rm ca}$ almost everywhere. Both sets are convex, so $\Nd$ contains a segment $[X,X+e]$ as soon as it contains its end points, and $\mathfrak C$ contains $x+\varepsilon e$ whenever $|\varepsilon|\,\|e'\|_{L^\infty}\le2r_1$. In the rest of the paper, $x=x^m$ and $\Nd=\Nd(x^m)$, and the radii $r_0$, $h_0$, $r_1$ and the constants $C$ and $c$ depend only on the constants of $\mathscr K$, so they are uniform in $m$ and in $\tau\le\tau_V$.

On $\Tor$, $\mathcal P_1$ denotes the continuous piecewise linear functions, with the hat functions $\varphi_j$, and $\mathcal P_0$ the piecewise constant functions on the grid $\xi_j=jh$. A vector of vertex values, also in loads $\Om(X)f$ and in $\bXh(f,\varphi)$, is measured through its piecewise linear interpolant. For a load vector $y\in\R^N$, we use the dual norms
\[
\dnorm{y}=\sup_{w\in\mathcal P_1\setminus\{0\}}\frac{y^\T w}{\|w\|_{H^{1/2}}},\qquad
\|y\|_{(H^1)^*}=\sup_{w\in\mathcal P_1\setminus\{0\}}\frac{y^\T w}{\|w\|_{H^1}} .
\]
For P1 fields $f$ and $\varphi$, we set $z_f=\Lhp(X)\Om(X)f$ and $z_\varphi=\Lhp(X)\Om(X)\varphi$, which we call discrete Neumann solutions, and for $w\in H^{1/2}$, we write $\|w\|_{H^{1/2}/\R}=\inf_{\beta\in\R}\|w-\beta\|_{H^{1/2}}$.

\subsection{\texorpdfstring{Boundary element bounds from \cite{EtoPaperI}}{Boundary element bounds from the earlier paper}}\label{sec:input-results}
The main input of the estimates of \cref{sec:fd} is the following proposition, which restates \cite[Theorem~3.3]{EtoPaperI}, \cite[Theorem~3.4]{EtoPaperI} and \cite[Theorem~3.6]{EtoPaperI}. It bounds the first variation of the mass, its second-order remainder, and the consistency error of $\Lh$ on the interpolating polygon. Parts (a) and (b) are used in the bound of the nonlinear remainder (\cref{lem:fd-rem}), and part (c) is used in the one-step consistency (\cref{lem:fd-cons}) and in the comparison of the linearized steps (\cref{lem:fd-comp}).

\begin{proposition}\label{prop:in-import}
Assume \ref{K1} and \ref{K3}, and let $h\le h_0$.
\begin{enumerate}[label=(\alph*),ref=\theproposition(\alph*)]
\item\label[proposition]{thm:in-c1p} Assume \ref{K2} with $r=1$. Let $X\in\Nd$, and let $e$, $f$ and $\varphi$ be piecewise linear. Then,
\begin{equation}\label{eq:in-C1}
\big|z_\varphi^\T D_X\Lh(X)[e]\,z_f\big|\le C\,\|e'\|_{L^2}\,\|f\|_{L^2}\,\|\varphi\|_{L^2}.
\end{equation}
If moreover the segment $[X,X+e]$ lies in $\Nd$, then
\[
\big|b_{X+e,h}(f,\varphi)-\bXh(f,\varphi)\big|\le C\,\|e'\|_{L^2}\big(\|f\|_{L^\infty}\|\varphi\|_{L^2}+\|f\|_{L^2}\|\varphi\|_{L^\infty}\big).
\]
\item\label[proposition]{thm:in-c3p} Assume \ref{K2} with $r=1$. Let $X$, $e$, $f$ and $\varphi$ be piecewise linear, and let the segment $[X,X+e]$ lie in $\Nd$. Then, the second-order Taylor remainder of $X\mapsto\bXh(f,\varphi)$ in the direction $e$,
\[
\delta^2\bXh(f,\varphi)=b_{X+e,h}(f,\varphi)-\bXh(f,\varphi)-D_X\bXh(f,\varphi)[e],
\]
satisfies
\begin{equation}\label{eq:in-C3}
|\delta^2\bXh(f,\varphi)|\le C\,\|e'\|_{L^\infty}\,\|e'\|_{L^2}\,\|f\|_{L^\infty}\,\|\varphi\|_{L^2}.
\end{equation}
\item\label[proposition]{lem:in-bweak} Let $X=\Ih x$, and let $\mathsf r=\Lh(X)\Ih g-\Mv\Ih(\Lam_\Gamma g)$, where $\Gamma=x(\Tor)$. If \ref{K2} holds with $r=2$, then
\[
\dnorm{\mathsf r}\le Ch\|g\|_{H^2}\qquad(g\in H^2),
\]
and if \ref{K2} holds with $r=3$, then
\[
\dnorm{\mathsf r}\le Ch^{3/2}\|g\|_{H^3},\qquad\|\mathsf r\|_{(H^1)^*}\le Ch^2\|g\|_{H^3}\qquad(g\in H^3).
\]
\end{enumerate}
\end{proposition}

Besides \cref{prop:in-import}, we use the following results of \cite{EtoPaperI}, in three groups: the continuity of the boundary integral operators in the curve, the gauge representation together with the coercivity of the single layer, and the discrete regularity on $\Nd$. They are proved there under \ref{K1}--\ref{K3}, where \ref{K2} with $r=1$ is used only through the definition of $\Nd$, and under the further hypotheses stated with each of them. The external inputs of \cite{EtoPaperI} are restated in \cref{sec:tools-ext,map:inputs}.

\subsubsection*{Continuity in the curve}

For every polygon $X\in\mathfrak C$, the formula of \cref{sec:scheme} for $\Lh(X)$ can be written in the parameter \cite[Proposition~3.1]{EtoPaperI}. In this form, $\Lh(X)$ depends on $X$ only through the Galerkin matrices of $\mathcal L_X$ and $\mathcal K_X$ on the panels. By the log form of \cite[Lemma~3.2]{EtoPaperI}, these matrices are twice differentiable in the vertices, so $\Lh$ is differentiable in the vertices.

For a direction $e$ and a curve $X$ with a positive chord-arc constant, $q_e(\xi,\eta)=(e(\xi)-e(\eta))/(X(\xi)-X(\eta))$, $\xi\ne\eta$, is the divided-difference quotient of $e$ relative to $X$. Let $M\ge c>0$ and $X,Y\in\mathfrak F(M,c)$. The Lipschitz dependence of \cite[Lemma~A.7]{EtoPaperI} states that, for $\varsigma\in L^2$,
\[
\|(\mathcal L_X-\mathcal L_Y)\varsigma\|_{H^1}\le C\|X'-Y'\|_{L^\infty}\|\varsigma\|_{L^2},\qquad\|(\mathcal L_X-\mathcal L_Y)\varsigma\|_{H^{1/2}}\le C\|X'-Y'\|_{L^\infty}\|\varsigma\|_{H^{-1/2}}.
\]
By the same lemma, for a Lipschitz direction $e$ and $\varsigma,\psi\in L^2$, the function $\varepsilon\mapsto\langle\mathcal L_{X+\varepsilon e}\varsigma,\psi\rangle$ is continuously differentiable near $\varepsilon=0$, and its derivative at $\varepsilon=0$ is $\langle f_{\mathcal L},\psi\rangle$, written as an integral over $\Tor\times\Tor$, where
\[
f_{\mathcal L}(\xi)=-\frac1{2\pi}\int_\Tor\operatorname{Re}q_e(\xi,\eta)\,\varsigma(\eta)\,\dd\eta .
\]
On the same classes, $\|\partial_\xi\mathcal L_X\varsigma\|_{L^2}\le C\|\varsigma\|_{L^2}$ for real $\varsigma\in L^2$ \cite[Lemma~A.5]{EtoPaperI}. In this paragraph, $C$ depends only on $M$ and $c$, and \ref{K1}--\ref{K3} are not needed. The two lemmas quoted here are stated in \cite{EtoPaperI} for arbitrary $M\ge c>0$, and their proofs use neither a family of reference curves nor \ref{K1}--\ref{K3}. At the curve $x$, with $q_e$ relative to $x$, the function $f_{\mathcal L}$ represents the shape derivative of the single layer in the sense of \cite[Lemma~A.24]{EtoPaperI}.

\subsubsection*{Gauge representation and coercivity}

The operators are uniformly invertible on $\mathfrak C$ \cite[Lemma~A.8]{EtoPaperI}. For all $X\in\mathfrak C$ and $\varsigma\in L^2$, we have $\langle\mathcal L_X\varsigma,\varsigma\rangle\ge c\|\varsigma\|_{H^{-1/2}}^2$ and
\[
C^{-1}\|\varsigma\|_{L^2}\le\|\mathcal L_X\varsigma\|_{H^1}\le C\|\varsigma\|_{L^2},\qquad C^{-1}\|\varsigma\|_{H^{-1/2}}\le\|\mathcal L_X\varsigma\|_{H^{1/2}}\le C\|\varsigma\|_{H^{-1/2}},
\]
and $\mathcal L_X$ maps $L^2$ onto $H^1$. The constants depend only on $M_*$, $c_*$, $\mu_1$, $B_1$ and, through $r_1$, on $C_1$. The gauge representation of \cite[Lemma~A.10]{EtoPaperI} gives, for $Y\in\mathfrak C$ and $g\in H^1$, a unique $(\varsigma,c)\in L^2\times\R$ with $\mathcal L_Y\varsigma=g+c$ and $\int\varsigma\,\dd\xi=0$, and it satisfies $\|\varsigma\|_{L^2}+|c|\le C\|g\|_{H^1}$. Let $Y$ be $x$ or $x+\varepsilon e$ with a piecewise linear $e$ and $|\varepsilon|\,\|e'\|_{L^\infty}\le2r_1$. Its image $\Gamma_Y=Y(\Tor)$ is a counter-clockwise Jordan curve that bounds a Lipschitz domain \cite[Proposition~A.9]{EtoPaperI}, and $\Lam_Yg=(\Lam_{\Gamma_Y}(g\circ Y^{-1}))\circ Y$ denotes its two-phase Dirichlet-to-Neumann map written in the parameter, where $\Lam_{\Gamma_Y}$ is defined as in \cref{sec:ms-flow} with the normal derivatives taken in the sense of traces on the Lipschitz domain, the exterior function being bounded. If $g\in H^2$, then $\Lam_Yg=\varsigma/|Y'|$ with the gauge density $\varsigma$ of $g$ on $Y$ \cite[Proposition~A.9]{EtoPaperI}.

\subsubsection*{Discrete regularity on $\Nd$}

By the discrete regularity of \cite[Lemma~A.20]{EtoPaperI}, for $X\in\Nd$ and the load $y$ of $f$, either consistent, $y_j=\int f\varphi_j\,|X'|\,\dd\xi$ with $f\in L^2$, or lumped, $y_j=\mu_jf(\xi_j)$ with $f$ P1, we have $\|\Lhp y\|_{H^1}\le C\|f\|_{L^2}$. By \cite[Lemma~A.21]{EtoPaperI}, for $X\in\Nd$, a piecewise linear direction $e$ and P1 functions $u$ and $w$ whose vertex values sum to zero, we have $|w^\T D_X\Lh(X)[e]\,u|\le C\|e'\|_{L^2}\|u\|_{H^1}\|w\|_{H^1}$. Finally, the shape derivative of $\Lh$ is consistent with that of the Dirichlet-to-Neumann map \cite[Theorem~3.7]{EtoPaperI}. Assume \ref{K2} with $r=2$, let $X=\Ih x$ with $h\le h_0$, $g\in H^2$ and $e$ piecewise linear, and let
\[
\mathsf r_j=\big(D\Lh(X)[e]\,\Ih g\big)_j-\frac{\dd}{\dd\varepsilon}\int_\Tor(\Lam_{x+\varepsilon e}\,g)\,\varphi_j\,|x'+\varepsilon e'|\,\dd\xi\Big|_{\varepsilon=0},
\]
where $\Lam_{x+\varepsilon e}$ is the two-phase Dirichlet-to-Neumann map of the curve $x+\varepsilon e$ and the derivative exists by \cite[(3.5)]{EtoPaperI} and \cite[Lemma~A.23]{EtoPaperI}. Then, we have $\|\mathsf r\|_{(H^1)^*}\le Ch\|e'\|_{L^\infty}\|g\|_{H^2}$ and $\dnorm{\mathsf r}\le Ch^{1/2}\|e'\|_{L^\infty}\|g\|_{H^2}$.

\subsection{One step of the scheme and the consistency of the pseudo-inverse}\label{sec:input-step}
The four results of this subsection follow from those of \cref{sec:input-results}. The proofs of \cref{lem:in-ellip,lem:in-loads} are short and follow their statements, and the proofs of \cref{cor:bem-acoerc,lem:bem-invcons} are given in \cref{app:bem-step}.

\begin{lemma}[Uniform ellipticity of $\Lh$]\label{lem:in-ellip}
Assume \ref{K1}, \ref{K2} with $r=1$, and \ref{K3}. Let $X\in\Nd$ and $h\le h_0$. Then, for all P1 functions $w$,
\begin{equation}\label{eq:in-ellip}
c\,\|w\|_{H^{1/2}/\R}^2\le w^\T\Lh(X)w\le C\|w\|_{H^{1/2}/\R}^2 .
\end{equation}
For every $y\in\R^N$ with $\one^\T y=0$, we have
\begin{equation}\label{eq:in-dual}
c\,\dnorm{y}^2\le y^\T\Lhp(X)y\le C\dnorm{y}^2,\qquad\|\Lhp(X)y\|_{H^{1/2}/\R}\le C\dnorm{y},
\end{equation}
and for every $y\in\R^N$, we have
\begin{equation}\label{eq:in-Lhp}
\|\Lhp(X)y\|_{H^{1/2}}\le C\dnorm{y}.
\end{equation}
\end{lemma}

\begin{proof}
For P1 functions, $w^\T\Lh(X)w$ is the discrete Schur form of \cite{EtoPaperI}, which vanishes on constants and satisfies the upper bound of \cite[Lemma~A.15]{EtoPaperI} and the lower bound on functions of mean zero recalled at the beginning of \cref{app:bem-step}. Applied to $w-\beta$ with constants $\beta$, these bounds give \eqref{eq:in-ellip}. For $\one^\T y=0$ and $u=\Lhp y$, we have $\Lh u=y$, since $\ker\Lh=\operatorname{span}\{\one\}$ \cite[Lemma~A.16]{EtoPaperI}, and \eqref{eq:in-dual} follows from $y^\T\Lhp y=u^\T\Lh u$, \eqref{eq:in-ellip}, the definition of $\dnorm{\cdot}$ and the Cauchy--Schwarz inequality for the positive semidefinite matrix $\Lh$. Finally, \eqref{eq:in-Lhp} is \eqref{eq:in-dual} applied to $Py$, because $\Lhp y=\Lhp Py$, $\dnorm{Py}\le2\dnorm{y}$, since $|\one^\T y|\le\dnorm{y}$ and $\dnorm{\one}\le N$, and $\Lhp Py$ has mean zero.
\end{proof}

\begin{lemma}[The constraint loads]\label{lem:in-loads}
For every polygon $X$ and every piecewise linear $\varphi$, the constraint functional is an integral in the parameter,
\begin{equation}\label{eq:use-area}
\one^\T\Om(X)\varphi=\int_\Tor\varphi\cdot\Jm X'\,\dd\xi .
\end{equation}
Assume \ref{K1}, \ref{K2} with $r=1$, and \ref{K3}, and let $X\in\Nd$ and $h\le h_0$. Then, for every piecewise linear $\varphi$ and every $y\in\R^N$,
\begin{equation}\label{eq:use-Om}
\dnorm{\Om(X)\varphi}\le C\|\varphi\|_{L^2},
\end{equation}
\begin{equation}\label{eq:use-load}
\big|(\Om(X)\varphi)^\T\Lhp(X)\,y\big|\le C\|\varphi\|_{L^2}\,\dnorm{y}.
\end{equation}
In particular, $|\bXh(u,\varphi)|\le C\|u\|_{L^2}\|\varphi\|_{L^2}$ for piecewise linear $u$ and $\varphi$.
\end{lemma}

\begin{proof}
On the panel $[\xi_p,\xi_{p+1}]$, the derivative $X'$ is constant, so the panel contributes $\tfrac12\Jm(X_{p+1}-X_p)\cdot(\varphi(\xi_p)+\varphi(\xi_{p+1}))$ to the right side of \eqref{eq:use-area}, and collecting the terms with $\varphi(\xi_j)$ gives $\omega_j\cdot\varphi(\xi_j)$. For \eqref{eq:use-Om}, the load $\Om(X)\varphi$ is the lumped load $\Mv\psi$ of a P1 function $\psi$ with $\|\psi\|_{L^2}\le C\|\varphi\|_{L^2}$ \cite[Appendix~A.6]{EtoPaperI}, and $|w^\T\Mv\psi|\le C\|\psi\|_{L^2}\|w\|_{L^2}$ \cite[(A.22)]{EtoPaperI}. Then, \eqref{eq:use-load} follows from \eqref{eq:in-Lhp} and \eqref{eq:use-Om}, and the bound of $\bXh$ is \eqref{eq:use-load} with $y=\Om(X)u$.
\end{proof}

\begin{corollary}[{Coercivity, unique solvability and stability on $\Nd$}]\label{cor:bem-acoerc}
Assume \ref{K1}, \ref{K2} with $r=1$, and \ref{K3}. There is a constant $C_\sharp$, depending only on $M_*$, $c_*$ and $C_1$, such that the following holds if, in addition to the conditions \cite[(2.3)]{EtoPaperI} of \cite[Definition~2.3]{EtoPaperI},
\begin{equation}\label{eq:use-small}
C_\sharp\,(r_1+h_0)\le\frac{c_*^2}{M_*}.
\end{equation}
Let $X\in\Nd$, $h\le h_0$ and $0<\tau\le1$. Then,
\[
\aX(u,u)+\tau^{-1}\bXh(u,u)\ge c\,\|u\|_{H^1}^2\qquad\text{for all piecewise linear }u\text{ with }\one^\T\Om(X)u=0,
\]
and $\sup_{\varphi\ne0}\one^\T\Om(X)\varphi/\|\varphi\|_{H^1}\ge c$, where the supremum is taken over the nonzero piecewise linear $\varphi$. More precisely, $\varphi_*=\Ih\bn_x$, the interpolant of the unit normal of $x$, satisfies $\one^\T\Om(X)\varphi_*\ge c_*/2$ and $\|\varphi_*\|_{H^1}\le C$ with $C$ depending only on $M_*$, $c_*$ and $C_1$. Here, $c>0$ depends only on the constants of \ref{K1}--\ref{K3}, on $C_1$ and $C_I$, and on $r_0$ and $h_0$. Consequently, \eqref{eq:step-var} with $X^m=X$ has a unique solution $(X^{m+1},c_h)$, so $\Fht$ is defined on $\Nd$. Moreover, after $r_0$ and $h_0$ are made smaller if necessary, depending only on $M_*$, $c_*$, $\mu_1$, $B_1$ and $C_1$, the following holds. If a piecewise linear $U$ and numbers $\mathsf R,\delta\ge0$ satisfy $|\aX(U,\varphi)+\tau^{-1}\bXh(U,\varphi)|\le\mathsf R\|\varphi\|_{H^1}$ for all piecewise linear $\varphi$ with $\one^\T\Om(X)\varphi=0$, and $|\one^\T\Om(X)U|\le\delta$, then
\begin{equation}\label{eq:use-brezzi}
\|U\|_{H^1}\le C\big(\mathsf R+\tau^{-1}\delta\big).
\end{equation}
\end{corollary}

\begin{lemma}[Consistency of $\Lhp$ with $\Lp$]\label{lem:bem-invcons}
Assume \ref{K1}, \ref{K2} with $r=1$, and \ref{K3}. Let $X=\Ih x$ with $h\le h_0$, and let $y$ be the load of a function $f$, either consistent, $y_j=\int f\varphi_j\,\rho\,\dd\xi$ with $f\in L^2$ and $\rho\in\{|X'|,|x'|\}$, or lumped, $y=M_vf$ with $f$ a P1 function. If $\one^\T y=0$, then there is a constant $\bar c$ with $|\bar c|\le C\|f\|_{L^2}$ such that
\[
\|\Lhp y-\Lp f-\bar c\|_{H^{1/2}}\le Ch^{1/2}\|f\|_{L^2},
\]
where $\Lp f$ is the continuum operator $\Lam^\dagger_\Gamma$ of the curve applied to $f$ as a function of the parameter.
\end{lemma}

\section{The fully discrete scheme and the proof of the main theorem}\label{sec:fd}

Let $e^m=x_h^m-\Ih x^m$ be the error against
the interpolated semidiscrete curve. One step of the error is the
difference of two steps from nearby starting points, the polygon step
$\Fht$ and the semidiscrete step $\Ft$,
\begin{equation}\label{eq:err-split}
e^{m+1}=\underbrace{\Fht(\Ih x^m)-\Ih\Ft(x^m)}_{\text{consistency}}
+\underbrace{D\Fht(\Ih x^m)e^m}_{\text{linear part}}+\underbrace{\mathcal R^m}_{\text{remainder}} .
\end{equation}
The first term is the error of one fully discrete step from the interpolated semidiscrete curve (\cref{lem:fd-cons}). The second is the fully discrete linearized step applied to the error, which we compare with the interpolated semidiscrete linearized step (\cref{lem:fd-comp}). The third is the nonlinear remainder of the fully discrete step (\cref{lem:fd-rem}).

Let $x=x^m$ be a time-semidiscrete iterate with $(m+1)\tau\le T$, so that $\Ft(x^m)=x^{m+1}$, and let $\tau\le\tau_V$ as in
\cref{lem:Kfamily}. In the statements of this section, $\tau$ and $h$ are
not larger than thresholds that depend only on the data of \cref{thm:L2},
and the constants $C$ depend only on the same data. The linear part is estimated through the stability of the semidiscrete
step, which is measured in the adapted norm
\[
\xnorm{e}^2=\|e\cdot\bn\|_{H^1(\Gamma)}^2+\|e\cdot\bt-L_\kc^{-1}D_\kc(e\cdot\bn)\|_{H^1(\Gamma)}^2,
\]
with $\kc$ the curvature of $\Gamma=x(\Tor)$ and the norms taken in
normalized arclength. This is a variant of the norm \cite[(5.8)]{Rao2026}.
It satisfies the triangle inequality and is uniformly equivalent to the norm of
$H^1(\Tor;\R^2)$ (\cref{lem:p52-adnorm}). The next lemma adapts \cite[Prop.~5.2]{Rao2026}. Its new step is the cancellation of the multiplier (\cref{sec:app-p52-tan}).

\begin{lemma}[Stability of the linearized semidiscrete step]\label{prop:stab}
There are $\tau_P>0$ and $C$ such that, for $0<\tau\le\tau_P$ and
$(m+1)\tau\le T$, the derivative $D\Ft(x^m)$ extends from smooth directions to a bounded operator on $H^1(\Tor;\R^2)$, and
\begin{equation}\label{eq:four}
\xnorm[x^{m+1}]{D\Ft(x^m)e}\le(1+C\tau)\,\xnorm[x^m]{e}\qquad\text{for all }e\in H^1(\Tor;\R^2).
\end{equation}
\end{lemma}

Its proof is given at the end of \cref{sec:app-p52-tan}, after the cancellation of the multiplier. The tangential row is Rao's once the multiplier has cancelled, so the tangential part $\eta(e)=e\cdot\bt-L_\kc^{-1}D_\kc(e\cdot\bn)$ of the adapted norm changes by $O(\tau)\|e\|_{H^1}$. The normal row, after the mobility $2\SL$ is inverted with $\tfrac12\Lam_\Gamma$, is one resolvent step for the third-order operator $\mathcal H=\tfrac12\Lam_\Gamma\Sop$ on normal velocities, whose principal part is $\ell^{-3}|\partial_\xi|^3$ in normalized arclength, with a forcing of order $\tau\|e\|_{H^1}$ in $L^2$. One step dissipates $c\tau\|g\|_{H^{5/2}}^2$ in $H^1$, and this absorbs the forcing, because pairing an $L^2$ forcing with $g$ in $H^1$ costs $\|g\|_{H^2}$.

The next three lemmas adapt the stiffness part of \cite[Lemma~5.1, (5.2)--(5.4)]{Rao2026}, \cite[Lemma~5.3]{Rao2026} and \cite[Lemma~6.1]{Rao2026}, respectively. Their new steps are the SGBEM mass in the first, the grouping of the multiplier and the frame terms (\cref{lem:fd-cancel}) in the second, and the second variation of the SGBEM mass in the third.

\begin{lemma}[One-step consistency]\label{lem:fd-cons}
Let $x$ be the time-semidiscrete iterate fixed above, and let $0<h\le\tau\le1$. Then, the consistency
error of one fully discrete step satisfies
\[
\|\Fht(\Ih x)-\Ih\Ft(x)\|_{H^1}\le Ch^2 ,
\]
with $C$ independent of $\tau$.
\end{lemma}

\begin{lemma}[Comparison of the linearized steps]\label{lem:fd-comp}
Let $e_h$ be piecewise linear, and let $0<h\le\tau\le1$. Then, we have
\[
\|D\Fht(\Ih x)e_h-\Ih D\Ft(x)e_h\|_{H^1}\le C\big(h\tau^{-1}+h^{1/2}\big)\|e_h\|_{H^1}.
\]
\end{lemma}

\begin{lemma}[Nonlinear remainder]\label{lem:fd-rem}
There is $r_2\in(0,r_0]$, depending only on the constants of $\mathscr K$, such that the following holds. If $0<h\le\tau\le1$ and the segment $[\Ih x,\Ih x+e_h]$
lies in $\mathcal N_{r_2}$, then we have
\[
\big\|\Fht(\Ih x+e_h)-\Fht(\Ih x)-D\Fht(\Ih x)e_h\big\|_{H^1}\le C\big(h^{-1/2}+\tau^{-1}\big)\|e_h\|_{H^1}^2 .
\]
\end{lemma}

\Cref{lem:fd-cons,lem:fd-comp,lem:fd-rem} are proved in \cref{sec:app-fd}. The step contains $\tau^{-1}\bXh$, the mass divided by $\tau$, so an error in the mass is amplified by $\tau^{-1}$. The consistency has no loss in $\tau$, because the boundary element defect
appears through $\tau^{-1}\bXh$ applied to an increment of size $\tau$, and it is of order $h^2$ in $(H^1)^*$
(\cref{lem:in-bweak}) when it is paired with the discrete pseudo-inverse applied to the load of a test function, which is bounded in $H^1$. In the comparison, $h\tau^{-1}$ is the $O(h)$ error of the
mass amplified by $\tau^{-1}$, as in Rao's analysis, and also the interpolation error of the continuum derivative, whose regular part has second derivative of order $\tau^{-1}$ by \eqref{eq:fd-z}. The term $h^{1/2}$
comes from the consistency of the shape derivative of $\Lh$ (\cite[Theorem~3.7]{EtoPaperI}), which is
stated with $\|e'\|_{L^\infty}$ and followed by an inverse estimate,
and from the consistency of $\Lhp$ with $\Lp$
(\cref{lem:bem-invcons}). In the remainder, the first variation of
$\tau^{-1}\bXh$ and the multiplier increment cost $\tau^{-1}$ and no negative power of $h$, because the increment bound of \cref{thm:in-c1p} is measured by $\|e'\|_{L^2}$, while the
stiffness increments and the second variation of the mass cost
$h^{-1/2}$, because $\|e_h'\|_{L^\infty}$ appears once in each, as
in \eqref{eq:in-C3}. \Cref{sec:proof} combines the three rates with \cref{prop:stab} and shows which of them impose the coupling \eqref{eq:coupling}.

\subsection{Proof of \texorpdfstring{\cref{thm:L2}}{the main theorem}}\label{sec:proof}

The induction, which follows \cite[\S6.2, (6.5)--(6.8)]{Rao2026}, proves the following bound of the error against the
interpolated semidiscrete curve.

\begin{lemma}[Error against the interpolated semidiscrete curve]\label{prop:superclose}
Under the hypotheses of \cref{thm:L2}, there are $\tau_0>0$, $c_0>0$ and $C$ such that, for $0<\tau\le\tau_0$, \eqref{eq:coupling} and
$0\le m\le T/\tau$, the polygons satisfy
\begin{equation}\label{eq:superclose}
\|x_h^m-\Ih x^m\|_{H^1}\le C\,\frac{h^2}{\tau} .
\end{equation}
\end{lemma}

We prove \cref{prop:superclose} and \cref{thm:L2}(a) together, with the same $\tau_0$ and $c_0$, starting
from the decomposition \eqref{eq:err-split}, in which $e^0=0$ and
$\mathcal R^m$ is the remainder of \cref{lem:fd-rem}. Combining
\cref{prop:stab} with \cref{lem:fd-comp} as in \cite[(5.28)]{Rao2026}, and
including the interpolation error \eqref{eq:fd-interp} of $\Ih D\Ft(x^m)e$, the triangle
inequality and the equivalence of the adapted norm with the $H^1$ norm (\cref{lem:p52-adnorm}), we find that the linear term
satisfies
\begin{equation}\label{eq:fd-stab}
\xnorm[x^{m+1}]{D\Fht(\Ih x^m)e}\le\big(1+C\tau+Ch\tau^{-1}+Ch^{1/2}\big)\xnorm[x^m]{e}
\end{equation}
for piecewise linear $e$ and $0<h\le\tau\le1$. With \cref{lem:fd-cons,lem:fd-rem} and the
equivalence of $\xnorm\cdot$ with the $H^1$ norm, we obtain, as long as the
segments lie in $\mathcal N_{r_2}$,
\begin{equation}\label{eq:recursion}
\xnorm[x^{m+1}]{e^{m+1}}\le\big(1+C\tau+C(h\tau^{-1}+h^{1/2})\big)\xnorm[x^m]{e^m}+Ch^2+C\big(h^{-1/2}+\tau^{-1}\big)\xnorm[x^m]{e^m}^2 .
\end{equation}
We show by induction that $\xnorm[x^m]{e^m}\le\mathfrak Kh^2/\tau$ for
$0\le m\le T/\tau$, with $\mathfrak K$ chosen below. Under this hypothesis, the
perturbations of the linear factor in \eqref{eq:recursion} are bounded by
$C'\tau$ provided
the conditions
\begin{equation}\label{eq:conditions}
\underbrace{h\tau^{-1}+h^{1/2}\le\tau}_{\text{comparison}},\qquad
\underbrace{\mathfrak Kh^{3/2}\le\tau^2}_{\text{remainder in }h^{-1/2}},\qquad
\underbrace{\mathfrak Kh^2\le\tau^3}_{\text{remainder in }\tau^{-1}}
\end{equation}
hold. The comparison condition comes from \cref{lem:fd-comp}, and each of its two terms requires $h\lesssim\tau^2$. The sources of the two losses are described after \cref{lem:fd-rem}, and the loss $h\tau^{-1}$ also bounds $\tau$ from below. The two remainder conditions come from \cref{lem:fd-rem} and require only
$h\lesssim\tau^{4/3}$ and $h\lesssim\tau^{3/2}$. Therefore, the coupling
\eqref{eq:coupling}, which is the coupling $h\le c\tau^2$ of \cite{Rao2026},
is imposed by the comparison of the linearized steps alone. We do
not claim that it is sharp. A formal expansion as in \cite[\S4]{LiY2026} indicates, however, that the present argument cannot reach the regime $\tau\sim h^2$. The tangential row of \eqref{eq:step} is that of the classical BGN scheme, whose tangential velocity depends at leading order on $\tau/h^2$ when $\tau\sim h^2$ \cite[Remark~4.10]{LiY2026}, so the polygons do not stay $H^1$-close to the interpolated semidiscrete curves, which correspond to the limit $\tau/h^2\to\infty$. This regime requires a comparison with a modified flow whose tangential velocity depends on $\tau/h^2$, as in \cite{LiY2026}.

The conditions \eqref{eq:conditions} hold under \eqref{eq:coupling} for
$\tau\le\tau_0$, with $c_0$ depending on $\mathfrak K$ and $\tau_0$ depending on
$c_0$, and so do $h\le\tau$ and the smallness of $h$ required by
\cref{lem:fd-cons,lem:fd-comp,lem:fd-rem}. Indeed, let $c_0\le\tfrac14$,
$\mathfrak Kc_0\le1$ and $\tau\le1$. Then, $h\le c_0\tau^2\le1$, and we have
\[
h\tau^{-1}+h^{1/2}\le(c_0+c_0^{1/2})\tau\le\tau,\qquad
\mathfrak Kh^{3/2}\le\mathfrak Kc_0\tau^2h^{1/2}\le\tau^2,\qquad
\mathfrak Kh^2\le\mathfrak Kc_0^2\tau^4\le\tau^3 .
\]
The constants are chosen in the
order $\mathfrak K$, then $c_0$, then $\tau_0$. The threshold $\tau_0$ is also
taken not larger than $1$, than the threshold $\tau_B$ of \cref{thm:semi-B}, than
the threshold $\tau_P$ of \cref{prop:stab}, than the threshold $\tau_V$ of
\cref{lem:Kfamily}, than the thresholds of \cref{lem:fd-cons,lem:fd-comp,lem:fd-rem,lem:p52-adnorm}, and than the thresholds of \cref{lem:semi-material} for
the finitely many orders used. It is also taken so small that the bound of $\|e^m\|_{W^{1,\infty}}$ below is not larger than $r_2/2$, with $r_2$ the radius of \cref{lem:fd-rem}. Then, \eqref{eq:recursion} gives
$\xnorm{e^{m+1}}\le(1+C''\tau)\xnorm{e^m}+Ch^2$. With $\bar C=\max(C,C'')$, the discrete Gronwall
lemma, \cref{lem:tool-gronwall}(a) applied to $\tau\xnorm{e^m}/h^2$, gives $\xnorm{e^m}\le \bar C\exp(\bar CT)\,Th^2/\tau$, since $e^0=0$. We choose
$\mathfrak K=2\bar C\exp(\bar CT)\,T$. The segments stay in $\mathcal N_{r_2}$ because the inverse
estimate, the norm equivalence and the remainder condition in $h^{-1/2}$ of \eqref{eq:conditions} give
\[
\|e^m\|_{W^{1,\infty}}\le Ch^{-1/2}\|e^m\|_{H^1}\le Ch^{-1/2}\,\frac{\mathfrak Kh^2}\tau\le C\tau\le\frac{r_2}2 .
\]
The existence and unique solvability of all steps, part of \cref{thm:L2}(a), are proved inside the same induction, since an iterate $x_h^m$ that satisfies the induction hypothesis lies in $\mathcal N_{r_2}\subset\Nd$, where \cref{cor:bem-acoerc} applies. The
polygons are regular, embedded and quasi-uniform because they are
$W^{1,\infty}$-close to the smooth embedded curves $x^m$. This proves
\cref{thm:L2}(a) and \cref{prop:superclose}.

For \eqref{eq:main}, we have
$\|x_h^m-x^m\|_{W^{1,\infty}}\le\|e^m\|_{W^{1,\infty}}+\|\Ih x^m-x^m\|_{W^{1,\infty}}\le C(h^{3/2}\tau^{-1}+h)$,
which is bounded by $C\tau$ under \eqref{eq:coupling}. Reparametrization by
normalized arclength is Lipschitz in $W^{1,\infty}$ near the smooth curve
$x=x^m$, since $|x'|$ is bounded below (\cref{lem:semi-material}). Indeed, the normalized arclength maps $\phi$ of $x$ and $\phi_h$ of $x_h=x_h^m$ (\cref{sec:ms-flow}) satisfy $\|\phi_h-\phi\|_{W^{1,\infty}}\le C\delta$ with $\delta=\|x_h-x\|_{W^{1,\infty}}$, so do their inverses $\chi$ and $\chi_h$, and splitting $x_h'(\chi_h)\chi_h'-x'(\chi)\chi'$ at $x'(\chi_h)$ gives $\|\hat x_h-\hat x\|_{W^{1,\infty}}\le C\delta$.
Hence, $\|\hat x_h^m-\hat x^m\|_{W^{1,\infty}}\le C\tau$,
and \cref{thm:semi-B} gives $\|\hat x^m-\hat X(t_m)\|_{C^1}\le C\tau$. The
Hausdorff distance is bounded by the $L^\infty$ distance of the
parametrizations. \qed

\section{Formal verification and proof audit}\label{sec:formal}

\subsection{Scope of the formalization}\label{sec:formal-scope}
The formalization consists of two projects in Lean~4 \cite{MU21}, version 4.28.0-rc1, with Mathlib \cite{MC20} at revision 789108e1. The first is the formalization of \cite{EtoPaperI}, which is archived at \url{https://doi.org/10.5281/zenodo.24103886}. The second depends on the first, contains the proofs of the present paper, and is archived at \url{https://doi.org/10.5281/zenodo.24305694}. The results stated as top-level theorems of the formalization are \cref{thm:L2}, \cref{thm:semi-B}, \cref{lem:semi-material} and \cref{prop:superclose}, and the other results of the paper, except the external results of \cref{sec:formal-ext}, appear there or in the formalization of \cite{EtoPaperI} as intermediate theorems, as combinations of such theorems or, for some identities of \cref{lem:fd-cancel}, as steps of the formal proofs (\cref{app:map}). Some formal statements carry hypotheses that the printed statements leave implicit, such as the smallness of $\tau$, the embeddedness of a curve or the membership in the admissible class, and these hypotheses hold wherever the results are used.

\subsection{Formal statement and correspondence}\label{sec:formal-stmt}
The formal counterpart of \cref{thm:L2} is the theorem in \cref{lst:thm11}, in which we quote the statement and its conclusion verbatim, with the Unicode symbols of Lean typeset as mathematical symbols.
\begin{listing}[H]
\begin{lrbox}{\leanbox}\begin{minipage}{\dimexpr\linewidth-2\fboxsep-2\fboxrule\relax}\small
\begin{alltt}
theorem thm_1_1
    (hD : Ext.DavidThm6_L2) (hC : Calderon.CalderonParamStmt)
    (hE : Stage2.Cons.EDtNStmt) (hPa : Stage2.EPotCoerStmt)
    (hDS : Stage3.S5.EDirSolvCStmt) : Thm11Stmt

def Thm11Stmt : Prop :=
  \(\forall\) (T \(\eta\)\(\sb0\) : \(\mathbb{R}\)) (Xh : \(\mathbb{R}\to\mathbb{R}\to\mathbb{C}\)), 0 < T \(\to\) 0 < \(\eta\)\(\sb0\) \(\to\) HypE T \(\eta\)\(\sb0\) Xh \(\to\)
    \(\exists\) \(\tau\)\(\sb0\) c\(\sb0\) C Cq : \(\mathbb{R}\), 0 < \(\tau\)\(\sb0\) \(\wedge\) 0 < c\(\sb0\) \(\wedge\) 0 < C \(\wedge\) 0 < Cq \(\wedge\)
      \(\forall\) \(\tau\) : \(\mathbb{R}\), 0 < \(\tau\) \(\to\) \(\tau\) \(\le\) \(\tau\)\(\sb0\) \(\to\) \(\forall\) (N : \(\mathbb{N}\)) [NeZero N], (1 : \(\mathbb{R}\)) / N \(\le\) c\(\sb0\) * \(\tau\) ^ 2 \(\to\)
        \(\exists\) Xs : \(\mathbb{N}\to\) Fin N \(\to\mathbb{C}\), Thm11Concl T \(\tau\) Xh C Cq Xs

def Thm11Concl (T \(\tau\) : \(\mathbb{R}\)) (Xh : \(\mathbb{R}\to\mathbb{R}\to\mathbb{C}\)) (C Cq : \(\mathbb{R}\)) \{N : \(\mathbb{N}\)\} [NeZero N]
    (Xs : \(\mathbb{N}\to\) Fin N \(\to\mathbb{C}\)) : Prop :=
  IsBemIterates T \(\tau\) Xh Xs \(\wedge\)
    \(\forall\) m \(\le\) \(\lfloor\)T / \(\tau\)\(\rfloor\sb{+}\),
      (PolyRegular (Xs m) \(\wedge\) PolyEmbedded (Xs m) \(\wedge\) PolyQU Cq (Xs m)) \(\wedge\)
      (IsUnit (Vmat (Xs m)).det \(\wedge\) IsUnit (Vbord (Xs m)).det) \(\wedge\)
      W1inf (xhat (p1 (Xs m)) - slice Xh (m * \(\tau\))) \(\le\) ENNReal.ofReal (C * \(\tau\)) \(\wedge\)
      Metric.hausdorffDist (Set.range (p1 (Xs m))) (Set.range (slice Xh (m * \(\tau\))))
        \(\le\) C * \(\tau\)
\end{alltt}\end{minipage}\end{lrbox}
\noindent\fcolorbox{black!35}{black!4}{\usebox{\leanbox}}
\caption{The formal statement of \cref{thm:L2} and its conclusion.}\label{lst:thm11}
\end{listing}
\noindent Here, \lean{HypE} states the hypotheses of \cref{thm:L2} on the exact solution, including \ref{item:E-cap}, \lean{IsBemIterates} states that the polygons \texttt{Xs m} are the iterates of the scheme \eqref{eq:step} started from $\Ih\hat X(\cdot,0)$ and that every step is uniquely solvable, and \lean{Thm11Concl} states the geometric properties and the error bounds of \cref{thm:L2}. The invertibility of the matrices \lean{Vmat} and \lean{Vbord} ensures that $\Lh$, which the formalization defines with a matrix inverse, is the matrix of the scheme. The mathematical statements and their Lean counterparts correspond as follows.
\begin{center}\small
\begin{tabular}{@{}ll@{}}
\toprule
Mathematical statement & Lean declaration\\
\midrule
\cref{thm:L2} & \lean{thm_1_1}\\
Hypotheses \ref{item:E-smooth}--\ref{item:E-cap} on the exact solution & \lean{HypE}\\
Iteration of the scheme \eqref{eq:step} & \lean{IsBemIterates}\\
Conclusion of \cref{thm:L2} & \lean{Thm11Concl}\\
$W^{1,\infty}$ error in \eqref{eq:main} & \lean{W1inf}\\
\cref{thm:semi-B} & \lean{thm_semi_B}\\
\cref{lem:semi-material} & \lean{lem_3_3}\\
\cref{prop:superclose} & \lean{thm_1_1_superclose_h2}\\
\bottomrule
\end{tabular}
\end{center}
The kernel check alone does not guarantee that the Lean definitions and statements agree with \cref{thm:L2} and with the sources of the five hypotheses. This agreement is a matter of reading, and a finer table, which matches each definition and statement of the paper with its Lean declaration, will accompany the code.

\subsection{External assumptions and kernel audit}\label{sec:formal-ext}
The five hypotheses of \lean{thm_1_1} are the external inputs of \cref{sec:tools-ext}, in the parameter forms stated there:
\begin{center}\small
\begin{tabular}{@{}llll@{}}
\toprule
\multicolumn{2}{@{}l}{Hypothesis of \lean{thm_1_1}} & Statement & Source\\
\midrule
\lean{hD} & \lean{Ext.DavidThm6_L2} & \cref{prop:ext-david} & \cite[Thm.~6, (66)]{David1984}\\
\lean{hC} & \lean{Calderon.CalderonParamStmt} & \cref{prop:ext-calderon} & \cite[Cor.~6.19 and Thm.~6.15]{Steinbach2008}\\
\lean{hDS} & \lean{Stage3.S5.EDirSolvCStmt} & \cref{prop:ext-dir} & \cite[Thm.~6.14]{GilbargTrudinger2001}, Kelvin transform\\
\lean{hE} & \lean{Stage2.Cons.EDtNStmt} & \cref{prop:ext-dtn} & \cite[(6.1), (6.22) and \S7.5]{Steinbach2008}\\
\lean{hPa} & \lean{Stage2.EPotCoerStmt} & \cref{prop:ext-slcoer} & \cite[Thm.~6.22]{Steinbach2008}\\
\bottomrule
\end{tabular}
\end{center}
The command \texttt{\#print axioms thm\_1\_1} reports the axioms in \cref{lst:axioms}, which are the standard axioms of Lean.
\begin{listing}[H]
\begin{lrbox}{\leanbox}\begin{minipage}{\dimexpr\linewidth-2\fboxsep-2\fboxrule\relax}\small
\begin{alltt}
\textquotesingle{}SGBEMConvergence.thm_1_1\textquotesingle{} depends on axioms: [propext, Classical.choice, Quot.sound]
\end{alltt}\end{minipage}\end{lrbox}
\noindent\fcolorbox{black!35}{black!4}{\usebox{\leanbox}}
\caption{The axioms on which the formal proof of \cref{thm:L2} depends.}\label{lst:axioms}
\end{listing}
\noindent The five hypotheses are arguments of the theorem and not axioms, so this output does not mean that they are proved. The maximal inequality of \cref{prop:ext-hl} and the Fredholm alternative of \cref{prop:ext-slfred} are proved in the formalization (\lean{hlMaximal_R_proof} and \lean{epotFred_proof}), \cref{prop:ext-layer} is derived there from \cref{prop:ext-david,prop:ext-calderon}, and \cref{prop:in-import} follows from the theorems of the formalization of \cite{EtoPaperI}, whose hypotheses are discharged by \lean{hD}, \lean{hC} and the formal proof of \cref{prop:ext-hl}, and the same holds for the consistency of the shape derivative of $\Lh$ in \cref{sec:input-results}, which is \cite[Theorem~3.7]{EtoPaperI}. The formal proofs of \cref{prop:ext-hl,prop:ext-slfred} are self-contained, and the classical sources cited after these propositions record their origin and are not used as hypotheses. Similarly, the estimates of Rao \cite{Rao2026} and of \cite[Thm.~8.3 and Thm.~10.9]{Steinbach2008} that the printed proofs quote are not additional hypotheses, since the formal proof of \cref{thm:L2} uses no hypotheses besides the five above. The counterparts of \cref{thm:semi-B,lem:semi-material} take the four hypotheses other than \lean{hC}. The theorem \lean{Audit.hypE_satisfiable} shows that a static circle satisfies \lean{HypE}, so that the hypothesis on the exact solution can be met.

\section*{Use of Generative-AI Tools Declaration}
The author used Claude Code (Anthropic), Codex (OpenAI), and ChatGPT (OpenAI) in the preparation of this manuscript.
These tools assisted in drafting and editing the mathematical exposition and proofs, as well as in developing the Lean 4 formalization mentioned above.
AI assistance was used throughout the manuscript, including its mathematical arguments and presentation.
The author reviewed the resulting mathematical statements, proofs, references, and formalization, and assumes full responsibility for the accuracy and integrity of the manuscript.

\section*{Acknowledgments}
This research is supported by the ANR (Agence Nationale de la Recherche), grant ANR-23-PEIA-0004 (PEPR IA, PDE-AI).

\appendix

\section{Conventions and external inputs}\label{app:tools}

This appendix fixes the notation on the torus that the main text and \cref{app:new} use, and it states the five results from the literature on which the main theorem depends as hypotheses, each with its source and the argument by which the source applies.

\subsection{Sobolev spaces on the torus}\label{sec:tools-fourier}

A function on $\Tor=\R/\Z$ is expanded as $g=\sum_n\hat g_n\mathbf e_n$ with $\mathbf e_n(\xi)=\eu^{2\pi\iu n\xi}$. With $\langle n\rangle=(1+|n|^2)^{1/2}$, we set $\|g\|_{H^r}^2=\sum_n\langle n\rangle^{2r}|\hat g_n|^2$ for $r\in\R$, and $\wnorm{g}=\sum_n|\hat g_n|$ is the Wiener norm. On $\Tor^2$, we use the same definitions with $\langle(p,q)\rangle=(1+p^2+q^2)^{1/2}$. For an integer $k\ge0$, we set $\|f\|_{W^{k,\infty}}=\sum_{j\le k}\|f^{(j)}\|_{L^\infty}$. Bounded operators from $H^a$ to $H^b$ form the space $\mathcal B(H^a,H^b)$.

For a kernel $k$ on $\Tor^2$, we write $\mathcal I_k\varphi(\xi)=\int_\Tor k(\xi,\eta)\varphi(\eta)\,\dd\eta$. For a closed curve $Z:\Tor\to\C$, the difference quotient is
\[
D[Z](\xi,\eta)=\frac{Z(\xi)-Z(\eta)}{\mathbf e_1(\xi)-\mathbf e_1(\eta)}\qquad(\xi\ne\eta),
\]
an $L^2$ function on $\Tor^2$ if $Z\in H^{1/2}$, and a continuous one with $D[Z](\xi,\xi)=Z'(\xi)/(2\pi\iu\,\mathbf e_1(\xi))$ if $Z$ is $C^1$.

\subsection{External inputs}\label{sec:tools-ext}

We state each of the five propositions below in the form in which the proofs use it. Throughout this subsection, a smooth curve is a smooth closed curve $x:\Tor\to\C$ that is regular and embedded, so that $\Gamma=x(\Tor)$ is a smooth Jordan curve. When $x$ is counter-clockwise, $\Omega_{\rm int}$ and $\Omega_{\rm ext}$ denote the bounded and the unbounded component of $\C\setminus\Gamma$, and $\bn(\xi)=\Jm x'(\xi)/|x'(\xi)|$ is the outward unit normal of $\Omega_{\rm int}$ at $x(\xi)$. A function $g$ on $\Tor$ is identified with the function $g\circ x^{-1}$ on $\Gamma$. For a smooth $g:\Tor\to\R$, a pair $(u_{\rm int},u_{\rm ext})$ of Dirichlet extensions of $g$ consists of a harmonic function $u_{\rm int}$ in $\Omega_{\rm int}$ and a bounded harmonic function $u_{\rm ext}$ in $\Omega_{\rm ext}$, both continuous up to $\Gamma$ and equal to $g$ on $\Gamma$. For such a pair, we consider the limit
\begin{equation}\label{eq:ext-lambda}
\lambda(\xi)=\lim_{\varepsilon\downarrow0}\frac{2g(\xi)-u_{\rm int}(x(\xi)-\varepsilon\bn(\xi))-u_{\rm ext}(x(\xi)+\varepsilon\bn(\xi))}{\varepsilon} .
\end{equation}

The main theorem depends on the five propositions below as hypotheses (\cref{sec:formal-ext}). \Cref{prop:ext-david} is an external input of \cite{EtoPaperI} and is used through the results of \cite{EtoPaperI} listed in \cref{sec:bem}. \Cref{prop:ext-calderon} is used through the layer representation of the Dirichlet-to-Neumann map, \cref{prop:ext-layer}. The last three propositions are used directly in the present proofs.

\begin{proposition}[David's theorem {\cite[Proposition~A.2]{EtoPaperI}}]\label{prop:ext-david}
There are constants $C$ and $\beta$ with the following property. Let $\Gamma$ be a rectifiable curve with arclength measure $\mu$, and let $\bar\mu(z)$ be the supremum of $\mu(D)/r(D)$ over the discs $D$ containing $z$, where $r(D)$ is the radius of $D$. Then, the maximal chordal truncation
\[
T^*_\mu g(z)=\sup_{\varepsilon>0}\Big|\int_{|z-w|\ge\varepsilon}\frac{g(w)}{z-w}\,\dd\mu(w)\Big|
\]
satisfies $\|T^*_\mu g\|_{L^2(\mu)}\le C\|(1+\bar\mu)^\beta g\|_{L^2(\mu)}$ for every continuous $g$ with compact support.
\end{proposition}

This is \cite[Thm.~6, (66)]{David1984} for the kernel $1/z$ and the exponent $2$, with the function $\bar\mu$ of \cite[(10)]{David1984}.

\begin{proposition}[Calder\'on identities]\label{prop:ext-calderon}
Let $x$ be a curve of a family that satisfies \ref{K1}--\ref{K3}, and let $\mathfrak C$ and $r_1$ be as in \cref{sec:bem}. Let $Y$ be a polygon of $\mathfrak C$, or the curve $x$, or a curve $x+\varepsilon e$ with a P1 direction $e$ and $|\varepsilon|\,\|e'\|_{L^\infty}\le2r_1$.
\begin{enumerate}
\item For $\varsigma,\vartheta\in L^2$, we have $\langle\mathcal K_Y\mathcal L_Y\varsigma,\vartheta\rangle=\langle\mathcal K_Y\mathcal L_Y\vartheta,\varsigma\rangle$.
\item Let $u$ be a periodic $C^1$ function or a P1 function, let $\varsigma\in L^2$, and let $v$ be a Lipschitz function with $\mathcal L_Y\varsigma=v$ almost everywhere. Then, $\langle\mathcal L_Y\partial_\xi u,\partial_\xi v\rangle=\langle\tfrac14u-\mathcal K_Y^2u,\varsigma\rangle$.
\item We have $\mathcal K_x1=-\tfrac12$ almost everywhere.
\end{enumerate}
\end{proposition}

These are the identities (6.26), $VK'=KV$, and (6.24), $VD=(\frac12I+K)(\frac12I-K)$, of \cite[Cor.~6.19, p.~138]{Steinbach2008}, together with Maue's formula \cite[Thm.~6.15, p.~131]{Steinbach2008} and the relation $(\frac12I+K)1=0$, which is \cite[(6.22)]{Steinbach2008} for $u\equiv1$, since the coefficient of the identity in \cite[(6.22) and (6.24)]{Steinbach2008} equals $\frac12$ almost everywhere \cite[(6.11), p.~123]{Steinbach2008}. In the parameter, $(V\psi)\circ Y=\mathcal L_Y(\psi\circ Y\,|Y'|)$, $(Ka)\circ Y=\mathcal K_Y(a\circ Y)$ and $\dd s=|Y'|\,\dd\xi$, so (a) is $(KV\psi,\phi)_{L^2(\Gamma)}=(\psi,KV\phi)_{L^2(\Gamma)}$, and (b) is $(Du,V\psi)_{L^2(\Gamma)}=(VDu,\psi)_{L^2(\Gamma)}=((\frac14-K^2)u,\psi)_{L^2(\Gamma)}$, where Maue's formula gives the left-hand side. The normal of \cite{Steinbach2008} is the outward one for the counter-clockwise curve $x$, which gives the sign in (c). Maue's formula is stated for functions $v$ that are continuous and continuously differentiable on the parts of $\Gamma$, and (b) extends it to Lipschitz $v$, since the piecewise linear interpolants of $v$ converge in $H^1$ and both sides are continuous in $v\in H^1$. Chapter~6 of \cite{Steinbach2008} assumes a Lipschitz boundary, and every curve $Y$ of the proposition bounds a Lipschitz domain \cite[Proposition~A.9]{EtoPaperI}(i).

\begin{proposition}[The two-phase Dirichlet problem]\label{prop:ext-dir}
Let $x$ be a smooth counter-clockwise curve, and let $g:\Tor\to\R$ be smooth. There is a pair of Dirichlet extensions of $g$ for which the limit \eqref{eq:ext-lambda} exists for every $\xi$ and is a continuous function of $\xi$.
\end{proposition}

The interior problem is solved by the Schauder theory. The domain $\Omega_{\rm int}$ is a $C^{2,\alpha}$ domain \cite[p.~94]{GilbargTrudinger2001}, the boundary values extend to $C^{2,\alpha}(\overline{\Omega}_{\rm int})$ \cite[Lemma~6.38]{GilbargTrudinger2001}, and \cite[Thm.~6.14]{GilbargTrudinger2001} gives a harmonic $u_{\rm int}\in C^{2,\alpha}(\overline{\Omega}_{\rm int})$ equal to $g$ on $\Gamma$. The exterior problem is reduced to an interior one by an inversion. We fix $a\in\Omega_{\rm int}$ and let $\iota(z)=a+(z-a)/|z-a|^2$. Composition with $\iota$ is the two-dimensional Kelvin transform and preserves harmonicity \cite[(4.31)--(4.32)]{GilbargTrudinger2001}, and $D=\iota(\Omega_{\rm ext})\cup\{a\}$ is a bounded domain with the smooth boundary $\iota(\Gamma)$. The interior result gives a harmonic $w\in C^{2,\alpha}(\overline D)$ with the boundary values $g\circ x^{-1}\circ\iota$, and $u_{\rm ext}=w\circ\iota$ is bounded, harmonic in $\Omega_{\rm ext}$ and of class $C^{2,\alpha}$ up to $\Gamma$. Since both functions are $C^1$ up to $\Gamma$, the limit \eqref{eq:ext-lambda} exists and equals $\partial_{\bn}u_{\rm int}-\partial_{\bn}u_{\rm ext}=\Lam_\Gamma g$, which is continuous.

Every pair of Dirichlet extensions of $g$ coincides with this one. For $u_{\rm int}$, this is the maximum principle. For $u_{\rm ext}$, the difference with $w\circ\iota$, composed with $\iota$, is a bounded harmonic function in $D\setminus\{a\}$ that vanishes on $\partial D$, so the singularity at $a$ is removable and the difference vanishes. Since $w$ is $C^1$ near $a$ and $|D\iota(y)|=|y-a|^{-2}$, we also have $u_{\rm ext}(y)-u_\infty=O(|y|^{-1})$ with $u_\infty=w(a)$, and $|\nabla u_{\rm ext}(y)|=O(|y|^{-2})$, as $|y|\to\infty$.

\begin{proposition}[Single layer equation for the Dirichlet-to-Neumann density]\label{prop:ext-dtn}
Let $x$ be a smooth counter-clockwise curve, let $g:\Tor\to\R$ be smooth, and let $(u_{\rm int},u_{\rm ext})$ be a pair of Dirichlet extensions of $g$ for which the limit \eqref{eq:ext-lambda} exists for every $\xi$ and is a continuous function of $\xi$. Then, $\int_\Gamma\lambda\,\dd s=0$, and there is $c\in\R$ with $\SL\lambda=g+c$ almost everywhere on $\Gamma$, that is, $\mathcal L_x(\lambda|x'|)=g+c$ almost everywhere on $\Tor$.
\end{proposition}

The single layer $\SL$ and the normal $\bn$ of \cref{sec:setting} are those of \cite{Steinbach2008}, with the fundamental solution $-\frac1{2\pi}\log|y-y'|$ \cite[p.~111]{Steinbach2008}. By the uniqueness after \cref{prop:ext-dir}, $\lambda=\gamma_1^{\rm int}u_{\rm int}-\gamma_1^{\rm ext}u_{\rm ext}$, and $u_{\rm ext}$ has the decay stated there. The interior representation formula \cite[(6.1), (6.22)]{Steinbach2008}, in which the coefficient of the identity equals $\frac12$ on a smooth boundary \cite[p.~123]{Steinbach2008}, gives $g=\SL\gamma_1^{\rm int}u_{\rm int}+\frac12g-\mathcal Kg$. The exterior representation formula \cite[\S7.5, pp.~181--182]{Steinbach2008}, whose limit over large circles is justified by the decay, gives $g=u_\infty-\SL\gamma_1^{\rm ext}u_{\rm ext}+\frac12g+\mathcal Kg$ with the same normal $\bn$. Adding the two identities, which hold in $H^{1/2}(\Gamma)$, gives $\SL\lambda=g-u_\infty$ almost everywhere. The interior normal derivative has mean zero by Gauss's theorem, and Green's formula in $\Omega_{\rm ext}\cap\{|y|<R\}$ and $|\nabla u_{\rm ext}|=O(R^{-2})$ give the same for the exterior one.

Since $\lambda$ is continuous, $\SL\lambda$ is continuous on $\Gamma$, so the identity of \cref{prop:ext-dtn} holds at every point of $\Gamma$, and $\mathcal L_x(\lambda|x'|)=g+c$ holds at every $\xi$. This is the form in which the proposition is used, always for the pair of \cref{prop:ext-dir}, whose $\lambda$ is continuous.

The last proposition is the coercivity of the single layer. A curve $Y$ parametrized proportionally to arclength is one with $|Y'|\equiv\ell$, and, for it, the parameter density of an arclength density $w$ is $\varsigma=\ell\,w\circ Y$.
\begin{proposition}[Coercivity of the single layer]\label{prop:ext-slcoer}
Let $Y$ be a smooth curve parametrized proportionally to arclength. There is $c>0$ such that $\langle\mathcal L_Y\varsigma,\varsigma\rangle\ge c\|\varsigma\|_{H^{-1/2}}^2$ for every $\varsigma\in L^2(\Tor)$ with $\int_\Tor\varsigma\,\dd\xi=0$.
\end{proposition}

This is \cite[Thm.~6.22]{Steinbach2008} for $d=2$, on the space $H^{-1/2}_*(\Gamma)$ of densities $w$ with $(w,1)_{L^2(\Gamma)}=0$ \cite[(6.30)]{Steinbach2008}. In two dimensions, the theorem puts no condition on the capacity of $\Gamma$. With $\varsigma=\ell\,w\circ Y$, we have $\int\varsigma\,\dd\xi=(w,1)_{L^2(\Gamma)}$ and $\langle\mathcal L_Y\varsigma,\varsigma\rangle=(\SL w,w)_{L^2(\Gamma)}$, and the norms of $H^{-1/2}(\Gamma)$ and of $H^{-1/2}(\Tor)$ in the parameter are equivalent, with constants that depend on the curve.

\section{The new arguments}\label{app:new}

This appendix proves the arguments described in the strategy in \cref{sec:intro} that have no counterpart for curve shortening, and, in \cref{app:bem-step}, the one-step results of \cref{sec:bem}.

\subsection{The splitting of the Dirichlet-to-Neumann map}\label{sec:app-cont}

\begin{definition}\label{def:cont-class}
Let $\sigma>\tfrac12$, $M>0$, $c_{\rm ch}>0$ and $0<\ell_{\min}\le\ell_{\max}$.
The class $\mathcal C=\mathcal C(\sigma;M,c_{\rm ch},\ell_{\min},\ell_{\max})$
consists of the closed $C^1$ curves $x:\Tor\to\R^2\cong\C$ for which there are
a length $\ell\in[\ell_{\min},\ell_{\max}]$, a real function $\kc$ with
$\|\kc\|_{H^\sigma}\le M$ and a $C^1$ function $y:\Tor\to\C$ with $|y|=1$ such that
\[
x'=\ell y,\qquad y'=\iu\ell\kc y,\qquad |x(\xi)-x(\eta)|\ge c_{\rm ch}\,\ell\,d(\xi,\eta),
\]
where $d$ is the distance on $\Tor$.
\end{definition}

As for the curves of \cref{sec:ms-flow}, $\ell$ is the length, $y$ the unit tangent and $\kc$ the curvature in normalized arclength, and they are determined by $x$. When $M$, $c_{\rm ch}$, $\ell_{\min}$ and $\ell_{\max}$ are fixed, we write $\mathcal C(\sigma)$ for $\mathcal C(\sigma;M,c_{\rm ch},\ell_{\min},\ell_{\max})$. Constants $C$ depend only on the constants of the class, on the constant $\mu$ of \cref{thm:cont-N1} below, and on the Sobolev indices.

Since $|\mathbf e_1(\eta)-\mathbf e_1(\xi)|=2|\sin\pi(\xi-\eta)|$, we have
$\log|x(\xi)-x(\eta)|=\log|2\sin\pi(\xi-\eta)|+\log|D[x](\xi,\eta)|$ off the
diagonal. We set
\[
G(\xi,\eta)=\log|D[x](\xi,\eta)|-\log\frac{\ell}{2\pi},\qquad
\SL_G=-\frac{\ell}{2\pi}\,\mathcal I_G .
\]
The function $D[x]$ is continuous on $\Tor^2$. By \cref{lem:cont-kernelreg}, $G$ lies in $H^{\sigma+3/2}(\Tor^2)$ for $\sigma>\frac52$, and $\mathcal I_G$ maps $H^{-a}$ into $H^b$ for $a_++b_+\le\sigma+\frac32$. The index $\sigma+\tfrac32$ comes from the bound $\|D[Z]\|_{H^a(\Tor^2)}\le2^a\|Z\|_{H^{a+1/2}}$ of \cref{lem:tool-kernel}, in which each coefficient $\hat Z_n$ appears in $|n|$ coefficients of $D[Z]$. It gives $D[x]\in H^{\sigma+3/2}(\Tor^2)$ for a curve whose curvature lies in $H^\sigma$, so that only half a derivative of the curve is lost.

For a curve $x$ of the class, the single layer acting on densities with
respect to arclength is
$\SL\varphi(\xi)=-\frac{\ell}{2\pi}\int_\Tor\log|x(\xi)-x(\eta)|\,\varphi(\eta)\,\dd\eta$,
so that $\SL\varphi=\mathcal L_x(\ell\varphi)$ with the parameter form
$\mathcal L_x$ of \cref{sec:input-K}. Since the parametrization has constant
speed, the arclength mean is the mean in $\xi$. Let $\Pz$ remove the mean, let
$\Lp=\Pz\SL\Pz$, and let $\Lam_0=2\ell^{-1}|\partial_\xi|$ be the Fourier
multiplier $4\pi|n|/\ell$. By $\log|2\sin\pi\upsilon|=-\sum_{n\ge1}n^{-1}\cos2\pi
n\upsilon$, the logarithmic part of $\SL$ is the multiplier $\SL_0=\ell/(4\pi|n|)$
on functions of mean zero, which is the inverse of $\Lam_0$ there. The
constant $\log(\ell/2\pi)$ of the kernel acts only on the mean. Hence, we have
\begin{equation}\label{eq:cont-SLsplit}
\Pz\SL\varphi=\SL_0\varphi+\Pz\SL_G\varphi\qquad\text{if }\int_\Tor\varphi\,\dd\xi=0 .
\end{equation}

\begin{lemma}\label{thm:cont-N1}\label{lem:cont-gauge}
Let $\sigma>\tfrac52$, and let $x\in\mathcal C(\sigma)$ be smooth. Assume
that $\langle\mathcal L_x\varsigma,\varsigma\rangle\ge\mu\|\varsigma\|^2_{H^{-1/2}}$
for all $\varsigma\in L^2$, with $\mu>0$.
\begin{enumerate}[label=(\alph*),ref=\thetheorem(\alph*)]
\item For $g\in H^1$, there is a unique $\varphi\in L^2$ with
$\int_\Tor\varphi\,\dd\xi=0$ and $\Pz\SL\varphi=\Pz g$, and
$\|\varphi\|_{L^2}\le C\|g\|_{H^1}$. We write $\Lam_\Gamma g=\varphi$. If the
curve is counter-clockwise and $g$ is smooth, then $\Lam_\Gamma g$ is the
two-phase Dirichlet-to-Neumann map of $g$. Moreover, $\Lam_\Gamma\Lp\varphi=\varphi$
for every $\varphi\in L^2$ of mean zero, and $\Lp\Lam_\Gamma g=\Pz g$.
\item On $H^1$, the remainder $Q_\Gamma=\Lam_\Gamma-\Lam_0$ satisfies
\begin{equation}\label{eq:cont-gauge-id}
Q_\Gamma=-\Lam_0\SL_G\Lam_\Gamma=-\Lam_\Gamma\SL_G\Lam_0 .
\end{equation}
\item For $r\ge\tfrac12-\sigma$, the remainder $Q_\Gamma$ extends to an
operator in $\mathcal B\bigl(H^r,H^{\sigma+1/2-(1-r)_+}\bigr)$ with norm
bounded by $C$. Hence, $\Lam_\Gamma=\Lam_0+Q_\Gamma$ lies in $\mathcal
B(H^r,H^{r-1})$ for $\tfrac12-\sigma\le r\le\sigma+\tfrac32$, with norm bounded
by $C$. In particular, $Q_\Gamma$ maps $H^r$ into $H^{r+1}$ for
$\tfrac12-\sigma\le r\le\sigma-\tfrac12$.
\end{enumerate}
\end{lemma}

The hypothesis on $\mathcal L_x$ is the coercivity in \ref{K3}, on all of $L^2$, and by \cref{lem:cont-near} it holds with a uniform constant near the exact curves. Smoothness of the curve is used only through \cref{prop:ext-slfred,prop:ext-dtn,prop:ext-dir}. The two identities in (b) are used in different directions. The first one places $\SL_G$, which gains $\sigma+\frac32$ derivatives in total, after $\Lam_\Gamma$, and it gives the smoothing of $Q_\Gamma$. The second one places $\Lam_\Gamma$ last, and it gives the bound of $Q_\Gamma$ on rough data.

\begin{proof}
For part~(a), let $\varphi\in L^2$ have mean zero. Then, we have
$\langle\Pz\SL\varphi,\varphi\rangle=\ell^{-1}\langle\mathcal L_x(\ell\varphi),\ell\varphi\rangle
\ge\mu\ell\|\varphi\|^2_{H^{-1/2}}$. This gives the uniqueness and, for a solution, $\mu\ell\|\varphi\|_{H^{-1/2}}\le\|g\|_{H^{1/2}}$. For the existence, the coercivity makes $\mathcal L_x$ injective on $L^2$, so \cref{prop:ext-slfred} gives $\psi_g,\psi_1\in L^2$ with $\mathcal L_x\psi_g=g$ and $\mathcal L_x\psi_1=1$. Since $\int_\Tor\psi_1\,\dd\xi=\langle\mathcal L_x\psi_1,\psi_1\rangle>0$, the density $\varphi=\ell^{-1}\bigl(\psi_g-(\int\psi_g/\int\psi_1)\psi_1\bigr)$ is a solution. By \eqref{eq:cont-SLsplit}, the equation reads
$\SL_0\varphi=\Pz g-\Pz\SL_G\varphi$, and applying $\Lam_0$ we obtain
\begin{equation}\label{eq:cont-gauge-eq}
\varphi=\Lam_0g-\Lam_0\SL_G\varphi .
\end{equation}
\Cref{lem:cont-kernelreg}(b) with $a=\tfrac12$ and $b=1$ gives
$\|\SL_G\varphi\|_{H^1}\le C\|\varphi\|_{H^{-1/2}}$, and therefore
$\|\varphi\|_{L^2}\le C(\|g\|_{H^1}+\|\varphi\|_{H^{-1/2}})\le C\|g\|_{H^1}$.
For a smooth counter-clockwise curve and smooth $g$, the two-phase
Dirichlet-to-Neumann map $\lambda$ of $g$ has mean zero and
$\mathcal L_x(\ell\lambda)=g+c$ with a constant $c$, by
\cref{prop:ext-dir,prop:ext-dtn}, so $\lambda=\Lam_\Gamma g$ by the uniqueness. The last two identities of (a) also follow from the uniqueness, since $\Lam_\Gamma g$ depends only on $\Pz g$ and a density $\varphi$ of mean zero solves the problem for $g=\SL\varphi$.

For part~(b), the first identity of \eqref{eq:cont-gauge-id} is
\eqref{eq:cont-gauge-eq}. For the second, let $g\in H^1$ and $q=\Lam_0g$. By \eqref{eq:cont-SLsplit}, we have $\Pz\SL q=\Pz(g+\SL_Gq)$,
and $\SL_Gq\in H^1$ by \cref{lem:cont-kernelreg}(b) with $a=0$ and $b=1$. Part~(a)
gives $q=\Lam_\Gamma(g+\SL_Gq)$, which is the second identity.

For part~(c), let first $\tfrac12-\sigma\le r\le1$. In the second identity, $\Lam_0$ maps $H^r$ to $H^{r-1}$, $\SL_G$ maps $H^{r-1}$ to $H^1$ by \cref{lem:cont-kernelreg}(b) with $a=1-r$ and $b=1$, and $\Lam_\Gamma$ maps $H^1$ to $L^2$. Hence, $Q_\Gamma$ extends to $\mathcal B(H^r,L^2)$, and $\Lam_\Gamma\in\mathcal B(H^r,H^{r-1})$. In the first identity, extended to $H^r$ by continuity, $\SL_G$ then maps $H^{r-1}$ to $H^{\sigma+1/2+r}$ by \cref{lem:cont-kernelreg}(b) with $a=1-r$ and $a+b=\sigma+\tfrac32$, which gives $Q_\Gamma\in\mathcal B(H^r,H^{\sigma-1/2+r})$. For $r\ge1$, the case $r=1$ gives $Q_\Gamma\in\mathcal B(H^r,H^{\sigma+1/2})$, and $H^{\sigma+1/2}\subset H^{r-1}$ if $r\le\sigma+\tfrac32$. The last claim is the inequality $\sigma+\tfrac12-(1-r)_+\ge r+1$. For $r\le1$, it reads $\sigma\ge\tfrac32$, and for $r\ge1$, it reads $r\le\sigma-\tfrac12$.
\end{proof}

\subsection{The absorption of the derivative lost in the metric step}\label{sec:app-semi-geom}

Let $\sigma\in(5/2,3)$ be fixed. Throughout this subsection, a closed curve $x$ of length $\ell$ and curvature $\kc$ is written in its normalized arclength $\xi\in\Tor$, the hats being dropped when no confusion arises, so that
\begin{equation}\label{eq:semi-frame}
x'=\ell y,\qquad y'=\iu\ell\kc\,y,\qquad |y|=1,
\end{equation}
with $'=\dd/\dd\xi$, $\bt=y$, $\bn=-\iu y$, $\partial_s\bt=-\kc\bn$ and $\partial_s\bn=\kc\bt$. We write $\|\cdot\|_r$ for the norm of $H^r(\Tor)$ and $\wnorm{\cdot}$ for the Wiener norm of \cref{sec:tools-fourier}. By \cref{lem:tool-frame}(d), every smooth, regular and embedded closed curve of length not larger than $\ell_{\max}$ satisfies
\begin{equation}\label{eq:semi-kmin}
\|\kc\|_{L^2}\ge\kc_{\min}:=\sqrt2/\ell_{\max}.
\end{equation}

For $M$, $c_{\rm ch}$, $\ell_{\min}\le\ell_{\max}$ and $\mu>0$, let $\mathcal U_{\rm ad}(M,c_{\rm ch},\ell_{\min},\ell_{\max},\mu)$ be the set of smooth, embedded, counter-clockwise closed curves in normalized arclength that belong to the class $\mathcal C(\sigma;M,c_{\rm ch},\ell_{\min},\ell_{\max})$ of \cref{def:cont-class} and whose single layer satisfies the hypothesis of \cref{thm:cont-N1} with the constant $\mu$. The constants are fixed by \cref{lem:semi-close}. The continuum estimates hold uniformly on $\mathcal U_{\rm ad}$, and a constant $C$ depends only on $\sigma$, on the constants of $\mathcal U_{\rm ad}$, on (E) and on the orders indicated.

\subsubsection*{One step and the defect}
For a curve $x$ in normalized arclength and a smooth velocity $v=v_n\bn+v_t\bt$ on it, we set
\[
a=\bt\cdot\partial_sv=\partial_sv_t+\kc v_n,\qquad b=\bn\cdot\partial_sv=\partial_sv_n-\kc v_t,\qquad w=x+\tau v,\qquad \zeta=\big((1+\tau a)^2+\tau^2b^2\big)^{1/2},
\]
so that $|w'|=\ell\zeta$, and $\mathfrak P$ is formed with $\ell$ as in \cref{sec:semi}.

Let $\kc=\kc^m$, $\kt=\kt^m$ and $v=v^m$, written in the normalized arclength of $x^m$. The step \eqref{eq:semi} reads $v_n=-\frac12\Lam^m\kt$ and $w''=-\ell^2\kt\bn$ with $w=x^{m+1}$. By \cref{lem:semi-geom}(b) and~(c), the normalized-arclength representative of the new curve has the curvature \eqref{eq:semi-newcurv}, and the expansion \eqref{eq:semi-expand} holds at a common $\xi$. These identities, without the bound of $\rho$, use only $w''=-\ell^2\kt\bn$ and $|w'|>0$, so they hold here for every $\tau>0$, with $\rho$ defined by \eqref{eq:semi-expand}.

For every smooth function $V$ on a smooth closed curve, with $\alpha=L_\kc^{-1}D_\kc V$ and with $\mathcal T_\kc$ formed with this $V$, we have
\begin{equation}\label{eq:semi-cancel}
\alpha-\mathfrak P[\partial_s\alpha+\kc V]=\mathcal T_\kc ,
\end{equation}
because the mean of $\partial_s\alpha$ vanishes by periodicity, so that $\mathfrak P[\partial_s\alpha]=\alpha-\alpha(0)$. This extends \cite[(3.30)--(3.31)]{Rao2026}, where $V=-f$.

We apply \eqref{eq:semi-cancel} with $\kt$ in place of $\kc$ and $V=v_n$. Let $\tilde\alpha=L_{\kt}^{-1}D_{\kt}v_n$ and $\tilde a=\partial_s\tilde\alpha+\kt v_n$. By the second form of $N_{\rm MS}$ in \cref{sec:semi} and \eqref{eq:semi-cancel}, we have $N_{\rm MS}(\kt;\Gamma^m)=-2\kt\tilde a-\mathfrak P[\tilde a]\,\partial_s\kt$, and \eqref{eq:semi-expand} becomes
\begin{equation}\label{eq:semi-closure}
\frac{\kc^{m+1}-\kt^m}\tau=N_{\rm MS}(\kt^m;\Gamma^m)\underbrace{-2\kt(a-\tilde a)-\mathfrak P[a-\tilde a]\,\partial_s\kt}_{=:d_2^m}+\underbrace{\rho}_{=:d_1^m},
\end{equation}
where $a=\partial_sv_t+\kc v_n$ is formed with the coefficient $\kc$ and $\tilde a$ with the coefficient $\kt$. Let now $\kt^{m+1}$ solve \eqref{eq:semi-reseq} at level $m+1$. Eliminating $\kc^{m+1}$ with \eqref{eq:semi-closure} and adding $\tau(\mathcal P^{m+1}-\mathcal P^m)\kt^{m+1}$ and $\tau\mathcal E^m\kt^m$ to both sides, we obtain \eqref{eq:recurrence} with
\begin{equation}\label{eq:semi-drows}
d^m=d_1^m+d_2^m+\underbrace{\mathcal E^m(\kt^m-\kt^{m+1})}_{=:d_3^m}+\underbrace{(\mathcal E^m-\mathcal E^{m+1})\kt^{m+1}}_{=:d_4^m}\underbrace{-(\mathcal P^{m+1}-\mathcal P^m)\kt^{m+1}}_{=:d_5^m}.
\end{equation}
This is the defect of \cref{prop:closure}.

\subsubsection*{The two estimates of one step}
The next lemma replaces \cite[Lemmas~4.2 and 4.3]{Rao2026}. Its two ingredients, \eqref{eq:semi-res} and \eqref{eq:semi-metric}, are compared at the end of \cref{sec:semi}.

\begin{lemma}[Uniform higher-order bounds]\label{lem:semi-Hr}
Assume (E). For every integer $r\ge3$, there are $\tau_H(r)$ and $C_r^{\rm H}$, with $0<\tau_H(r)\le\tau_A$ for the threshold $\tau_A$ of \cref{prop:semi-A}, such that, for $0<\tau\le\tau_H(r)$, the curves of \cref{prop:semi-A} satisfy
\[
\max_{m\tau\le T}\bigl(\|\kc^m\|_r+\|\kt^m\|_r\bigr)\le C_r^{\rm H},\qquad\tau\sum_{(m+1)\tau\le T}\|\kt^m\|_{r+3/2}^2\le C_r^{\rm H} .
\]
\end{lemma}

\begin{proof}
The curves are smooth by \cref{lem:semi-step}(e) and \cref{lem:semi-geom}. We argue by induction on $r$, the level $2$ being given by \cref{prop:semi-A} and \cref{lem:semi-step}(a). Let $r\ge3$, and let $K$ bound $\|\kc^m\|_{r-1}+\|\kt^m\|_{r-1}$ for $\tau\le\tau_H(r-1)$. Let $X=\|\kt^m\|_r$, $Y=\|\kt^m\|_{r+3/2}$ and $k=\|\kc^m\|_r$. By \cref{lem:semi-step}(d), we have
\begin{equation}\label{eq:semi-Hr-energy}
X^2+c\tau Y^2\le(1+C\tau)k^2+C\tau .
\end{equation}
By \cref{lem:semi-vel}(b) with $s=r-2$ and $s=r-1$, there is $A$, depending only on $K$, $r$ and the constants of $\mathcal U_{\rm ad}$, with
\[
\|b\|_{r-2}\le A(1+\|\kt^m\|_r),\qquad\|b\|_{r-1}\le A(1+\|\kt^m\|_{r+1})+\varpi,\qquad\varpi=A(1+k),
\]
since $\|\kc^m\|_{r-3/2}\le K$ and $\|\kc^m\|_{r-1/2}\le k$.
If moreover $\tau(1+X)\le1$, then \eqref{eq:semi-metric} gives $\|\kc^{m+1}\|_r\le(1+C\tau)X+C\tau Y+C\tau(1+A(1+k))$ with $C$ independent of $k$ and $X$. We square this bound. Young's inequality bounds the cross term $2(1+C\tau)C\tau XY$ by $\frac c4\tau Y^2+C\tau X^2$, the other cross terms by $C\tau(X^2+1+k^2)$, and the square $C^2\tau^2(Y+1+A(1+k))^2$ by $\frac c4\tau Y^2+C\tau(1+k^2)$ for $\tau\le c/(12C^2)$. Inserting \eqref{eq:semi-Hr-energy}, we obtain
\begin{equation}\label{eq:semi-Hr-rec}
\|\kc^{m+1}\|_r^2+\tfrac c2\tau\|\kt^m\|_{r+3/2}^2\le(1+C_\#\tau)\|\kc^m\|_r^2+C_\#\tau ,
\end{equation}
where $C_\#$, enlarged to dominate the constant of \eqref{eq:semi-Hr-energy}, depends only on $K$, $r$ and the constants of $\mathcal U_{\rm ad}$, and not on $\|\kc^m\|_r$.

Let $K_r^2=\eu^{C_\#T}(\|\kc(0)\|_r^2+C_\#T)$ and $\tilde K_r^2=(1+C_\#)K_r^2+C_\#$, and take $\tau_H(r)\le\tau_H(r-1)$ so small that $\tau_H(r)(1+\tilde K_r)\le1$ and the thresholds above hold, which is possible because $C_\#$ does not depend on $K_r$. If \eqref{eq:semi-Hr-rec} holds for the steps $j<m$, then \cref{lem:tool-gronwall}(a) gives $\|\kc^m\|_r\le K_r$, \eqref{eq:semi-Hr-energy} gives $\|\kt^m\|_r\le\tilde K_r$, so that $\tau(1+\|\kt^m\|_r)\le1$, and \eqref{eq:semi-Hr-rec} holds for the step $m$. Hence, \eqref{eq:semi-Hr-rec} holds for every step, and adding it over the steps bounds the sum of the dissipation terms.
\end{proof}

\subsection{Proof of \texorpdfstring{\cref{thm:semi-B}}{Theorem 3.2}}\label{sec:app-semi-first}

For a coefficient $f$ and a length $\ell$, we write $L_{f,\ell}=-\ell^{-2}\partial_\xi^2+f^2$ when the length matters, so that $L_f=L_{f,\ell}$ in normalized arclength. Let $(f,\ell)$ and $(g,\tilde\ell)$ satisfy $\|f\|_\sigma\le M$, $\|g\|_\sigma\le M$, $\|f\|_{L^2}\ge\frac12\kc_{\min}$, $\|g\|_{L^2}\ge\frac12\kc_{\min}$ and $\ell,\tilde\ell\in[\ell_{\min},\ell_{\max}]$. The difference $L_{f,\ell}^{-1}-L_{g,\tilde\ell}^{-1}$ equals $L_{f,\ell}^{-1}\bigl(g^2-f^2-(\tilde\ell^{-2}-\ell^{-2})\partial_\xi^2\bigr)L_{g,\tilde\ell}^{-1}$. By \cref{lem:setting-L}, $L_{f,\ell}^{-1}$ is bounded from $H^{-1}$ to $H^1$ uniformly in $(f,\ell)$. For $0<s\le\sigma-1$, the solution $u$ of $L_{f,\ell}u=F$ with $F\in H^{s-1}$ satisfies $\ell^{-2}u''=f^2u-F$, where $f^2u\in H^1\subset H^{s-1}$ because $H^\sigma$ is an algebra and $s<2$, so the bound extends to $\mathcal B(H^{s-1},H^{s+1})$. Since $\|g^2-f^2\|_\sigma\le C\|f-g\|_\sigma$, the resolvent identity gives, for $0\le s\le\sigma-1$,
\begin{equation}\label{eq:semi-resolvent}
\|L_{f,\ell}^{-1}-L_{g,\tilde\ell}^{-1}\|_{\mathcal B(H^{s-1},H^{s+1})}\le C\bigl(\|f-g\|_\sigma+|\ell-\tilde\ell|\bigr).
\end{equation}

Let $\Gamma^m$, $\kt^m$, $v^m$, $\Gamma^{m+1}$ and $\kt^{m+1}$ be as in the definition of the defect before \cref{lem:semi-stab}, where the defect $d^m$ of \eqref{eq:semi-drows} and the length defect $d_\ell^m$ are defined.

\begin{lemma}\label{lem:semi-defect}
In this setting, there is $\tau_D$ with $0<\tau_D\le\tau_R$, depending only on the constants of $\mathcal U_{\rm ad}$, such that $d^m\in L^2$ and the following hold for $\tau\le\tau_D$.
\begin{enumerate}
\item We have $\|d^m\|_{L^2}+|d_\ell^m|\le C\tau^{1/3}$.
\item If $\|\kt^m\|_{\sigma+3}+\|\kt^{m+1}\|_{\sigma+3}+\|\kc^m\|_{\sigma+3/2}\le K$, then $\|d^m\|_{L^2}+|d_\ell^m|\le C_K\tau$.
\end{enumerate}
\end{lemma}

\begin{proof}
The terms are bounded in the following table, where $\Delta\ell=\ell^{m+1}-\ell^m$ satisfies $|\Delta\ell|\le C\tau$ by \eqref{eq:semi-length}. The column (a) uses \cref{lem:semi-step}(a),(b), and the column (b) uses \cref{lem:semi-step}(c) with the hypothesis of (b).
\begin{center}
\footnotesize\setlength{\tabcolsep}{4pt}
\begin{tabular}{@{}llll@{}}
\toprule
Term & Bound & (a) & (b)\\
\midrule
$d_1^m$, \eqref{eq:semi-expand} & $C\tau$ & $\tau$ & $\tau$\\
$d_2^m$, coefficients $\kc^m\to\kt^m$ & $C\|\kc^m-\kt^m\|_1$ & $\tau^{(\sigma-1)/3}$ & $\tau$\\
$d_3^m=\mathcal E^m(\kt^m-\kt^{m+1})$ & $C\|\kt^m-\kt^{m+1}\|_1$ & $\tau^{(\sigma-1)/3}$ & $\tau$\\
$d_4^m=(\mathcal E^m-\mathcal E^{m+1})\kt^{m+1}$ & $C(|\Delta\ell|+\|\kc^{m+1}-\kc^m\|_{\sigma-1})$ & $\tau^{1/3}$ & $\tau$\\
$d_5^m=-(\mathcal P^{m+1}-\mathcal P^m)\kt^{m+1}$ & $C|\Delta\ell|\,\|\kt^{m+1}\|_3$ & $\tau^{\sigma/3}$ & $\tau$\\
$d_\ell^m$, \eqref{eq:semi-length} & $C\|\kc^m-\kt^m\|_{L^2}+C\tau\|b\|_{L^2}^2$ & $\tau^{\sigma/3}$ & $\tau$\\
\bottomrule
\end{tabular}
\end{center}
We take $\tau_D$ not larger than the thresholds of \cref{lem:semi-geom} for the cases used below, which depend only on the constants of $\mathcal U_{\rm ad}$, and so small that $\tau(1+\|\kt^j\|_3)\le1$ for $j=m,m+1$, which is possible by \cref{lem:semi-step}(b). For $d_2^m$, the factors $\kt$, $\partial_s\kt$ are bounded and $\mathfrak P$ maps $L^2$ to $H^1$, and $L_\kc(v_t-\tilde\alpha)=(D_\kc-D_{\kt})v_n-(\kc^2-\kt^2)\tilde\alpha$ is bounded in $H^{-1}$ by $C\|\kc-\kt\|_1$, so that $\|a-\tilde a\|_{L^2}\le C\|\kc^m-\kt^m\|_1$. The rows $d_3^m$ and $d_4^m$ carry the change of the coefficients in one step, and we split
\[
\kt^{m+1}-\kt^m=\underbrace{\kt^{m+1}-\kc^{m+1}}_{\text{new step}}+\underbrace{\kc^{m+1}-\kt^m}_{\text{metric step}},\qquad\kc^{m+1}-\kc^m=\underbrace{\kc^{m+1}-\kt^m}_{\text{metric step}}+\underbrace{\kt^m-\kc^m}_{\text{old step}} .
\]
The metric step is bounded in $H^1$ by $C\tau$ by \eqref{eq:semi-G3} with $k=0$ and $\theta=0$, whose hypotheses hold by \cref{lem:semi-vel}, and in $H^\sigma$ by a constant, hence in $H^{\sigma-1}$ by $C\tau^{1/(\sigma-1)}$. The new and the old step are bounded by \cref{lem:semi-step}(a), and the old step is bounded in $H^{\sigma-1}$ by $C\tau^{1/3}$, which is the rate of $d_4^m$ by \eqref{eq:semi-Eshape} with $g=\kt^{m+1}$. For $d_5^m$, the multiplier of $\mathcal P^{m+1}-\mathcal P^m$ is bounded by $C|\Delta\ell|\langle n\rangle^3$, and $\|\kt^{m+1}\|_3\le C\tau^{-(3-\sigma)/3}$. In the column (b), \cref{lem:semi-step}(c) with $s=\sigma-1$ and $s=1$ bounds the old and the new step by $C_K\tau$, and \eqref{eq:semi-G3} with $k=0$ and $\theta=\sigma-2$ bounds the metric step in $H^{\sigma-1}$ by $C_K\tau$.
\end{proof}

\begin{proof}[Proof of \cref{thm:semi-B}]
Let $\tau_B=\min(\tau_A,\tau_H(6),\tau_D)$ and $M_\tau=\lfloor T/\tau\rfloor$. The solvability, regularity and embeddedness of every step with $(m+1)\tau\le T$ are \cref{prop:semi-A}, which also identifies the quantities of \cref{thm:semi-B} with those of the curves $\hat x^m$. The bounds for the curvatures and the lengths are the stability estimate of \cref{lem:semi-stab}, now with forcings of order $\tau$.

For the curvatures and the lengths, \cref{lem:semi-Hr} with $r=6$ gives $\|\kc^m\|_6+\|\kt^m\|_6\le C$. Since $\sigma+3<6$, the finer form \cref{lem:semi-defect}(b) and \cref{lem:semi-exact}(e) give the defect bound $D=C\tau$, and \cref{lem:semi-step}(c) with $s=\sigma$ gives $\|\kc^m-\kt^m\|_\sigma\le C\tau$, hence $a_0=C\tau$, $\|\epsilon^0\|_\sigma\le C\tau$ and $\lambda_0=0$. \Cref{lem:semi-stab} with $m_*=M_\tau$ gives $\|\kt^m-\kc(t_m)\|_\sigma+|\ell^m-\ell(t_m)|\le C\tau$, and $\|\kc^m-\kt^m\|_\sigma\le C\tau$ gives the bound for $\kc^m$.

For the $C^1$ bound, we follow \cite[(4.21)--(4.26)]{Rao2026}. Let $\delta_m=|\hat x^m(0)-\hat X(0,t_m)|+|y^m(0)-Y(0,t_m)|$, where $\partial_\xi\hat X=\ell Y$. By \cref{lem:semi-geom}(b), $\hat x^{m+1}(0)=\hat x^m(0)+\tau v^m(0)$ and $y^{m+1}(0)=y^m(0)(1-\iu\tau b(0))+O(\tau^2)$. The exact marked point and tangent satisfy $\partial_t\hat X(0)=(V\bn+\alpha\bt)(0)$ and $\partial_tY=-\iu\omega Y$ with $\omega=\partial_sV-\kc\beta$ for $0<t\le T$. Indeed, $\partial_t\hat X=V\bn+\beta\bt$, and $\beta(0,t)=\alpha(0,t)$ by \eqref{eq:marked-point}. Differentiating $|\partial_\xi\hat X|^2=\ell^2$ in $t$ gives $\ell'=\partial_\xi\beta+\ell\kc V$, and therefore $\partial_t\bt=\ell^{-1}\partial_\xi(V\bn+\beta\bt)-\ell^{-1}\ell'\bt=\omega\bn$. Since $\bt=Y$ and $\bn=-\iu Y$, this is $\partial_tY=-\iu\omega Y$. Unlike Rao, we expand the exact quantities by Taylor's theorem at $t_{m+1}$ and compare the velocities of the step with $V$, $\alpha$ and $\omega$ at $t_{m+1}$. By \cref{thm:cont-N1}(c), \cref{prop:cont-tame}(d) and \eqref{eq:semi-resolvent}, these differences at $\xi=0$ are bounded by $C(\|\kt^m-\kc(t_{m+1})\|_\sigma+\|\kc^m-\kc(t_{m+1})\|_\sigma+|\ell^m-\ell(t_{m+1})|)\le C\tau$. Hence, $\delta_{m+1}\le(1+C\tau)\delta_m+C\tau^2$ with $\delta_0=0$, and \cref{lem:tool-gronwall}(a), applied to $\delta_m/\tau$, gives $\delta_m\le C\tau$. Finally, \cref{lem:tool-frame}(b) gives $\|\hat x^m-\hat X(t_m)\|_{C^1}\le C(\delta_m+|\ell^m-\ell(t_m)|+\|\kc^m-\kc(t_m)\|_{L^2})\le C\tau$.
\end{proof}

\subsection{Coercivity of the single layer under the capacity condition}\label{sec:app-cap}

The single curve case uses two results on the single layer of a smooth curve parametrized proportionally to arclength, the external input \cref{prop:ext-slcoer} and the Fredholm alternative \cref{prop:ext-slfred}.

\begin{lemma}\label{lem:K3-one}
Let $Y$ be a smooth embedded closed curve with $|Y'|\equiv\ell>0$, and let $\Gamma=Y(\Tor)$.
\begin{enumerate}[label=(\alph*)]
\item If $\langle\mathcal L_Y\varsigma,\varsigma\rangle\ge\mu\|\varsigma\|_{H^{-1/2}}^2$ for all $\varsigma\in L^2$ with some $\mu>0$, then $\mathcal L_Y$ maps $L^2$ boundedly onto $H^1$, and there is $B$ with $\|\varsigma\|_{L^2}\le B\|\mathcal L_Y\varsigma\|_{H^1}$ for all $\varsigma\in L^2$.
\item If there are a continuous function $w_{\rm eq}$ on $\Gamma$ and a number $\lambda_{\rm eq}>0$ with $\int_\Gamma w_{\rm eq}\,\dd s=1$ and $\SL w_{\rm eq}=\lambda_{\rm eq}$ on $\Gamma$, then there is $\mu>0$ with $\langle\mathcal L_Y\varsigma,\varsigma\rangle\ge\mu\|\varsigma\|_{H^{-1/2}}^2$ for all $\varsigma\in L^2$.
\end{enumerate}
In particular, under the hypothesis of (b), $Y$ satisfies \ref{K3} with constants $\mu$ and $B$ that depend on $Y$.
\end{lemma}

\begin{proof}
For (a), the curve lies in $\mathfrak F(\ell,c)$ for some $c>0$. The bound $\|\partial_\xi\mathcal L_Y\varsigma\|_{L^2}\le C\|\varsigma\|_{L^2}$ of \cite[Lemma~A.5]{EtoPaperI} and the square integrability of the kernel $\log|Y(\xi)-Y(\eta)|$ show that $\mathcal L_Y$ is bounded from $L^2$ to $H^1$. The coercivity makes $\mathcal L_Y$ injective on $L^2$, and \cref{prop:ext-slfred} gives the rest.

For (b), the parameter density $\varsigma_{\rm eq}=\ell\,w_{\rm eq}\circ Y\in L^2$ satisfies $\int_\Tor\varsigma_{\rm eq}\,\dd\xi=1$ and $\mathcal L_Y\varsigma_{\rm eq}=\lambda_{\rm eq}$. We split $\varsigma\in L^2$ as $\varsigma=\varsigma_0+a_0\varsigma_{\rm eq}$ with $a_0=\int_\Tor\varsigma\,\dd\xi$. Since $\mathcal L_Y$ is symmetric and $\langle\mathcal L_Y\varsigma_{\rm eq},\varsigma_0\rangle=\lambda_{\rm eq}\int_\Tor\varsigma_0\,\dd\xi=0$, we have
\begin{equation}\label{eq:cap-split}
\langle\mathcal L_Y\varsigma,\varsigma\rangle=\langle\mathcal L_Y\varsigma_0,\varsigma_0\rangle+\lambda_{\rm eq}\,a_0^2 .
\end{equation}
By \cref{prop:ext-slcoer}, the form $\langle\mathcal L_Y\varsigma_0,\varsigma_0\rangle$ is bounded below by $c_{\rm sl}\|\varsigma_0\|_{H^{-1/2}}^2$ with some $c_{\rm sl}>0$, and $\|\varsigma\|_{H^{-1/2}}^2\le2\|\varsigma_0\|_{H^{-1/2}}^2+2a_0^2\|\varsigma_{\rm eq}\|_{H^{-1/2}}^2$. Hence, \eqref{eq:cap-split} gives (b) with $\mu=\frac12\min\bigl(c_{\rm sl},\lambda_{\rm eq}/(1+\|\varsigma_{\rm eq}\|_{H^{-1/2}}^2)\bigr)$. The last claim combines (a) and (b).
\end{proof}

\subsection{The multiplier in the tangential row}\label{sec:app-p52-tan}

With the notation of \cref{sec:app-p52}, constants $C$ depend only on the constants of $\mathscr K\ni x$ (\cref{lem:Kfamily}) and on the uniform bounds of finitely many derivatives of $x$, $v$, $\kt$ and of $|x'|^{-1}$ (\cref{lem:semi-material}). We use the frame of \cref{sec:ms-flow}, the operators $\Lop=\Lop_\kc$ and $\Dop=\Dop_\kc$ of \cref{sec:setting}, and write $e_n=e\cdot\bn$, $e_t=e\cdot\bt$. The tangential part $\eta(e)=e_t-\Lop^{-1}\Dop e_n$ is that of the adapted norm of \cref{sec:fd}. The normal components of $\partial_s^2F=-\kt\bn$ and $\partial_s^2x=-\kc\bn$ give the analogue of \cite[(5.10)]{Rao2026},
\begin{equation}\label{eq:p52-curv}
\kc-\kt=\tau\,\partial_s^2v\cdot\bn .
\end{equation}
A smooth solution of the step \eqref{eq:semi} satisfies $\kt=c_{\SL}-2\SL v_n$ with $c_{\SL}=c+2\avg{\SL v_n}$ (\cref{lem:p52-gauge}). With $q_e(\xi,\eta)=(e(\xi)-e(\eta))/(x(\xi)-x(\eta))$ in complex notation, let
\[
f_{\mathcal L}[e]\varsigma(\xi)=-\frac1{2\pi}\int_\Tor\operatorname{Re}q_e(\xi,\eta)\,\varsigma(\eta)\,\dd\eta
\]
be the derivative of the parameter form of the single layer in the direction $e$ (\cref{lem:p52-dSL}), and let $\varsigma_v=v_n|x'|$ and $\Jm$ be as in \cref{sec:scheme}. By \cref{lem:p52-diff}, for a direction $e\in H^1(\Tor;\R^2)$, the derivative $u=D\Ft(x)[e]$ and some $\delta c\in\R$ satisfy \eqref{eq:p52-lin} for all $\varphi\in H^1(\Tor;\R^2)$, with $w=(u-e)/\tau$ and with the forms $\bx^{\SL}$, $\ax'[e]$ and $\bx^{\SL\prime}[e]$ defined there, and $\int w_n\,\dd s=-\int v\cdot\Jm e'\,\dd\xi$.

By \cref{lem:p52-cov}, the step and its derivative commute with the change to normalized arclength, and we compute for a curve with $|x'|\equiv\ell$. Then, $\partial_s=\ell^{-1}\partial_\xi$, $\dd s=\ell\,\dd\xi$ and $\mathcal L_x(v\cdot\Jm e')=\SL g_e$ with $g_e=\ell^{-1}v\cdot\Jm e'$, and the bounds of $x$, $v$ and $\kt$ remain uniform by \cref{lem:semi-material}. The variations of the geometry are $\delta\bn=-(\bn\cdot\partial_se)\bt$ and $\delta(\dd s)=(\bt\cdot\partial_se)\,\dd s$, with $\bt\cdot\partial_se=\partial_se_t+\kc e_n$ and $\bn\cdot\partial_se=\partial_se_n-\kc e_t$. We test \eqref{eq:p52-lin} with $\varphi=\psi\bt$ and $\varphi=\psi\bn$ for real smooth $\psi$, first for smooth $e$ and then, by continuity in $H^{-1}$, for $e\in H^1$.

For $\varphi=\psi\bt$, one has $\varphi_n=0$, $\varphi\cdot\Jm e'=-(\bn\cdot\partial_se)\psi\,\ell$ and $\partial_s(\psi\bt)=(\partial_s\psi)\bt-\kc\psi\bn$. Hence, the kernel and datum variations vanish, and $\ax(w,\psi\bt)=\langle\Lop w_t-\Dop w_n,\psi\rangle$. With $2\SL v_n=c_{\SL}-\kt$, the test variation is $-\int(c_{\SL}-\kt)(\bn\cdot\partial_se)\psi\,\dd s$, while the right side of \eqref{eq:p52-lin} is $-c_{\SL}\int(\bn\cdot\partial_se)\psi\,\dd s$. The multiplier $c_{\SL}$ thus appears with the same sign on both sides and cancels exactly. With $u=e+\tau w$ and $\partial_sF=\bt+\tau\partial_sv$, the remaining terms are
\[
\begin{aligned}
&\tau\langle\Lop w_t-\Dop w_n,\psi\rangle+\underbrace{\ax(e,\psi\bt)-\int(\bt\cdot\partial_se)\,\bt\cdot\partial_s(\psi\bt)\,\dd s}_{\text{frame part}}\\
&\qquad+\int\kt(\bn\cdot\partial_se)\psi\,\dd s
-\tau\int(\bt\cdot\partial_se)\,\partial_sv\cdot\partial_s(\psi\bt)\,\dd s=0 .
\end{aligned}
\]
Since $\partial_se=(\bt\cdot\partial_se)\bt+(\bn\cdot\partial_se)\bn$, the frame part equals $-\int\kc(\bn\cdot\partial_se)\psi\,\dd s$, and by \eqref{eq:p52-curv} it combines with the term with $\kt$ to $-\tau\int(\partial_s^2v\cdot\bn)(\bn\cdot\partial_se)\psi\,\dd s$. This is the cancellation of \cite[(5.13)]{Rao2026}, where $\kt$ is replaced by $-v_n$. Every remaining term carries a factor $\tau$, and division by $\tau$ gives
\begin{equation}\label{eq:p52-tan}
\Lop w_t-\Dop w_n=\mathcal G_xe,\qquad
\eta(u)=\eta(e)+\tau\Lop^{-1}\mathcal G_xe,
\end{equation}
with
\[
\langle\mathcal G_xe,\psi\rangle=\int(\partial_s^2v\cdot\bn)
(\bn\cdot\partial_se)\psi\,\dd s+\int(\bt\cdot\partial_se)
\bigl[(\partial_sv\cdot\bt)\partial_s\psi-\kc(\partial_sv\cdot\bn)\psi\bigr]\,\dd s .
\]
Since $v$ is smooth, $\mathcal G_x$ is bounded from $H^1$ to $H^{-1}$, and $\|\eta(u)-\eta(e)\|_{H^1}\le C\tau\|e\|_{H^1}$.

The normal row is the resolvent step \eqref{eq:p52-res} for the operator $\mathcal H=\frac12\Lam_\Gamma\Sop$, and the proof of \cref{prop:stab} combines it with \eqref{eq:p52-tan}.

\begin{proof}[Proof of \cref{prop:stab}]
Let $\tau_P\le\min(1,\tau_X)$, with the threshold $\tau_X$ of \cref{lem:p52-adnorm}, be so small that $(1+c_b)\tau_P^{2/3}\|\mathcal H_1\|_{H^1\to L^2}\le\frac12$ for all iterates, which is possible because $b$ is bounded below and $\|\mathcal H_1\|_{H^1\to L^2}$ is bounded uniformly in $m$, and that $C\tau_P\le\frac12$ for the constant $C$ of the energy bound for $\vartheta$ below. Let first $\hat x$ be the curve of constant speed, and let $u=D\Ft(\hat x)e$. By \eqref{eq:p52-res}, \cref{lem:p52-resolvent} and \cref{lem:p52-forcing}, we have $u_n=\vartheta+r$ with $\vartheta=(I+\tau\mathcal H)^{-1}(e_n+\tau f[e])$ and $\|r\|_{H^1}\le C\tau\|e\|_{H^1}$. We write $(I+\tau\mathcal H_0)\vartheta=e_n+\tau q$ with $q=f[e]-\mathcal H_1\vartheta$, so that $\|q\|_{L^2}\le C(\|e\|_{H^1}+\|\vartheta\|_{H^1})$. The forcing $q$ is only in $L^2$, and the dissipation of the multiplier absorbs it through the energy bound \eqref{eq:tool-energy} of \cref{lem:tool-mult}(c), which gives $\|\vartheta\|_{H^1}^2\le\|e_n\|_{H^1}^2+C\tau(\|\vartheta\|_{H^1}^2+\|e\|_{H^1}^2)$. With $(a+b)^2\le(1+\tau)a^2+(1+\tau^{-1})b^2$, this bound and \eqref{eq:p52-tan} give
\[
\|u_n\|_{H^1}^2\le(1+C\tau)\|e_n\|_{H^1}^2+C\tau\|e\|_{H^1}^2,\qquad\|\eta(u)\|_{H^1}^2\le(1+\tau)\|\eta(e)\|_{H^1}^2+C\tau\|e\|_{H^1}^2 ,
\]
and the sum of the two bounds and \cref{lem:p52-adnorm} give $\xnorm[\hat x]{u}\le(1+C\tau)\xnorm[\hat x]{e}$. By \cref{lem:p52-cov}, the same bound holds for $D\Ft(x^m)$ in the adapted norm of $x^m$.

It remains to change the level of the adapted norm. Since $x^{m+1}=x^m+\tau v^m$ with $v^m$ bounded in $W^{3,\infty}$, the frames and the curvature change by $O(\tau)$ in $W^{1,\infty}$ at corresponding material points, and the identity $\Lop_{m+1}^{-1}\Dop_{m+1}-\Lop_m^{-1}\Dop_m=\Lop_{m+1}^{-1}(\Dop_{m+1}-\Dop_m)+\Lop_{m+1}^{-1}(\Lop_m-\Lop_{m+1})\Lop_m^{-1}\Dop_m$, with $\Lop_j=\Lop_{\kc^j}$ and $\Dop_j=\Dop_{\kc^j}$, bounds the change of $\Lop^{-1}\Dop$ on $H^1$ by $C\tau$. The normalized arclengths of the two curves differ by a reparametrization $\phi$ with $\|\phi-\mathrm{id}\|_{C^1}\le C\tau$, which costs a further factor $1+C\tau$. Hence, $\xnorm[x^{m+1}]{u}\le(1+C\tau)\xnorm[x^m]{u}+C\tau\|u\|_{H^1}$, which is the norm comparison after \cite[(5.20)]{Rao2026}. We apply it to $u=D\Ft(x^m)e$, together with the bound at the level of $x^m$, and use \cref{lem:p52-adnorm} only on the term of order $\tau$, which gives \eqref{eq:four}. The extension of $D\Ft(x^m)$ to $H^1$ is that of \cref{lem:p52-diff}.
\end{proof}

\subsection{Consistency, comparison of the linearized steps and the nonlinear remainder}\label{sec:app-fd}

This subsection proves \cref{lem:fd-cons,lem:fd-comp,lem:fd-rem}, following \cite[\S\S5--6]{Rao2026}. The mass is the boundary element mass $\bXh(u,\varphi)=2(\Om\varphi)^\T\Lhp(\Om u)$. The linearized continuum step is rewritten in the form of the mass $\bx$ with $\Lp$ (\cref{lem:fd-lin}), with the derivative of the Dirichlet-to-Neumann map along piecewise linear directions (\cref{lem:fd-shape}).

Throughout, $x=x^m$ is a time-semidiscrete iterate in its own parameter, with $(m+1)\tau\le T$, $F=\Ft(x)$, $v=(F-x)/\tau$, and all constants are uniform in $m$ by \cref{lem:semi-material}, which also gives $|x'|\ge(C_1^{\rm s})^{-1}$. We write $X=\Ih x$, $F_h=\Fht(X)$ for the solution of \eqref{eq:step-var} with $X^m=X$, and $E=F_h-\Ih F$. Norms are taken in the parameter. The boundary element results are applied with $x\in\mathscr K$ (\cref{lem:Kfamily}) and the radii $r_0$, $h_0$ of \cref{sec:input-nbhd}. For degree one, Rao's Gauss--Lobatto nodes are the vertices $\xi_j$, and by \cref{lem:semi-material}, $x$ satisfies his hypotheses on the time-semidiscrete solution. We use his estimates only for the stiffness form and for interpolation. For smooth $F$ and P1 functions $f$, $g$, $e_h$, $\varphi$, $E$, they give
\begin{gather}
|a_{\Ih x}(\Ih F,\varphi)-\ax(F,\varphi)|\le Ch^2\|\varphi\|_{H^1},\label{eq:fd-rao-a}\\
\begin{gathered}
|a_{\Ih x}(f,g)-\ax(f,g)|\le Ch\|f\|_{H^1}\|g\|_{H^1},\qquad
|a'_{\Ih x}[e_h](\Ih F,\varphi)-a'_x[e_h](F,\varphi)|\le Ch\|e_h\|_{H^1}\|\varphi\|_{H^1},\\
|a'_{\Ih x}[e_h](E,\varphi)|\le C\|e_h'\|_{L^2}\|E'\|_{L^\infty}\|\varphi'\|_{L^2},
\end{gathered}\label{eq:fd-rao-b}
\end{gather}
where the first is the stiffness part of the proof of \cite[Lemma~5.1]{Rao2026}, which rests on \cite[(5.3)--(5.4)]{Rao2026}, proved there including degree one, and the others are the stiffness parts of \cite[(5.25)]{Rao2026} and of the bounds after \cite[(5.27)]{Rao2026}.

The following facts are used for P1 tests $\varphi$. Loads $y\in\R^N$ are measured in the norms $\dnorm{y}\ge\|y\|_{(H^1)^*}$ of \cref{sec:input-nbhd}. The lumped load $M_vf$ of a P1 function $f$ and the consistent load of $f\in L^2$ satisfy $\dnorm{y}\le C\|f\|_{L^2}$, $\dnorm{\Om\varphi}\le C\|\varphi\|_{L^2}$ by \eqref{eq:use-Om}, and $\|\Lhp\Om\varphi\|_{H^1}\le C\|\varphi\|_{L^2}$ by \cite[Lemma~A.20]{EtoPaperI} (\cref{sec:input-results}).
\begin{enumerate}[label=(F\arabic*)]
\item\label{item:fd-conv-kernel} The estimates are taken on the discrete constraint kernel $\{\one^\T\Om(X)\varphi=0\}$. A constant added to a trace that is paired with $\Om\varphi$ does not contribute there, and a constant $\bar c$ paired with $\int\varphi\cdot\bn\,\dd s=\int\varphi\cdot\Jm x'\,\dd\xi$ contributes a term bounded by $Ch^2|\bar c|\,\|\varphi\|_{H^1}$, by \eqref{eq:use-area} and an integration by parts.
\item\label{item:fd-conv-p1} For a P1 vector $r_h$, $|(\Om\varphi)^\T r_h|\le C\|\varphi\|_{L^2}\|r_h\|_{L^2}$, since $|\omega_j|\le\mu_j$. For a P1 function $w$ and for $w\in W^{2,\infty}$, respectively,
\[
\Big|(\Om\varphi)^\T w-\int w\,\varphi\cdot\bn\,\dd s\Big|\le Ch\|w\|_{L^2}\|\varphi\|_{H^1},\qquad
\Big|(\Om\varphi)^\T\Ih w-\int w\,\varphi\cdot\bn\,\dd s\Big|\le Ch^2\|w\|_{W^{2,\infty}}\|\varphi\|_{L^2}.
\]
The first follows from $(\Om\varphi)^\T w=\int\Ih(w\varphi)\cdot\Jm X'\,\dd\xi$ and an inverse estimate, and the second is the vertex rule, since the first moments of the hat functions about their vertices vanish.
\item\label{item:fd-conv-load} For every load $y$, $|(\Om\varphi)^\T\Lhp y|\le C\|\varphi\|_{L^2}\dnorm{y}$ by \eqref{eq:use-load}, and $|(\Om\varphi)^\T\Lhp y|=|y^\T\Lhp\Om\varphi|\le C\|\varphi\|_{L^2}\|y\|_{(H^1)^*}$ by the symmetry of $\Lhp$.
\item\label{item:fd-conv-inv} Let $y$ be the consistent load, with $\rho=|x'|$, of $f\in L^2$, and let $\one^\T y=0$. Then, we have
\[
\Big|(\Om\varphi)^\T\Lhp y-\int(\Lp f)\,\varphi\cdot\bn\,\dd s\Big|\le Ch^{1/2}\|f\|_{L^2}\|\varphi\|_{H^1}.
\]
This follows from \cref{lem:bem-invcons}, which gives $\Lhp y=\Lp f+\bar c+O_{H^{1/2}}(h^{1/2}\|f\|_{L^2})$ with $|\bar c|\le C\|f\|_{L^2}$, from \ref{item:fd-conv-p1} for the P1 function $\Lhp y$, and from \ref{item:fd-conv-kernel} for $\bar c$.
\end{enumerate}

The consistency and the comparison use the boundary element defect. The rows of $\Om=\Om(X)$ are central-difference normals, and Taylor expansion at $\xi_j$ gives
\begin{equation}\label{eq:fd-mu}
\omega_j/\mu_j=\bn(\xi_j)+O(h^2),\qquad\mu_j=h|x'(\xi_j)|+O(h^3).
\end{equation}
Hence, $\Om\Ih v=M_v\Ih v_n+g_h$, where $g_h$ is the lumped load of a P1 function of size $O(h^2)$, so that $\dnorm{g_h}\le Ch^2$. Let $d_h=\Lh\Ih\kt-M_v\Ih(\Lam_\Gamma\kt)$ be the boundary element defect, where $\Lam_\Gamma$ is the Dirichlet-to-Neumann map of \eqref{eq:semi}, which is the map of \cref{lem:in-bweak} by \cref{prop:ext-dtn,prop:ext-layer}. By \cref{lem:in-bweak} and $\|\kt\|_{H^3}\le C$ (\cref{lem:semi-material}),
\[
\|d_h\|_{(H^1)^*}\le Ch^2,\qquad\dnorm{d_h}\le Ch^{3/2}.
\]
Since $v_n=-\tfrac12\Lam_\Gamma\kt$ and $\Lhp\Lh=P$, we have $\Om\Ih v=-\tfrac12(\Lh\Ih\kt-d_h)+g_h$, and \eqref{eq:in-Lhp} gives
\begin{equation}\label{eq:fd-LOmv}
\Lhp\Om\Ih v=-\tfrac12(\Ih\kt-c'\one)+\mathsf y_h,\qquad \mathsf y_h=\Lhp(\tfrac12d_h+g_h),\qquad\|\mathsf y_h\|_{H^{1/2}}\le Ch^{3/2},
\end{equation}
where $c'$ is the vertex mean of $\Ih\kt$.

\begin{proof}[Proof of \cref{lem:fd-cons}]
As in the proof of \cite[Lemma~5.1]{Rao2026}, we insert $\Ih F=X+\tau\Ih v$ into \eqref{eq:step-var} and subtract the weak form $\ax(F,\varphi)=\int\kt\,\varphi\cdot\bn\,\dd s$ of \eqref{eq:semi}. For P1 $\varphi$ on the kernel, the multiplier term vanishes, the factors $\tau^{-1}$ and $\tau$ of the mass cancel, and by \eqref{eq:fd-LOmv} the term $c'\one$ vanishes as well. Hence, the residual is
\[
r(\varphi)=\underbrace{a_X(\Ih F,\varphi)-\ax(F,\varphi)}_{\text{stiffness}}+\underbrace{\int\kt\,\varphi\cdot\bn\,\dd s-(\Om\varphi)^\T\Ih\kt}_{\text{vertex rule}}+\underbrace{2(\Om\varphi)^\T \mathsf y_h}_{\text{defect}} .
\]
The stiffness is bounded by $Ch^2\|\varphi\|_{H^1}$ by \eqref{eq:fd-rao-a}, and the vertex rule by $Ch^2\|\varphi\|_{L^2}$ by \ref{item:fd-conv-p1}, since $\|\kt\|_{W^{2,\infty}}\le C\|\kt\|_{H^3}$. The defect is bounded by $C\|\varphi\|_{L^2}(\|d_h\|_{(H^1)^*}+\|g_h\|_{(H^1)^*})\le Ch^2\|\varphi\|_{L^2}$ by the second bound of \ref{item:fd-conv-load}.

For the constraint, $\one^\T\Om E=-\tau\one^\T\Om\Ih v=-\tau(\tfrac12\one^\T d_h+\one^\T g_h)$, since $\one^\T\Lh=0$, and $\|\one\|_{H^1}=1$ gives $|\one^\T\Om E|\le C\tau h^2$. The error satisfies $(a_X+\tau^{-1}b_{X,h})(E,\varphi)=-r(\varphi)$ on the kernel, and the stability bound \eqref{eq:use-brezzi} with $\mathsf R=Ch^2$ and $\delta=C\tau h^2$ gives $\|E\|_{H^1}\le Ch^2$.
\end{proof}

\subsubsection*{The derivative of the Dirichlet-to-Neumann map}
The comparison uses the derivative of the Dirichlet-to-Neumann map along piecewise linear directions $e_h$. For a curve $Y\in\mathfrak C$ and $\mathfrak g\in H^1$, let $(\varsigma_Y,c_Y)$ be the gauge pair of $\mathfrak g$ on $Y$ (\cref{sec:input-results}), that is, $\mathcal L_Y\varsigma_Y=\mathfrak g+c_Y$ and $\int\varsigma_Y\,\dd\xi=0$. For $\mathfrak g\in H^2$ and $|\varepsilon|\,\|e_h'\|_{L^\infty}\le2r_1$, we have $x+\varepsilon e_h\in\mathfrak C$ and $\Lam_{x+\varepsilon e_h}\mathfrak g=\varsigma_{x+\varepsilon e_h}/|x'+\varepsilon e_h'|$ as functions of the parameter. The derivative therefore comes from two places, the single layer, whose derivative $f_{\mathcal L}$ is bounded in $H^1$ along piecewise linear directions by \cref{lem:p52-dSL}, and the speed.

\begin{lemma}[The derivative of the Dirichlet-to-Neumann map]\label{lem:fd-shape}
Let $e_h$ be piecewise linear and $\mathfrak g\in H^2$. Let $\varsigma_0$ be the gauge density of $\mathfrak g$ on $x$, and let $\varsigma_1 $ be the gauge density of $-f_{\mathcal L}[e_h]\varsigma_0$ on $x$, with $f_{\mathcal L}$ as in \cref{lem:p52-dSL}. Then, the difference quotients of $\varepsilon\mapsto\Lam_{x+\varepsilon e_h}\mathfrak g$ converge in $L^2$ as $\varepsilon\to0$ to
\[
D\Lam[e_h]\mathfrak g=\frac{\varsigma_1}{|x'|}-\varsigma_0\,\frac{x'\cdot e_h'}{|x'|^3}.
\]
Moreover, the following hold, with constants independent of $h$ and $m$.
\begin{enumerate}
\item\label{item:fd-shape-lin} $D\Lam[e_h](\alpha \mathfrak g+\beta)=\alpha D\Lam[e_h]\mathfrak g$ for $\alpha,\beta\in\R$.
\item\label{item:fd-shape-bd} $\|D\Lam[e_h]\mathfrak g\|_{L^2}\le C\|e_h\|_{H^1}\|\mathfrak g\|_{H^2}$.
\item\label{item:fd-shape-load} Let $\hat f=D\Lam[e_h]\mathfrak g+(\Lam_\Gamma \mathfrak g)\,\bt\cdot\partial_se_h$. Then, the derivative of $\varepsilon\mapsto\int(\Lam_{x+\varepsilon e_h}\mathfrak g)\,\varphi_j\,|x'+\varepsilon e_h'|\,\dd\xi$ at $\varepsilon=0$ is $\int\hat f\varphi_j\,|x'|\,\dd\xi$ for every $j$, and $\int\hat f\,\dd s=0$.
\end{enumerate}
\end{lemma}

\begin{proof}
Let $x_\varepsilon=x+\varepsilon e_h$, and let $(\varsigma_\varepsilon,c_\varepsilon)$ be the gauge pair of $\mathfrak g$ on $x_\varepsilon$. Since $\Lam_\Gamma$ maps $H^2$ to $H^1$ (\cref{thm:cont-N1}(c)), we have $\varsigma_0=|x'|\Lam_\Gamma\mathfrak g\in H^1$ and $\|\varsigma_0\|_{H^1}\le C\|\mathfrak g\|_{H^2}$. Subtracting the gauge equations shows that $\varsigma_\varepsilon-\varsigma_0-\varepsilon\varsigma_1$, which has mean zero, is the gauge density on $x_\varepsilon$ of the sum of the Taylor remainder $-[(\mathcal L_{x_\varepsilon}-\mathcal L_x)\varsigma_0-\varepsilon f_{\mathcal L}[e_h]\varsigma_0]$ and the cross term $-\varepsilon(\mathcal L_{x_\varepsilon}-\mathcal L_x)\varsigma_1$. Both are $O(\varepsilon^2)$ in $H^1$, by \cref{lem:p52-dSL}\ref{item:p52-dSL-b} and by the Lipschitz bound of \cite[Lemma~A.7]{EtoPaperI}, so the stability of the gauge representation on $\mathfrak C$ gives $\|\varsigma_\varepsilon-\varsigma_0-\varepsilon\varsigma_1\|_{L^2}\le C_{e_h}\varepsilon^2$. The expansion $|x'+\varepsilon e_h'|^{-1}=|x'|^{-1}-\varepsilon(x'\cdot e_h')|x'|^{-3}+O(\varepsilon^2)$, which is uniform because $e_h'$ is bounded, then gives the formula for $D\Lam[e_h]\mathfrak g$.

Part~\ref{item:fd-shape-lin} holds because the gauge pair of the constant function $1$ is $(0,-1)$ on every curve. For part~\ref{item:fd-shape-bd}, the stability of the gauge representation and \cref{lem:p52-dSL}\ref{item:p52-dSL-a} give $\|\varsigma_1\|_{L^2}\le C\|f_{\mathcal L}[e_h]\varsigma_0\|_{H^1}\le C\|e_h\|_{H^1}\|\varsigma_0\|_{H^1}$, and the second term of $D\Lam[e_h]\mathfrak g$ is bounded by $C\|\varsigma_0\|_{L^\infty}\|e_h'\|_{L^2}$. For part~\ref{item:fd-shape-load}, the derivative of $|x'+\varepsilon e_h'|$ at $\varepsilon=0$ is $(\bt\cdot\partial_se_h)|x'|$, and the product rule applies because the difference quotients converge in $L^2$. Finally, $\int(\Lam_{x_\varepsilon}\mathfrak g)|x_\varepsilon'|\,\dd\xi=\int\varsigma_\varepsilon\,\dd\xi=0$ for every $\varepsilon$, which gives $\int\hat f\,\dd s=0$.
\end{proof}

\subsubsection*{The linearized steps}\label{sec:fd-steps}
\Cref{lem:p52-diff} gives the linearized continuum step in the $\SL$-form. The discrete step is written with $\Lhp$, so the comparison needs the continuum step in the form of the mass $\bx$ with $\Lp=\Pz\SL\Pz$. The next lemma makes this change and gives the multiplier in closed form. For a direction $e$, let $D\bn[e]=\Jm e'/|x'|-\bn\,\bt\cdot\partial_se$ be the derivative of the unit normal, recall that $D(\bn\,\dd s)[e]=\Jm e'\,\dd\xi$, and let $g_e=v\cdot\Jm e'/|x'|$, so that $v\cdot D\bn[e]=g_e-v_n\,\bt\cdot\partial_se$.

\begin{lemma}[The linearized continuum step]\label{lem:fd-lin}
Let $e_h$ be piecewise linear and $u=D\Ft(x)e_h$. Then, $\int(u-e_h)\cdot\bn\,\dd s+\tau\int v\cdot\Jm e_h'\,\dd\xi=0$, and for all $\varphi\in H^1(\Tor;\R^2)$,
\begin{equation}\label{eq:fd-lin}
\ax(u,\varphi)+\tau^{-1}\bx(u-e_h,\varphi)+\ax'[e_h](F,\varphi)+\bx'[e_h](v,\varphi)=(\delta c+c_P)\int\varphi\cdot\bn\,\dd s+c\int\varphi\cdot\Jm e_h'\,\dd\xi,
\end{equation}
where $c=\avg{\kt}$ is the multiplier of the continuum step, $\ax'[e_h](F,\varphi)=-\int(\bt\cdot\partial_se_h)\,\partial_sF\cdot\partial_s\varphi\,\dd s$ is the stiffness variation of \cite[(5.11)]{Rao2026},
\[
\bx'[e_h](v,\varphi)=\underbrace{-2\int\big(\Lp D\Lam[e_h]\Lp v_n\big)\varphi_n\,\dd s}_{\text{DtN term}}+\underbrace{2\int\Lp\big(v\cdot D\bn[e_h]\big)\varphi_n\,\dd s}_{\text{second-slot frame term}}+\underbrace{2\int(\Lp v_n)\,\varphi\cdot\Jm e_h'\,\dd\xi}_{\text{first-slot frame term}} ,
\]
and
\[
\ell\,\delta c=\ax(u,\bn)+\ax'[e_h](F,\bn)-\avg{\kt}\int\bn\cdot\Jm e_h'\,\dd\xi,\qquad c_P=\frac2\ell\int(\Lp v_n)\,\bn\cdot\Jm e_h'\,\dd\xi .
\]
Moreover, $|\delta c|+|c_P|\le C\|e_h\|_{H^1}$.
\end{lemma}

\begin{proof}
We start from \cref{lem:p52-diff} with $e=e_h$, which gives the linearized constraint and \eqref{eq:p52-lin} in the $\SL$-form, with the kernel variation given by the integral $f_{\mathcal L}[e_h]\varsigma_v$. Let $W=u-e_h$ and $g=g_{e_h}$. The passage from $\SL$ to $\Lp$ uses that $\Lp f$ and $\SL(f-\avg f)$ differ by a constant, and each such constant produces a term $\lambda\int\varphi_n\,\dd s$ with a number $\lambda$ independent of $\varphi$. The linearized constraint gives $\avg{W_n}=-\tau\avg g$, so $\tau^{-1}\bx(W,\varphi)$ is $\tau^{-1}\bx^{\SL}(W,\varphi)$ plus the mean term $2\avg g\int(\SL1)\varphi_n\,\dd s$. \Cref{lem:fd-shape} with $\mathfrak g=\Lp v_n$, whose gauge density is $\varsigma_v$ because $\int v_n\,\dd s=0$, writes the DtN term as the kernel variation of $\bx^{\SL\prime}$ plus $2\int\Lp(v_n\,\bt\cdot\partial_se_h)\varphi_n\,\dd s$. This term and the second-slot frame term add up to $2\int(\Lp g)\varphi_n\,\dd s$, which is the datum variation of $\bx^{\SL\prime}$ minus the mean term, so the mean terms cancel. The first-slot frame term is the test variation of $\bx^{\SL\prime}$ minus $2\avg{\SL v_n}\int\varphi\cdot\Jm e_h'\,\dd\xi$. With \eqref{eq:p52-lin} and $c_{\SL}-2\avg{\SL v_n}=\avg{\kt}$ (\cref{lem:p52-gauge}), we obtain \eqref{eq:fd-lin} with some number $\lambda$ in place of $\delta c+c_P$.

The test $\varphi=\bn$ identifies $\lambda=\delta c+c_P$, since $\bx(W,\bn)$, the DtN term and the second-slot frame term vanish for it, because $\Lp$ maps into functions of arclength mean zero. The bounds follow from $\|u\|_{H^1}\le C\|e_h\|_{H^1}$, which holds by \cref{prop:stab} and \cref{lem:p52-adnorm}, and from the boundedness of $\Lp v_n=-\frac12(\kt-\avg{\kt})$.
\end{proof}

The discrete linearized system seeks a P1 field $u_h$ and $\delta c_h\in\R$ with, for every P1 $\varphi$,
\begin{equation}\label{eq:fd-lin-h}
\begin{aligned}
&a_X(u_h,\varphi)+\tau^{-1}b_{X,h}(u_h-e_h,\varphi)+a_X'[e_h](F_h,\varphi)+\tau^{-1}b_{X,h}'[e_h](F_h-X,\varphi)\\
&\qquad=\delta c_h\one^\T\Om\varphi+c_h\one^\T\Om(e_h)\varphi,
\end{aligned}
\end{equation}
together with $\one^\T\Om(u_h-e_h)+\one^\T\Om(e_h)(F_h-X)=0$. The derivative $D\Lh[e_h]$ of $\Lh(X)$ in the direction $e_h$ exists because $\Lh$ is differentiable in the vertices (\cref{sec:input-results}), and the system has a unique solution, linear in $e_h$, by \cref{cor:bem-acoerc} and Brezzi's theorem. Since $\Om$ is linear in $X$ and $\Lhp D\Lh[e_h]$ vanishes on $\one$ and maps into $\one^\perp$, the mass variation is
\[
b_{X,h}'[e_h](f,\varphi)=\underbrace{2\big(\Om(e_h)\varphi\big)^\T\Lhp\Om f}_{\text{first-slot frame term}}+\underbrace{2(\Om\varphi)^\T\Lhp\Om(e_h)f}_{\text{second-slot frame term}}\underbrace{-2(\Om\varphi)^\T\Lhp D\Lh[e_h]\Lhp\Om f}_{\text{operator derivative}} .
\]

With $a_e=e_h\cdot\bt$ and $F_s=\partial_sF$, let
\[
z=u-a_eF_s=u-(e_h\cdot x'/|x'|^2)\,F' .
\]

\begin{lemma}[Regularity of the continuum derivative]\label{lem:fd-reg}
Let $e_h$ be piecewise linear, $u=D\Ft(x)e_h$, and $0<h\le\tau\le1$. Then, $z\in H^2$, and
\begin{gather}
\|u\|_{H^1}\le C\|e_h\|_{H^1},\qquad\|z\|_{H^2}\le C\tau^{-1}\|e_h\|_{H^1},\label{eq:fd-z}\\
\|u-\Ih u\|_{H^1}\le Ch\tau^{-1}\|e_h\|_{H^1}.\label{eq:fd-interp}
\end{gather}
\end{lemma}

\begin{proof}
The bound for $u$ in \eqref{eq:fd-z} follows from \cref{prop:stab} and \cref{lem:p52-adnorm}. For $z$, we follow the proof of \cite[Lemma~5.3]{Rao2026}. As there, the subtraction of $a_eF_s$ cancels the stiffness variation, and we have $\ax(u,\varphi)+\ax'[e_h](F,\varphi)=\ax(z,\varphi)+\int R_e\cdot\partial_s\varphi\,\dd s$ with $R_e=-\kc(e_h\cdot\bn)F_s+(e_h\cdot\bt)\,\partial_s^2F$, which is a smooth coefficient times $e_h$, so that its derivative lies in $L^2$. In Rao's equation for $z$, the display before \cite[(5.23)]{Rao2026}, the remaining terms are his local mass $\tau^{-1}\bx$ and its variation $b_x'[e_h](v,\varphi)$. We replace them by the remaining terms of \eqref{eq:p52-lin}. Each of them is $\int G\cdot\varphi\,\dd\xi$ with $\|G\|_{L^2}\le C\tau^{-1}\|e_h\|_{H^1}$, because the densities $\SL(u_n-e_n)$, $f_{\mathcal L}[e_h]\varsigma_v$, $\mathcal L_x(v\cdot\Jm e_h')$ and $(\SL v_n)\Jm e_h'$ are bounded in $L^2$ by $C\|e_h\|_{H^1}$ (\cref{lem:p52-dSL}\ref{item:p52-dSL-a}), and because the test $\varphi=\bn$ bounds the multiplier of the $\SL$-form by $C\tau^{-1}\|e_h\|_{H^1}$. Hence, $\ax(z,\varphi)=\int\tilde G\cdot\varphi\,\dd\xi$ with $\|\tilde G\|_{L^2}\le C\tau^{-1}\|e_h\|_{H^1}$, so $z'/|x'|$ has the weak derivative $-\tilde G$, and the bound for $z$ in \eqref{eq:fd-z} follows. For \eqref{eq:fd-interp}, we write $u=z+Be_h$ with $B=F'\otimes x'/|x'|^2\in W^{2,\infty}$, and we use the superapproximation $\|Be_h-\Ih(Be_h)\|_{H^1}\le Ch\|e_h\|_{H^1}$ of \cite[proof of Lemma~5.3]{Rao2026} and $\|z-\Ih z\|_{H^1}\le Ch\|z\|_{H^2}$.
\end{proof}

\subsubsection*{The three cancellations}
The residual of $\Ih u$ in \eqref{eq:fd-lin-h} contains terms of order one, namely the constant parts of the two data of the mass, the frame terms in both slots and the multiplier terms. They are not small one by one, and the next lemma removes them by three exact identities. The means of the data are removed by the linearized area constraint, and the multiplier terms and the frame terms in the first slot, which carry the same test factor $\Om\varphi$ or $\varphi\cdot\bn\,\dd s$, reduce together to the curvature tested against the variation of this factor. Let $\tilde u=u-e_h$, so that $\Ih u-e_h=\Ih\tilde u$, and let
\[
g=g_{e_h}=v\cdot D\bn[e_h]+v_n\,\bt\cdot\partial_se_h,\qquad
T[e_h](\varphi)=2\int\Lp\big(v_n\,\bt\cdot\partial_se_h\big)(\varphi\cdot\bn)\,\dd s .
\]
Let $m_u=\avg{\tilde u_n}$ and $m_g=\avg{g}$, let $\tilde u_0=P_0\tilde u_n$ and $\mathfrak u=\Lp\tilde u_0$, let $\hat y$ and $\hat y_1$ be the consistent loads, with $\rho=|x'|$, of $g$ and of $1$, and let $\hat y_0=\hat y-m_g\hat y_1$, the consistent load of $g-m_g$. Finally, let $\hat f=D\Lam[e_h]\kt+(\Lam_\Gamma\kt)\,\bt\cdot\partial_se_h$ be the function of \cref{lem:fd-shape}\ref{item:fd-shape-load} with $\mathfrak g=\kt$.

\begin{lemma}[The three cancellations]\label{lem:fd-cancel}
Let $d_h$, $g_h$, $c'$ and $\mathsf y_h$ be as in \eqref{eq:fd-LOmv} and the paragraph before it. The following identities hold.
\begin{enumerate}[label=(\roman*)]
\item\label{item:fd-cancel-mean} The means satisfy $m_u+\tau m_g=0$. Consequently, the discrete datum $\Phi=\Om\Ih\tilde u+\Om(e_h)(F_h-X)$ of the mass decomposes exactly as
\[
\Phi=M_v\Ih\tilde u_0+\tau\hat y_0+R,\qquad R=\rho_1+\tau\rho_2+\tau m_g(\hat y_1-M_v\one)+\Om(e_h)E,
\]
with $\rho_1=\Om\Ih\tilde u-M_v\Ih\tilde u_n$ and $\rho_2=\Om(e_h)\Ih v-\hat y$. Here, $\one^\T\hat y_0=0$.
\item\label{item:fd-cancel-T} With the loads above, we have
\[
\begin{aligned}
&2\int(\Lp g)(\varphi\cdot\bn)\,\dd s+\int(\Lp\hat f)(\varphi\cdot\bn)\,\dd s\\
&\qquad=2\int\Lp\big(v\cdot D\bn[e_h]\big)(\varphi\cdot\bn)\,\dd s-2\int\big(\Lp D\Lam[e_h]\Lp v_n\big)(\varphi\cdot\bn)\,\dd s .
\end{aligned}
\]
Moreover, $\one^\T\hat y^f=0$ for the consistent load $\hat y^f$ of $\hat f$, and on the discrete side, we have
\[
-2(\Om\varphi)^\T\Lhp D\Lh[e_h]\Lhp\Om\Ih v=(\Om\varphi)^\T\Lhp D\Lh[e_h]\Ih\kt-2(\Om\varphi)^\T\Lhp D\Lh[e_h]\,\mathsf y_h,
\]
with $\mathsf y_h$ as in \eqref{eq:fd-LOmv}.
\item\label{item:fd-cancel-mult} The multiplier term and the first-slot frame term combine, on the discrete and on the continuum side, to
\begin{gather}
c_h\one^\T\Om(e_h)\varphi-\tfrac2\tau\big(\Om(e_h)\varphi\big)^\T\Lhp\Om(F_h-X)=\big(\Om(e_h)\varphi\big)^\T\kh ,\label{eq:mult-disc}\\
c\int\varphi\cdot\Jm e_h'\,\dd\xi-2\int(\Lp v_n)\,\varphi\cdot\Jm e_h'\,\dd\xi=\int\kt\,\varphi\cdot\Jm e_h'\,\dd\xi .\label{eq:mult-cont}
\end{gather}
Moreover, $\kh=\Ih\kt+\theta_h$ with
\[
\theta_h=\varpi_h\one-\Lhp(d_h+2g_h)-\tfrac2\tau\Lhp\Om E,\qquad\varpi_h=c_h-c' .
\]
\end{enumerate}
\end{lemma}

\begin{proof}
For \ref{item:fd-cancel-mean}, $\int v\cdot\Jm e_h'\,\dd\xi=\int g\,\dd s$, so the linearized constraint of \cref{lem:fd-lin} is $m_u+\tau m_g=0$. With $F_h-X=\tau\Ih v+E$, we have $\Om\Ih\tilde u=M_v\Ih\tilde u_0+m_uM_v\one+\rho_1$ and $\tau\Om(e_h)\Ih v=\tau\hat y_0+\tau m_gM_v\one+\tau m_g(\hat y_1-M_v\one)+\tau\rho_2$, and the two multiples of $M_v\one$ cancel. Moreover, $\one^\T\hat y=\int g\,\dd s$ and $\one^\T\hat y_1=\int\dd s$.

For \ref{item:fd-cancel-T}, the definition of $D\bn[e_h]$ shows that $2\int(\Lp g)(\varphi\cdot\bn)\,\dd s$ is the second-slot frame term plus $T[e_h](\varphi)$. The relations $\kt=c-2\Lp v_n$ and $\int v_n\,\dd s=0$ give $\Lam_\Gamma\kt=-2v_n$, and \cref{lem:fd-shape}\ref{item:fd-shape-lin} gives $-2\Lp D\Lam[e_h]\Lp v_n=\Lp D\Lam[e_h]\kt=\Lp\hat f+2\Lp(v_n\,\bt\cdot\partial_se_h)$, so $\int(\Lp\hat f)(\varphi\cdot\bn)\,\dd s$ is the DtN term minus $T[e_h](\varphi)$. Adding the two gives the identity. The mean of $\hat f$ vanishes by \cref{lem:fd-shape}\ref{item:fd-shape-load}, and the discrete identity follows from \eqref{eq:fd-LOmv} and $D\Lh[e_h]\one=0$.

For \ref{item:fd-cancel-mult}, $D\Om[e_h]=\Om(e_h)$, and \eqref{eq:mult-disc} is the formula for $\kh$ in \cref{sec:scheme}. The relation $2\Lp v_n=c-\kt$ gives \eqref{eq:mult-cont}. The same formula for $\kh$, $F_h-X=\tau\Ih v+E$ and \eqref{eq:fd-LOmv} give $\theta_h$.
\end{proof}

\begin{proof}[Proof of \cref{lem:fd-comp}]
By \cref{lem:fd-diff-h}, $DF_{h,\tau}(X)e_h=u_h$ is the solution of \eqref{eq:fd-lin-h}. We bound the residual of $\Ih u$ in \eqref{eq:fd-lin-h} against P1 tests $\varphi$ on the discrete constraint kernel, by subtracting \eqref{eq:fd-lin} for the same $\varphi$. We split $F_h=\Ih F+E$ in the stiffness variation and $F_h-X=\tau\Ih v+E$ in the operator derivative. By \cref{lem:fd-cancel}, in which the terms $\pm T[e_h]$ cancel, and since $\delta c_h\one^\T\Om\varphi$ vanishes on the kernel, the residual is
\[
\begin{aligned}
&\underbrace{a_X(\Ih u,\varphi)-\ax(u,\varphi)+a_X'[e_h](\Ih F,\varphi)-\ax'[e_h](F,\varphi)}_{\text{stiffness part}}\\
&+\underbrace{2\tau^{-1}(\Om\varphi)^\T\Lhp\Phi-2\tau^{-1}\int \mathfrak u\,\varphi\cdot\bn\,\dd s-2\int(\Lp g)(\varphi\cdot\bn)\,\dd s}_{\text{velocity mass part}}\\
&+\underbrace{\Big[(\Om\varphi)^\T\Lhp D\Lh[e_h]\Ih\kt-2(\Om\varphi)^\T\Lhp D\Lh[e_h]\,\mathsf y_h-\int(\Lp\hat f)(\varphi\cdot\bn)\,\dd s\Big]}_{\text{operator derivative part}}\\
&+\underbrace{a_X'[e_h](E,\varphi)-2\tau^{-1}(\Om\varphi)^\T\Lhp D\Lh[e_h]\Lhp\Om E}_{\text{replacement part}}\\
&\underbrace{-\big(\Om(e_h)\varphi\big)^\T\kh+\int\kt\,\varphi\cdot\Jm e_h'\,\dd\xi+(\delta c+c_P)\int\varphi\cdot\bn\,\dd s}_{\text{multiplier part}} ,
\end{aligned}
\]
with $\mathfrak u=\Lp\tilde u_0$ as above. We bound the five parts in units of $\|e_h\|_{H^1}\|\varphi\|_{H^1}$, and the constraint residual in units of $\|e_h\|_{H^1}$. The terms $\tilde u_0$, $\hat y_0$ and $R$ are those of the decomposition of $\Phi$ in \cref{lem:fd-cancel}\ref{item:fd-cancel-mean}.

The stiffness part is Rao's. By \eqref{eq:fd-rao-b}, \eqref{eq:fd-interp} and $\|\Ih u\|_{H^1}\le C\|e_h\|_{H^1}$, which follows from \eqref{eq:fd-z}, \eqref{eq:fd-interp} and $h\le\tau$, it is bounded by $C(h+h\tau^{-1})$.

For the velocity mass part, we use the decomposition of $\Phi$ in \cref{lem:fd-cancel}\ref{item:fd-cancel-mean}. For $M_v\Ih\tilde u_0$, we have $\Lam_\Gamma \mathfrak u=\tilde u_0$ by \cref{lem:cont-gauge}(a), and $\|\mathfrak u\|_{H^2}\le C\|e_h\|_{H^1}$ by \eqref{eq:cont-SLsplit} and \cref{lem:cont-kernelreg}(b) in the normalized arclength of $x$, transported to the parameter of $x$ by \cref{lem:semi-material}, and \cref{lem:in-bweak} with the datum $\mathfrak u$ gives $\Lhp M_v\Ih\tilde u_0=P\Ih \mathfrak u-\Lhp\mathsf r$ with $\dnorm{\mathsf r}\le Ch\|e_h\|_{H^1}$. The constant removed by $P$ does not contribute by \ref{item:fd-conv-kernel}, the first bound of \ref{item:fd-conv-p1} for the P1 function $\Ih \mathfrak u$ and $\|\mathfrak u-\Ih \mathfrak u\|_{L^2}\le Ch^2\|\mathfrak u\|_{H^2}$ compare $(\Om\varphi)^\T\Ih \mathfrak u$ with $\int \mathfrak u\,\varphi\cdot\bn\,\dd s$ up to $Ch\|e_h\|_{H^1}\|\varphi\|_{H^1}$, and the first bound of \ref{item:fd-conv-load} bounds the contribution of $\mathsf r$. Hence, this part differs from $2\tau^{-1}\int \mathfrak u\,\varphi\cdot\bn\,\dd s$ by $C\tau^{-1}h$, which is the loss $h\tau^{-1}$ of Rao's mass. For $\tau\hat y_0$, \ref{item:fd-conv-inv} with $\Lp(g-m_g)=\Lp g$ and $\|g\|_{L^2}\le C\|e_h\|_{H^1}$ gives the error $Ch^{1/2}$. For $R$, the loads $\rho_1$, $\hat y_1-M_v\one$ and $\Om(e_h)E$ are lumped loads of P1 functions of sizes $O(h^2)|\tilde u|$, $O(h^2)$ and $Ch^2\|e_h'\|_{L^2}$, by \eqref{eq:fd-mu} and \cref{lem:fd-cons}, and $\dnorm{\rho_2}\le Ch^{3/2}\|e_h\|_{H^1}$ follows from a summation by parts in the entries $\tfrac12\Jm(e_{h,j+1}-e_{h,j-1})\cdot v(\xi_j)$ and $\int\varphi_j\,v\cdot\Jm e_h'\,\dd\xi$ and an inverse estimate. Hence, $\tau^{-1}\dnorm{R}\le C(h^2\tau^{-1}+h^{3/2})\|e_h\|_{H^1}$, and \ref{item:fd-conv-load} bounds the velocity mass part by $C(h\tau^{-1}+h^{1/2})$.

For the operator derivative part, \cite[Theorem~3.7]{EtoPaperI} (\cref{sec:input-results}) with the datum $\kt$ and \cref{lem:fd-shape}\ref{item:fd-shape-load} give $D\Lh[e_h]\Ih\kt=\hat y^f+\mathsf r$ with $\|\mathsf r\|_{(H^1)^*}\le Ch\|e_h'\|_{L^\infty}\le Ch^{1/2}\|e_h\|_{H^1}$ by the inverse estimate, and $\|\Lhp\Om\varphi\|_{H^1}\le C\|\varphi\|_{L^2}$ bounds the contribution of $\mathsf r$. Since $\one^\T\hat y^f=0$ and $\|\hat f\|_{L^2}\le C\|e_h\|_{H^1}$ by \cref{lem:fd-shape}\ref{item:fd-shape-bd}, \ref{item:fd-conv-inv} compares $(\Om\varphi)^\T\Lhp\hat y^f$ with $\int(\Lp\hat f)(\varphi\cdot\bn)\,\dd s$ up to $Ch^{1/2}$. The term with $\mathsf y_h$ is bounded by \cite[Lemma~A.21]{EtoPaperI} and an inverse estimate,
\[
\begin{aligned}
2|(\Lhp\Om\varphi)^\T D\Lh[e_h]\,\mathsf y_h|&\le C\|e_h'\|_{L^2}\|\varphi\|_{L^2}\|\mathsf y_h\|_{H^1}\\
&\le C\|e_h\|_{H^1}\|\varphi\|_{L^2}h^{-1/2}\|\mathsf y_h\|_{H^{1/2}}\le Ch\|e_h\|_{H^1}\|\varphi\|_{L^2},
\end{aligned}
\]
by \eqref{eq:fd-LOmv}. Hence, the operator derivative part is bounded by $Ch^{1/2}$.

For the replacement part, $\|E'\|_{L^\infty}\le Ch^{-1/2}\|E\|_{H^1}\le Ch^{3/2}$ and \eqref{eq:fd-rao-b} bound the stiffness term by $Ch^{3/2}$, and \cite[Lemma~A.21]{EtoPaperI} with $\|\Lhp\Om E\|_{H^1}\le Ch^2$ (\cite[Lemma~A.20]{EtoPaperI}) bounds the mass term by $Ch^2\tau^{-1}$.

For the multiplier part, $(\Om(e_h)\varphi)^\T\Ih\kt$ is the trapezoidal rule for $\int\Ih\kt\,\varphi\cdot\Jm e_h'\,\dd\xi$, so its difference from $\int\kt\,\varphi\cdot\Jm e_h'\,\dd\xi$ is bounded by $Ch$. The part $\theta_h$ contributes $C\|\theta_h\|_{L^2}\|\varphi\|_{L^\infty}\|e_h'\|_{L^2}$. By \eqref{eq:in-Lhp}, $\|\Lhp(d_h+2g_h)\|_{H^{1/2}}\le C(\dnorm{d_h}+\dnorm{g_h})\le Ch^{3/2}$, where the power $3/2$ comes from $\dnorm{d_h}$, and $\tfrac2\tau\|\Lhp\Om E\|_{H^{1/2}}\le Ch^2\tau^{-1}$. For the constant $\varpi_h$, we test the discrete step in the form $\aX(F_h,\varphi)=(\Om\varphi)^\T\kh$ with $\varphi_*=\Ih\bn_x$ of \cref{cor:bem-acoerc}, for which $\one^\T\Om\varphi_*\ge c>0$, and compare with $\ax(F,\varphi_*)=\int\kt\,\varphi_*\cdot\bn\,\dd s$ through \eqref{eq:fd-rao-a}, the vertex rule \ref{item:fd-conv-p1} and $|\aX(E,\varphi_*)|\le Ch^2$. This gives $\varpi_h\,\one^\T\Om\varphi_*=(\Om\varphi_*)^\T\Lhp(d_h+2g_h)+\tfrac2\tau(\Om\varphi_*)^\T\Lhp\Om E+O(h^2)$, and hence
\begin{equation}\label{eq:fd-gamma}
|c_h-c'|=|\varpi_h|\le C(h^{3/2}+h^2\tau^{-1}),\qquad\|\theta_h\|_{L^2}\le C(h^{3/2}+h^2\tau^{-1}).
\end{equation}
Finally, $|\delta c+c_P|\le C\|e_h\|_{H^1}$ by \cref{lem:fd-lin}, and \ref{item:fd-conv-kernel} bounds $(\delta c+c_P)\int\varphi\cdot\bn\,\dd s$ by $Ch^2$. Therefore, the multiplier part is bounded by $C(h+h^2\tau^{-1})$.

Since $h^{3/2}\le h\le h^{1/2}$ and $h^2\tau^{-1}\le h\tau^{-1}$, the residual is bounded by $C(h\tau^{-1}+h^{1/2})\|e_h\|_{H^1}\|\varphi\|_{H^1}$. The constraint residual of $\Ih u$ is $\one^\T\Phi=\one^\T M_v\Ih\tilde u_0+\one^\T R$, in which $\one^\T M_v\Ih\tilde u_0$ is the vertex rule for $\int\tilde u_0\,\dd s=0$ and $|\one^\T R|\le\dnorm{R}\le C(h^2+\tau h^{3/2})\|e_h\|_{H^1}$ by the bound of $R$ above, so $|\one^\T\Phi|\le Ch\|e_h\|_{H^1}$. The stability bound \eqref{eq:use-brezzi} with $\mathsf R=C(h\tau^{-1}+h^{1/2})\|e_h\|_{H^1}$ and $\delta=Ch\|e_h\|_{H^1}$, in which the constraint residual enters as $\tau^{-1}\delta$, gives the claim.
\end{proof}

\subsubsection*{The nonlinear remainder}
\begin{proof}[Proof of \cref{lem:fd-rem}]
Let $u_h$ solve \eqref{eq:fd-lin-h} at $X$, let $X_1=X+e_h$ and $R_h=F_{h,\tau}(X_1)-F_{h,\tau}(X)-u_h$. Let $r_2\in(0,r_0]$, depending only on the constants of $\mathscr K$, be so small that the polygons of $\mathcal N_{r_2}$ have chord-arc constant at least $c_*/2$ and the results of \cref{sec:bem} used below hold on $\mathcal N_{r_2}$, and let $[X,X_1]\subset\mathcal N_{r_2}$. As in \cite[(6.3)]{Rao2026}, subtracting the two step equations and the linearized system gives, for every P1 $\varphi$ with $\one^\T\Om(X_1)\varphi=0$,
\[
\big(a_{X_1}+\tau^{-1}b_{X_1,h}\big)(R_h,\varphi)=-\delta_1a(u_h,\varphi)-\tau^{-1}\delta_1b(u_h-e_h,\varphi)-\delta_2a(F_h,\varphi)-\tau^{-1}\delta_2b(F_h-X,\varphi)+\mathcal M(\varphi),
\]
where $\delta_1$ and $\delta_2$ are the first- and second-order increments in the geometry from $X$ to $X_1$, and $\mathcal M$ collects the multiplier terms, which Rao's scheme does not have. The constraint residual of $R_h$ is $-\one^\T\Om(e_h)(u_h-e_h)$. The bounds below use $u_h$ only through $\|u_h\|_{H^1}+\|u_h-e_h\|_{H^1}$, and they use neither \cref{lem:fd-diff-h} nor the conclusion of \cref{lem:fd-comp}. The bound \eqref{eq:fd-gamma} of $c_h$, taken from the proof of \cref{lem:fd-comp}, concerns only the step $F_h$ and does not use \cref{lem:fd-diff-h}. For fixed $N$, \eqref{eq:fd-lin-h} is a uniquely solvable linear system in finite dimension, so $\|u_h\|_{H^1}\le C_N\|e_h\|_{H^1}$, and the bounds below give $\|R_h\|_{H^1}\le C_N'\|e_h\|_{H^1}^2$, which is the differentiability of \cref{lem:fd-diff-h}. With \cref{lem:fd-diff-h}, the estimate \eqref{eq:fd-stab}, which rests on \cref{prop:stab}, \cref{lem:fd-comp}, \eqref{eq:fd-interp} and \cref{lem:p52-adnorm} and not on the present lemma, and \cref{lem:p52-adnorm} give $\|u_h\|_{H^1}+\|u_h-e_h\|_{H^1}\le C\|e_h\|_{H^1}$, and the bounds below are applied again with this uniform bound.

As in the proof of \cite[Lemma~6.1]{Rao2026}, with $\|F_h'\|_{L^\infty}\le C$ by \cref{lem:fd-cons}, the stiffness increments are bounded by $C\|e_h'\|_{L^\infty}\|e_h\|_{H^1}\|\varphi\|_{H^1}$, hence by $Ch^{-1/2}\|e_h\|_{H^1}^2\|\varphi\|_{H^1}$. For the mass increments, the boundary element bounds of \cref{prop:in-import} replace Rao's quadrature bounds. The increment bound of \cref{thm:in-c1p}, which includes the terms with $D\Om$, bounds the first-order mass increment by $C\tau^{-1}\|e_h\|_{H^1}^2\|\varphi\|_{H^1}$, and \cref{thm:in-c3p} with $f=(F_h-X)/\tau$, which is bounded in $L^\infty$ by \cref{lem:fd-cons}, bounds the second-order one by $C\|e_h'\|_{L^\infty}\|e_h'\|_{L^2}\|\varphi\|_{L^2}\le Ch^{-1/2}\|e_h\|_{H^1}^2\|\varphi\|_{L^2}$.

The constraint functional is linear in $X$, so only its first-order increment $\one^\T\Om(e_h)\varphi$ appears in $\mathcal M$, multiplied by the multiplier increment. Testing \eqref{eq:fd-lin-h} with $\varphi_*=\Ih\bn_x$, for which $\one^\T\Om\varphi_*\ge c$, bounds this increment by $C\tau^{-1}\|e_h\|_{H^1}$. In this test, $|c_h|\le C$ by \eqref{eq:fd-gamma}, and $|D_X\bXh[e_h](F_h-X,\varphi_*)|\le C\tau\|e_h\|_{H^1}$ is the limit as $t\to0$ of the increment bound of \cref{thm:in-c1p} along $te_h$, divided by $t$. Hence, $|\mathcal M(\varphi)|\le C\tau^{-1}\|e_h\|_{H^1}^2\|\varphi\|_{H^1}$, and the constraint residual of $R_h$ is bounded by $C\|e_h\|_{H^1}^2$. The stability bound \eqref{eq:use-brezzi} at $X_1$, with $\mathsf R=C(h^{-1/2}+\tau^{-1})\|e_h\|_{H^1}^2$ and $\delta=C\|e_h\|_{H^1}^2$, gives the claim.
\end{proof}

\subsection{One step of the scheme}\label{app:bem-step}
The proofs of this subsection use the Schur form of $\Lh$ from \cite{EtoPaperI}. For $X\in\mathfrak C$, let $Q^X_{\rm int}=\tfrac12+\mathcal K_X$ and $Q^X_{\rm ext}=\tfrac12-\mathcal K_X$, which are operators on functions of the parameter and are not the matrices $Q_{\rm int}$ and $Q_{\rm ext}$ of \cref{sec:scheme}. For $g\in H^{1/2}$, the Galerkin inverse $\mathcal L_{X,h}^{-1}g\in\mathcal P_0$ is defined by $\langle\mathcal L_X\mathcal L_{X,h}^{-1}g,\chi\rangle=\langle g,\chi\rangle$ for all $\chi\in\mathcal P_0$, and the bordered Galerkin inverse $\mathcal L_{X,h,0}^{-1}g$ is the density $\varsigma\in\mathcal P_0$ of mean zero for which there is $\nu\in\R$ with $\langle\mathcal L_X\varsigma,\chi\rangle=\langle g,\chi\rangle+\nu\int\chi\,\dd\xi$ for all $\chi\in\mathcal P_0$. The discrete Schur form is
\[
\mathfrak s_{X,h}(v,w)=2\langle\mathcal L_X\partial_\xi v,\partial_\xi w\rangle+\langle\mathcal L_{X,h}^{-1}Q^X_{\rm int}v,Q^X_{\rm int}w\rangle+\langle\mathcal L_{X,h,0}^{-1}Q^X_{\rm ext}v,Q^X_{\rm ext}w\rangle
\]
for Lipschitz $v$ and $w$ \cite[(A.15)]{EtoPaperI}, and the continuum form $\langle\mathcal S_Xv,w\rangle$ is the same expression with the exact inverses $\mathcal L_X^{-1}$ and $\mathcal L_{X,0}^{-1}$, the latter being the inverse with the gauge of \cite[Lemma~A.10]{EtoPaperI} recalled in \cref{sec:input-results}. The form $\mathfrak s_{X,h}$ satisfies $|\mathfrak s_{X,h}(v,w)|\le C_{\mathfrak s}\|v\|_{H^{1/2}}\|w\|_{H^{1/2}}$ for Lipschitz $v$ and $w$ \cite[Lemma~A.15]{EtoPaperI}. If $X$ is a polygon, the double layer $\mathcal K_X$ is bounded on $L^2$, $H^{1/2}$ and $H^1$ \cite[Lemma~A.8]{EtoPaperI}(c), and the upper bound extends to $v\in H^1$ by continuity. For a polygon $X\in\mathfrak C$ and P1 functions $v$ and $w$, we have $\mathfrak s_{X,h}(v,w)=w^\T\Lh(X)v$ and $\mathfrak s_{X,h}(\one,w)=\mathfrak s_{X,h}(w,\one)=0$, and $\mathfrak s_{X,h}(w,w)\ge c_{\mathfrak s}\|w\|_{H^{1/2}}^2$ if the vertex values of $w$ sum to zero. By \cite[(A.17)]{EtoPaperI}, we have $|\mathfrak s_{X,h}(u,w)-\langle\mathcal S_Xu,w\rangle|\le Ch^{1/2}\|u\|_{H^1}\|w\|_{H^{1/2}}$ for a polygon $X\in\mathfrak C$, $u\in H^1$ and Lipschitz $w$, and also for $X\in\mathfrak C$ and Lipschitz $u$ and $w$. For a polygon $X\in\mathfrak C$, $g\in H^1$ with gauge pair $(\varsigma,c)$, and P1 functions $w$, we have $\langle\mathcal S_Xg,w\rangle=\int_\Tor\varsigma\,w\,\dd\xi$ by \cref{prop:ext-layer} for polygons. For a P1 function $w$, we have $\sum_jw_j=N\int_\Tor w\,\dd\xi$, so the vertex values sum to zero if and only if $w$ has mean zero.

\begin{proof}[Proof of \cref{cor:bem-acoerc}]
Let $\Gamma=x(\Tor)$ and $\bn_x=\Jm x'/|x'|$. By \ref{K1}, $\ell(x)\ge c_*$ and $c_*^2/M_*\le c_*$. We first show that the normal components of translations satisfy
\begin{equation}\label{eq:use-c0}
c_{\rm tr}:=\min_{|u|=1}\int_\Gamma(u\cdot\bn_x)^2\,\dd s\ge\frac{c_*^2}{M_*}.
\end{equation}
Let $|u|=1$ and $v=\Jm^\T u$, so that $u\cdot\bn_x\,\dd s=v\cdot x'\,\dd\xi$. The chord $D(\xi)=x(\xi+\tfrac12)-x(\xi)$ satisfies $|D|\ge c_*/2$ and $D(\xi+\tfrac12)=-D(\xi)$, so its direction passes through $v$ or $-v$, and there $|v\cdot D|\ge c_*/2$. Hence, $v\cdot x$ has oscillation at least $c_*/2$, its total variation over the closed curve is at least $c_*$, and $\int_\Gamma|u\cdot\bn_x|\,\dd s\ge c_*$. The Cauchy--Schwarz inequality and $\ell(x)\le M_*$ give \eqref{eq:use-c0}.

Let $u$ be piecewise linear with $\one^\T\Om(X)u=0$, and write $u=\bar u+u_0$ with $\bar u=\int_\Tor u\,\dd\xi$. Since $|\Gamma_p|\le M_1h$ on each panel, $\aX(u,u)\ge M_1^{-1}|u|_{H^1}^2$, and $\|u_0\|_{H^1}\le C|u|_{H^1}$ by the Poincar\'e inequality. By \eqref{eq:in-dual}, $\bXh(u,u)=2(\Om u)^\T\Lhp(\Om u)\ge c\,\dnorm{\Om u}^2$, and by \eqref{eq:use-Om}, $\dnorm{\Om u_0}\le C_\Omega|u|_{H^1}$. For the translation, let $w=\Ih(\bar u\cdot\bn_x)$, so that $\|w\|_{H^{1/2}}\le C|\bar u|$. By \eqref{eq:use-area} for the field $w\bar u$,
\[
w^\T\Om(X)\bar u=\int_\Tor w\,\bar u\cdot\Jm X'\,\dd\xi\ge\int_\Gamma(\bar u\cdot\bn_x)^2\,\dd s-\underbrace{C(h+r_1)|\bar u|^2}_{\text{replacement}},
\]
where the replacement term bounds the cost of replacing $w$ by $\bar u\cdot\bn_x$ and $X'$ by $x'$. With $C_\sharp$ large enough, \eqref{eq:use-small} and \eqref{eq:use-c0} give $w^\T\Om(X)\bar u\ge\tfrac12c_{\rm tr}|\bar u|^2$, hence $\dnorm{\Om\bar u}\ge c_\Omega|\bar u|$. If $c_\Omega|\bar u|\ge2C_\Omega|u|_{H^1}$, then $\dnorm{\Om u}\ge\tfrac12c_\Omega|\bar u|$, so $\bXh(u,u)\ge c|\bar u|^2$, and we use $\tau^{-1}\ge1$. Otherwise, $|\bar u|\le C|u|_{H^1}$, and the stiffness bound alone suffices, because $\bXh(u,u)\ge0$. In both cases, the left side is bounded below by $c(|u|_{H^1}^2+|\bar u|^2)\ge c\|u\|_{H^1}^2$.

For the constraint, $\varphi_*=\Ih\bn_x$ satisfies $\|\varphi_*\|_{H^1}\le C$, and by \eqref{eq:use-area}, $\one^\T\Om(X)\varphi_*$ differs from $\int\bn_x\cdot\Jm x'\,\dd\xi=\ell(x)\ge c_*$ by an amount bounded by $C(r_1+h)$, which is not larger than $c_*/2$ by \eqref{eq:use-small}.

For \eqref{eq:use-brezzi}, let $B_\tau(u,\varphi)=\aX(u,\varphi)+\tau^{-1}\bXh(u,\varphi)$. Since $|\Gamma_p|\ge c_{\rm ca}h$, \cref{lem:in-loads} gives $|B_\tau(u,\varphi)|\le C\tau^{-1}\|u\|_{H^1}\|\varphi\|_{H^1}$. Let $\beta=\one^\T\Om(X)U/\one^\T\Om(X)\varphi_*$, so that $|\beta|\le2\delta/c_*$, and $U_0=U-\beta\varphi_*$, so that $\one^\T\Om(X)U_0=0$. The coercivity and the assumption with $\varphi=U_0$ give
\[
c\|U_0\|_{H^1}^2\le B_\tau(U_0,U_0)=B_\tau(U,U_0)-\beta B_\tau(\varphi_*,U_0)\le\big(\mathsf R+C\tau^{-1}|\beta|\big)\|U_0\|_{H^1},
\]
and \eqref{eq:use-brezzi} follows from $\|U\|_{H^1}\le\|U_0\|_{H^1}+|\beta|\,\|\varphi_*\|_{H^1}$. Finally, \eqref{eq:step-var} with $X^m=X$ is a square linear system in $U=X^{m+1}-X$ and $c_h$, namely $B_\tau(U,\varphi)-c_h\one^\T\Om(X)\varphi=-\aX(X,\varphi)$ for all piecewise linear $\varphi$ and $\one^\T\Om(X)U=0$. For a solution of the homogeneous system, the test $\varphi=U$ gives $B_\tau(U,U)=0$, so the coercivity gives $U=0$, and the test $\varphi=\varphi_*$ gives $c_h=0$. Hence, the system is uniquely solvable.
\end{proof}

\begin{proof}[Proof of \cref{lem:bem-invcons}]
Let first $y$ be consistent with $\rho=|X'|$. Since $\one^\T y=0$, we have $\int_\Tor|X'|f\,\dd\xi=0$, so $u_X=\mathcal L_X(|X'|f)+c_u$, with $c_u$ chosen so that $u_X$ has mean zero, satisfies $\|u_X\|_{H^1}\le C\|f\|_{L^2}$ \cite[Lemma~A.8]{EtoPaperI}. Let $u_h=\Lhp y$, which has mean zero. The pair $(|X'|f,-c_u)$ is the gauge pair of $u_X$, so $\langle\mathcal S_Xu_X,w\rangle=\int_\Tor|X'|fw\,\dd\xi=w^\T y=\mathfrak s_{X,h}(u_h,w)$ for P1 functions $w$, where the last equality uses $\Lh(X)u_h=y$. Thus, $u_h$ is the Galerkin approximation of $u_X$ with the perturbed form $\mathfrak s_{X,h}$ on the P1 functions of mean zero, and we use the argument of Strang's lemma \cite[Thm.~8.3]{Steinbach2008}, with the difference that the consistency error is evaluated at $u_X$ and not at the unperturbed Galerkin solution. The coercivity of $\mathfrak s_{X,h}$, the approximation property $\inf_{v_h}\|u_X-v_h\|_{H^{1/2}}\le Ch^{1/2}\|u_X\|_{H^1}$ of \cite[Thm.~10.9]{Steinbach2008} and the consistency \cite[(A.17)]{EtoPaperI} give
\begin{equation}\label{eq:invcons-strang}
\|u_h-u_X\|_{H^{1/2}}\le Ch^{1/2}\|f\|_{L^2}.
\end{equation}
By the definition in \cref{sec:scheme}, $\Lp f=P_0^x\mathcal L_x(|x'|P_0^xf)$, where $P_0^x$ removes the arclength mean on the curve. Since $\int f\,\dd s_X=0$ and $\|X'-x'\|_{L^\infty}\le C_Ih$, we have $\||X'|f-|x'|P_0^xf\|_{L^2}\le Ch\|f\|_{L^2}$, and the outer projection and $c_u$ change only constants. Hence, \cite[Lemma~A.7]{EtoPaperI} and \cite[Lemma~A.8]{EtoPaperI} give a constant $\bar c$ with
\begin{equation}\label{eq:invcons-curve}
\|u_X-\Lp f-\bar c\|_{H^1}\le\|(\mathcal L_X-\mathcal L_x)(|X'|f)\|_{H^1}+\|\mathcal L_x(|X'|f-|x'|P_0^xf)\|_{H^1}\le Ch\|f\|_{L^2}.
\end{equation}
Adding \eqref{eq:invcons-strang} and \eqref{eq:invcons-curve} gives the claim, and $|\bar c|\le C\|f\|_{L^2}$ because $u_X$ and $\Lp f$ are bounded in $L^2$.

If $\rho=|x'|$, then $y$ is the consistent load with $\rho=|X'|$ of $\tilde f=f|x'|/|X'|$, where $\|\tilde f-f\|_{L^2}\le Ch\|f\|_{L^2}$, and $\Lp$ is bounded from $L^2$ to $H^1$ \cite[Lemma~A.8]{EtoPaperI}(b). If $y=M_vf$ is lumped, let $M_c$ be the consistent mass matrix, $(M_c)_{ij}=\int\varphi_i\varphi_j\,|X'|\,\dd\xi$. The trapezoidal rule is exact for linear functions on each panel, so $\one^\T M_cf=\one^\T M_vf=0$, and the consistent case applies to $M_cf$. On a panel, $M_v-M_c$ is $\tfrac{|\Gamma_p|}6\left(\begin{smallmatrix}1&-1\\-1&1\end{smallmatrix}\right)$, so the inverse estimates on $\mathcal P_1$ give
\[
|w^\T(M_v-M_c)f|\le Ch^2\|w'\|_{L^2}\|f'\|_{L^2}\le Ch^{1/2}\|w\|_{H^{1/2}}\|f\|_{L^2}.
\]
Hence, $\dnorm{(M_v-M_c)f}\le Ch^{1/2}\|f\|_{L^2}$, and \eqref{eq:in-Lhp} completes the proof.
\end{proof}

\section{Formally verified auxiliary results}\label{app:map}

This appendix collects the precise statements of auxiliary results used in the paper whose proofs are not printed here. They are proved in the Lean~4 formalization described in \cref{sec:formal}, from the five hypotheses of \cref{sec:formal-ext}. The same holds for the results of Sections~2 to~5 whose proofs are not printed in full, namely \cref{prop:closure,lem:semi-material,lem:Kfamily}, and for \cref{prop:in-import}, which is proved in \cite{EtoPaperI} and in its formalization. The proofs of \cref{lem:in-ellip,lem:in-loads} are given in \cref{sec:input-step}, and the proofs of \cref{thm:semi-B,cor:bem-acoerc,lem:bem-invcons,prop:stab,lem:fd-cons,lem:fd-comp,lem:fd-rem} and of \cref{lem:fd-shape,lem:fd-lin,lem:fd-reg} are given in \cref{app:new}.

\subsection{Standard tools}\label{map:tools}

\begin{lemma}[Embeddings]\label{lem:tool-sobolev}
Let $s>\frac12$, let $l\ge0$ be an integer, and let $k\in\R$ and $\theta\in[0,1]$. For $g\in H^s(\Tor)$ and $r\in\R$ with $s-r>\frac12$, we have $\sum_n\langle n\rangle^r|\hat g_n|\le C_{s-r}\|g\|_{H^s}$, and in particular $\wnorm{g}\le C_s\|g\|_{H^s}$. If $g$ is smooth, then
\[
\|g\|_{L^\infty}\le\wnorm{g},\qquad\|g^{(l)}\|_{L^\infty}\le C_l\|g\|_{H^{l+1}} .
\]
For $g\in H^{k+1}(\Tor)$, we have $\|g\|_{H^{k+\theta}}\le\|g\|_{H^k}^{1-\theta}\|g\|_{H^{k+1}}^{\theta}$. The constants depend only on the indicated indices.
\end{lemma}

\begin{lemma}[Kernel operators]\label{lem:tool-kernel}
For all real $a$ and $b$, with $a_+=\max(a,0)$ and $b_+=\max(b,0)$, we have
\begin{equation}\label{eq:tool-kernelop}
\|\mathcal I_k\|_{H^{-a}\to H^b}\le\|k\|_{H^{a_++b_+}(\Tor^2)} .
\end{equation}
For $a\ge0$ and $Z\in H^{a+1/2}$, we have $\|D[Z]\|_{H^a(\Tor^2)}\le2^a\|Z\|_{H^{a+1/2}}$.
\end{lemma}

A closed curve $x:\Tor\to\C$ of length $\ell$ in normalized arclength is written as $x'=\ell y$ with $|y|=1$. With the conventions of \cref{sec:ms-flow}, we have $y'=\iu\ell\kc y$. The signed area of a closed $C^1$ curve is $\Area(x)=\frac12\int_\Tor x\wedge x'\,\dd\xi$.

\begin{lemma}[Plane curves]\label{lem:tool-frame}
\begin{enumerate}
\item For a regular $C^2$ parametrization $w$ of the curve, the curvature is $\kc=\operatorname{Im}(\overline{w'}w'')/|w'|^3$.
\item Let $L,K\ge0$. There is $C$, depending only on $L$ and $K$, with the following property. Let $(x_j,\ell_j,\kc_j,y_j)$, $j=1,2$, satisfy $x_j'=\ell_jy_j$, $|y_j|=1$ and $y_j'=\iu\ell_j\kc_jy_j$ on $\Tor$, with closed $C^1$ curves $x_j$ and continuous $\kc_j$, $|\ell_2|\le L$ and $|\kc_1|\le K$. Then, for $0\le\xi\le1$,
\[
|y_1-y_2|+|x_1-x_2|+|x_1'-x_2'|\le C\big(|x_1(0)-x_2(0)|+|y_1(0)-y_2(0)|+|\ell_1-\ell_2|+\|\kc_1-\kc_2\|_{L^2}\big).
\]
\item Let $x$ and $z$ be closed $C^1$ curves, and let $u\in\C$ with $|u|=1$. If $|x'|\le M_1$, $|z|\le M_0$, $|ux'-z'|\le\delta$ and $\frac12\delta(M_1+M_0)<\Area(z)$, then $\Area(x)>0$.
\item For a smooth, regular and embedded closed curve of length $\ell$, we have $\int\kc^2\,\dd s\ge2/\ell$.
\end{enumerate}
\end{lemma}

\begin{lemma}[Resolvent multipliers]\label{lem:tool-mult}
Let $b>0$ and $0<\tau\le1$, and let $R_n=(1+\tau b|n|^3)^{-1}$.
\begin{enumerate}
\item For $0\le p\le3$, we have $\sup_n\tau\langle n\rangle^pR_n\le C\tau^{1-p/3}$, with $C$ depending only on $b$ and $p$.
\item For $0\le\theta\le1$, we have $0\le1-R_n\le(\tau b|n|^3)^\theta$.
\item There is $\beta>0$, depending only on $b$, with the following property. If $g,e\in H^1$ and $q\in L^2$ satisfy $(1+\tau b|n|^3)\hat g_n=\hat e_n+\tau\hat q_n$ for all $n$, then
\begin{equation}\label{eq:tool-energy}
\|g\|_{H^1}^2\le\|e\|_{H^1}^2+\tau\Big(2\beta\|g\|_{L^2}^2+\frac1{2\beta}\|q\|_{L^2}^2\Big).
\end{equation}
\end{enumerate}
\end{lemma}

\begin{lemma}[Discrete Gronwall inequalities]\label{lem:tool-gronwall}
Let $\tau>0$, $M\tau\le T$ and $E_m\ge0$ for $m\le M$.
\begin{enumerate}
\item If $C\ge0$ and $E_{m+1}\le(1+C\tau)E_m+C\tau$ for $m<M$, then $E_m\le\eu^{CT}(E_0+CT)$ for $m\le M$.
\item Let $0\le\beta<1$ and $K,a\ge0$. If $E_m\le a+K\sum_{j<m}\tau\big((m-j)\tau\big)^{-\beta}E_j$ for $m\le M$, then $E_m\le Ca$ for $m\le M$, with $C$ depending only on $\beta$, $K$ and $T$.
\end{enumerate}
\end{lemma}

Part (b) is, for $0<\beta<1$, the discrete Gronwall inequality with a weakly singular kernel \cite[Thm.~6.1]{DixonMcKee1986}. For $\beta=0$, part (b) follows from the case $\beta=\frac12$ with $K$ replaced by $KT^{1/2}$, since $1\le T^{1/2}((m-j)\tau)^{-1/2}$ for $j<m\le M$.

\subsection{Inputs proved in the formalization}\label{map:inputs}

\begin{proposition}[Maximal theorem {\cite[Proposition~A.1]{EtoPaperI}}]\label{prop:ext-hl}
For $f\in L^2(\R;\C)$ and $t\in\R$, let $\mathcal M_\R f(t)=\sup_{r>0}\frac1{2r}\int_{t-r}^{t+r}|f(s)|\,\dd s$. There is a constant $A$ such that $\|\mathcal M_\R f\|_{L^2(\R)}\le A\|f\|_{L^2(\R)}$ for every $f\in L^2(\R;\C)$.
\end{proposition}

This follows from \cite[Thm.~10]{HardyLittlewood30} with $k=2$, $a=-\infty$ and $b=\infty$, applied to $|f|$, since $\frac1{2r}\int_{t-r}^{t+r}|f(s)|\,\dd s$ is the mean of the averages of $|f|$ over $[t-r,t]$ and $[t,t+r]$, and both are bounded by the function $\Theta$ of \cite[Thm.~7]{HardyLittlewood30}.

\begin{lemma}[Layer representation of the Dirichlet-to-Neumann map {\cite[Proposition~A.9]{EtoPaperI}}]\label{prop:ext-layer}
Let $x$ be a curve of a family that satisfies \ref{K1}--\ref{K3}, and let $\mathfrak C$ and $r_1$ be as in \cref{sec:bem}. Let $Y$ be a polygon of $\mathfrak C$, or the curve $x$, or a curve $x+\varepsilon e$ with a P1 direction $e$ and $|\varepsilon|\,\|e'\|_{L^\infty}\le2r_1$. If $Y$ is a polygon, let $g\in H^1$ and let $w$ be a P1 function. Otherwise, let $g\in H^2$ and let $w$ be Lipschitz. Let $\varsigma,\varsigma_{\rm int},\varsigma_{\rm ext}\in L^2$ and $c,c'\in\R$ satisfy
\[
\mathcal L_Y\varsigma=g+c,\quad\int_\Tor\varsigma\,\dd\xi=0,\qquad\mathcal L_Y\varsigma_{\rm int}=(\tfrac12+\mathcal K_Y)g,\qquad\mathcal L_Y\varsigma_{\rm ext}=(\tfrac12-\mathcal K_Y)g+c',\quad\int_\Tor\varsigma_{\rm ext}\,\dd\xi=0 .
\]
Then, we have
\[
2\langle\mathcal L_Y\partial_\xi g,\partial_\xi w\rangle+\langle\varsigma_{\rm int},(\tfrac12+\mathcal K_Y)w\rangle+\langle\varsigma_{\rm ext},(\tfrac12-\mathcal K_Y)w\rangle=\int_\Tor\varsigma\,w\,\dd\xi .
\]
\end{lemma}

The lemma is derived from \cref{prop:ext-david,prop:ext-calderon}.

\begin{proposition}[Fredholm alternative for the single layer]\label{prop:ext-slfred}
Let $Y$ be a smooth curve parametrized proportionally to arclength, and assume that $\mathcal L_Y\varsigma=0$ implies $\varsigma=0$ for $\varsigma\in L^2(\Tor)$. Then, for every $g\in H^1(\Tor)$, there is $\varsigma\in L^2(\Tor)$ with $\mathcal L_Y\varsigma=g$, and there is $B$ such that $\|\varsigma\|_{L^2}\le B\|\mathcal L_Y\varsigma\|_{H^1}$ for every $\varsigma\in L^2(\Tor)$.
\end{proposition}

The single layer operator $A$ on a Lipschitz boundary is Fredholm of index zero from $H^{-1/2}(\Gamma)$ to $H^{1/2}(\Gamma)$ \cite[Thm.~2 and p.~616]{Costabel1988}, hence an isomorphism if it is injective \cite[\S5(i)]{Costabel1988}. By the regularity result \cite[Thm.~3 and (2.11)]{Costabel1988} with the regularity index $\frac12$, a density $\psi\in H^{-1/2}(\Gamma)$ with $A\psi\in H^1(\Gamma)$ lies in $L^2(\Gamma)$, with $\|\psi\|_{L^2}\le C(\|A\psi\|_{H^1}+\|\psi\|_{H^{-1/2}})$. Hence, injectivity on $L^2$ gives injectivity on $H^{-1/2}(\Gamma)$, and, for $\psi=A^{-1}g$ with $g\in H^1$, the bound of $A^{-1}$ gives $\|\psi\|_{L^2}\le C\|A\psi\|_{H^1}$. The proposition is this statement in the parameter, as for \cref{prop:ext-slcoer}.

\subsection{Continuum estimates}\label{map:cont}

For a curve $x$, we write $\mathring x=x-\hat x_0$ for the curve with its mean removed. In the next lemma, $G$, $G_1$ and $G_2$ are the kernel $G$ of \cref{sec:app-cont} formed with the curves $x$, $x_1$ and $x_2$ and their lengths $\ell$, $\ell_1$ and $\ell_2$.

\begin{lemma}\label{lem:cont-kernelreg}
Let $\sigma>\tfrac52$ and $x,x_1,x_2\in\mathcal C(\sigma)$.
\begin{enumerate}[label=(\alph*),ref=\thetheorem(\alph*)]
\item For $p\ge0$ with $\mathring x\in H^{p+1/2}$, we have $G\in H^p(\Tor^2)$ and
$\|G\|_{H^p(\Tor^2)}\le C\bigl(1+\|\mathring x\|_{H^{p+1/2}}\bigr)$.
In particular, $\|G\|_{H^{\sigma+3/2}(\Tor^2)}\le C$.
\item If $a_++b_+\le\sigma+\tfrac32$, then $\|\mathcal I_G\|_{H^{-a}\to H^b}\le C$.
\item For $1<p\le\sigma+\tfrac32$, we have
\[
\|G_1-G_2\|_{H^p(\Tor^2)}\le C\bigl(\|\mathring x_1-\mathring x_2\|_{H^{p+1/2}}+|\ell_1-\ell_2|\bigr).
\]
\end{enumerate}
\end{lemma}

\begin{lemma}\label{prop:cont-tame}
Let $\sigma>\tfrac52$. In parts (b)--(e), the curves are as in
\cref{thm:cont-N1}, with the same $\mu$.
\begin{enumerate}[label=(\alph*),ref=\thelemma(\alph*)]
\item Let $x\in\mathcal C(\sigma)$ and $b\ge0$ with $\mathring x\in
H^{b+1/2}$. Then, for $\varphi\in L^2$, we have
\[
\|\SL_G\varphi\|_{H^b}\le C\bigl(1+\|\mathring x\|_{H^{b+1/2}}\bigr)\|\varphi\|_{L^2}.
\]
\item Let $r\ge0$ with $\kc\in H^{r-1/2}$. Then, for $g\in H^1$, we have
\[
\|Q_\Gamma g\|_{H^r}\le C\bigl(1+\|\kc\|_{H^{r-1/2}}\bigr)\|g\|_{H^1},\qquad
\|\Lam_\Gamma g\|_{H^r}\le\frac{4\pi}{\ell}\|g\|_{H^{r+1}}+C\bigl(1+\|\kc\|_{H^{r-1/2}}\bigr)\|g\|_{H^1}.
\]
\item Let $\rho\ge-\tfrac12$ with $\kc\in H^\rho$. Then, for $g\in H^1$, we
have $\|Q_\Gamma g\|_{H^{\rho+1/2}}\le C(1+\|\kc\|_{H^\rho})\|g\|_{H^1}$.
\item\label[lemma]{thm:cont-N3} For two curves $\Gamma_1$ and $\Gamma_2$
and $1<r\le\sigma+\tfrac32$, we have
\[
\|(\Lam_{\Gamma_1}-\Lam_{\Gamma_2})g\|_{H^{r-1}}\le C\bigl(|\ell_1-\ell_2|
+\|\kc_1-\kc_2\|_{H^{\max(r-3/2,0)}}\bigr)\|g\|_{H^r}.
\]
\item\label[lemma]{thm:cont-N3plus} For two curves $\Gamma_1$ and
$\Gamma_2$, $\tfrac12\le t\le\sigma+\tfrac12$ and $1<r\le t+1$, we have
\[
\|(Q_{\Gamma_1}-Q_{\Gamma_2})g\|_{H^t}\le C\bigl(|\ell_1-\ell_2|
+\|\kc_1-\kc_2\|_{H^{t-1/2}}\bigr)\|g\|_{H^r}.
\]
In particular, for $t=\sigma-\tfrac12$, the curvature difference is measured in
$H^{\sigma-1}$.
\end{enumerate}
\end{lemma}

\subsection{The single layer near the exact curves}\label{sec:app-cont-K3}

Perturbations are controlled by the Lipschitz bound of \cite[Lemma~A.7]{EtoPaperI}(i), recalled in \cref{sec:input-results}. For $M\ge c>0$, let $C_L$ be its constant for the class $\mathfrak F(M,c)$, so that for $X,Y\in\mathfrak F(M,c)$ and $\varsigma\in L^2$,
\begin{equation}\label{eq:K-lip}
\|(\mathcal L_X-\mathcal L_Y)\varsigma\|_{H^{1/2}}\le C_L\|X'-Y'\|_{L^\infty}\|\varsigma\|_{H^{-1/2}},\qquad
\|(\mathcal L_X-\mathcal L_Y)\varsigma\|_{H^1}\le C_L\|X'-Y'\|_{L^\infty}\|\varsigma\|_{L^2}.
\end{equation}

\begin{lemma}\label{lem:K3-pert}
Let $X,Y\in\mathfrak F(M,c)$ and $\delta=\|X'-Y'\|_{L^\infty}$. Assume that $Y$ satisfies \ref{K3} with the constants $\mu,B>0$, that is, $\mathcal L_Y$ maps $L^2$ boundedly onto $H^1$, $\langle\mathcal L_Y\varsigma,\varsigma\rangle\ge\mu\|\varsigma\|_{H^{-1/2}}^2$ and $\|\varsigma\|_{L^2}\le B\|\mathcal L_Y\varsigma\|_{H^1}$ for all $\varsigma\in L^2$. If
\[
C_L\delta\le\tfrac12\mu,\qquad C_L\delta B\le\tfrac12,
\]
then $\langle\mathcal L_X\varsigma,\varsigma\rangle\ge\frac12\mu\|\varsigma\|_{H^{-1/2}}^2$ and $\|\varsigma\|_{L^2}\le2B\|\mathcal L_X\varsigma\|_{H^1}$ for all $\varsigma\in L^2$.
\end{lemma}

\begin{lemma}\label{lem:K3-reparam}
Let $y$ be a smooth regular closed curve, let $\phi=\phi_y$ be its normalized arclength map of \cref{sec:ms-flow}, and let $\hat y=y\circ\chi_y$, so that $y=\hat y\circ\phi$. Assume $a\le\phi'\le A$ with $0<a\le A$, and assume that $\hat y$ satisfies \ref{K3} with constants $\mu$ and $B$. Then, $y$ satisfies \ref{K3} with constants $\mu'$ and $B'$ that depend only on $\mu$, $B$, $a$ and $A$.
\end{lemma}

\begin{lemma}\label{lem:K3-exact}
Assume (E). There are $c_E>0$, $\mu_0>0$ and $B_0$ such that every exact curve $\hat X(\cdot,t)$, $0\le t\le T$, satisfies $|\hat X(\xi,t)-\hat X(\eta,t)|\ge c_E\,d(\xi,\eta)$ and \ref{K3} with the constants $\mu_0$ and $B_0$. These constants depend only on (E). Since $\ell(t)\le\ell_{\max}=\max_t\ell(t)$, the chord bound also reads $|\hat X(\xi,t)-\hat X(\eta,t)|\ge(c_E/\ell_{\max})\,\ell(t)\,d(\xi,\eta)$.
\end{lemma}

\begin{lemma}\label{lem:cont-near}
Assume (E), let $\mu_0$ be the constant of \cref{lem:K3-exact}, and let $M'$ and $c'>0$ be such that $\hat X(\cdot,t)\in\mathfrak
F(M',c')$ for $0\le t\le T$. There is $\delta>0$, depending only on (E), $M'$ and $c'$, with the following property.
If a smooth curve $x\in\mathfrak F(M',c')$ with $|x'|$ constant satisfies
$\|x'-\partial_\xi\hat X(\cdot,t)\|_{L^\infty}\le\delta$ for some $t\in[0,T]$,
then $\langle\mathcal L_x\varsigma,\varsigma\rangle\ge\tfrac12\mu_0\|\varsigma\|^2_{H^{-1/2}}$
for all $\varsigma\in L^2$.
\end{lemma}

The formal counterpart of this lemma asserts only some constant $\mu_0>0$, and only for $M'$ and $c'$ beyond bounds fixed by the exact solution. The statement above follows directly from \cref{lem:K3-exact} and \eqref{eq:K-lip}. With $Y=\hat X(\cdot,t)$ and the constant $C_L$ of \eqref{eq:K-lip} for $\mathfrak F(M',c')$, the coercivity of $\mathcal L_Y$ in \cref{lem:K3-exact} and the first bound of \eqref{eq:K-lip} give $\langle\mathcal L_x\varsigma,\varsigma\rangle\ge(\mu_0-C_L\delta)\|\varsigma\|_{H^{-1/2}}^2$, so $\delta=\mu_0/(2C_L)$ has the property.

\subsection{The time-semidiscrete scheme}\label{sec:app-semi}

By \cref{lem:semi-exact}(d) below and \cref{lem:K3-exact}, the exact curves $\hat X(\cdot,t)$ satisfy $\|\kc(t)\|_\sigma\le M_E$, $\ell(t)\in[\ell_E^-,\ell_E^+]$ with $\ell_E^->0$, $|\hat X(\xi,t)-\hat X(\eta,t)|\ge c_E\,d(\xi,\eta)\ge c_E'\,\ell(t)\,d(\xi,\eta)$ with $c_E'=c_E/\max_t\ell(t)$, and the hypothesis of \cref{thm:cont-N1} with the constant $\mu_0$ of \cref{lem:K3-exact}.

\begin{lemma}\label{lem:semi-close}
There are constants $M_+\ge M_E$, $0<c_+\le c_E'$, $0<\ell_-\le\ell_E^-$, $\ell_+\ge\ell_E^+$, $\mu_+>0$, $\varepsilon_0>0$ and $C$, depending only on $\sigma$ and (E), such that the exact curves satisfy the bounds of $\mathcal U_{\rm ad}=\mathcal U_{\rm ad}(M_+,c_+,\ell_-,\ell_+,\mu_+)$ and the following holds. Let $x$ be a smooth closed curve in normalized arclength with length $\ell_x$ and curvature $\kc_x$, and let $0\le t\le T$ with $\|\kc_x-\kc(t)\|_\sigma+|\ell_x-\ell(t)|\le\varepsilon_0$. Then, $x\in\mathcal U_{\rm ad}$, and there is a unit complex number $u$ with
\[
\|u\,x'-\partial_\xi\hat X(\cdot,t)\|_{L^\infty}\le C\bigl(|\ell_x-\ell(t)|+\|\kc_x-\kc(t)\|_{L^2}\bigr).
\]
\end{lemma}

In the next lemma, $\Gamma_j$ are curves of $\mathcal U_{\rm ad}$ with lengths $\ell_j$ and curvatures $\kc_j$, and $\mathcal E_j=\mathcal E_{\kc_j,\Gamma_j}$.

\begin{lemma}\label{lem:semi-E}
Let $\Gamma\in\mathcal U_{\rm ad}$, and let $f$ satisfy $\|f\|_\sigma\le M$ and $\|f\|_{L^2}\ge\frac12\kc_{\min}$. Then, the following hold.
\begin{enumerate}
\item For $1\le r\le\sigma$, we have
\begin{equation}\label{eq:semi-B}
\|\mathcal E_{f,\Gamma}g\|_{\min(r-1,\,\sigma-3/2)}\le C\|g\|_r .
\end{equation}
In particular, $\mathcal E_{f,\Gamma}$ maps $H^1$ to $L^2$ and $H^\sigma$ to $H^{\sigma-3/2}$.
\item Let $B\ge0$. For $\Gamma_1,\Gamma_2\in\mathcal U_{\rm ad}$ and $\|g\|_\sigma\le B$, we have
\begin{equation}\label{eq:semi-Eshape}
\|(\mathcal E_1-\mathcal E_2)g\|_{\sigma-5/2}\le C\bigl(|\ell_1-\ell_2|+\|\kc_1-\kc_2\|_{\sigma-1}\bigr),
\end{equation}
where $C$ depends also on $B$.
\end{enumerate}
\end{lemma}

Let $V=\partial_t\hat X\cdot\bn$ and $\beta=\partial_t\hat X\cdot\bt$ be the components of the velocity of the parametrization of (E). They are smooth on $\Tor\times[0,T]$, and, by (E), $V=-\frac12\Lam_\Gamma\kc$ for $0<t\le T$.

\begin{lemma}\label{lem:semi-exact}
Assume (E). Then, the following hold on $\Tor\times[0,T]$.
\begin{enumerate}
\item The length satisfies $\ell'(t)=\ell\int_\Tor\kc V\,\dd\xi$, and $\beta(\xi,t)=\beta(0,t)-\ell\int_0^\xi(\kc V-\langle\kc V\rangle)\,\dd\xi'$.
\item The curvature satisfies $\partial_t\kc=-\partial_s^2V-\kc^2V+\beta\,\partial_s\kc$.
\item For $0<t\le T$, we have $V=-\frac12\Lam_\Gamma\kc$ and $\beta=\mathcal T_\kc$. The equation \eqref{eq:exact-curv} holds for $0\le t\le T$.
\item For every $r\ge0$, there is $C_{E,r}$ with $\|\kc(t)\|_r\le C_{E,r}$ and $\|\kc(t)-\kc(t')\|_r+|\ell(t)-\ell(t')|\le C_{E,r}|t-t'|$ for $t,t'\in[0,T]$.
\item For $\tau>0$ and $(m+1)\tau\le T$, define the truncation error $d_E^m$ of the exact flow by
\begin{equation}\label{eq:semi-exact-rec}
\big(I+\tau\mathcal P_{\ell(t_m)}\big)\kc(t_{m+1})=\kc(t_m)+\tau\Big(N_{\rm MS}(\kc(t_m);\Gamma(t_m))-\mathcal E_{\kc(t_m),\Gamma(t_m)}\kc(t_m)+d_E^m\Big).
\end{equation}
Then, $\|d_E^m\|_\sigma\le C\tau$.
\end{enumerate}
\end{lemma}

In the next lemma, as in \cref{sec:app-semi-geom}, the velocity is $v=v_n\bn+v_t\bt$, the updated curve is $w=x+\tau v$, and $a=\bt\cdot\partial_sv$, $b=\bn\cdot\partial_sv$ and $\zeta$ are defined as there.

\begin{lemma}\label{lem:semi-geom}
Let $x\in\mathcal C(\sigma;M,c_{\rm ch},\ell_{\min},\ell_{\max})$ be smooth, and let $\kt$, $v_n$, $v_t$ be smooth functions with $\|\kt\|_\sigma+\|v_n\|_{\sigma-1}+\|v_t\|_\sigma\le B$ and $w''=-\ell^2\kt\bn$. There are $\tau_G>0$ and $C$, depending only on $\sigma$, $M$, $c_{\rm ch}$, $\ell_{\min}$, $\ell_{\max}$ and $B$, such that the following hold for $0<\tau\le\tau_G$, where in \textup{(d)}--\textup{(f)} the threshold may be smaller and depend also on the integer and on $A$ named there.
\begin{enumerate}
\item We have $w'=yP$ with $P=\ell(1+\tau a)-\iu\tau\ell b$, so $|w'|=\ell\zeta$ with $\zeta=((1+\tau a)^2+\tau^2b^2)^{1/2}\ge\frac12$. The equation $w''=-\ell^2\kt\bn$ is equivalent to the two rows $L_\kc v_t=D_\kc v_n$ and $\kt=\kc+\tau(L_\kc v_n+D_\kc v_t)$, and
\begin{equation}\label{eq:semi-rows}
(1+\tau a)'=-\tau\ell\kc b,\qquad (\tau b)'=\ell\bigl(\kc(1+\tau a)-\kt\bigr),\qquad \zeta'=-\tau\ell\kt b/\zeta .
\end{equation}
\item The new curve $w$ has length $\ell^+=\ell\int_\Tor\zeta\,\dd\xi$, and
\begin{equation}\label{eq:semi-length}
0\le\ell^+-\ell-\tau\ell\int_\Tor\kc v_n\,\dd\xi\le\tau^2\ell\int_\Tor b^2\,\dd\xi,\qquad|\ell^+-\ell|\le C\tau .
\end{equation}
Its normalized-arclength representative is $x^+=w\circ\psi$, where $\psi=\chi_w$ is the inverse of $\Phi=\phi_w$, $\Phi'=(\ell/\ell^+)\zeta$ and $\psi(0)=0$. The curvature of $x^+$ is
\begin{equation}\label{eq:semi-newcurv}
\kc^+=(\kt \Theta)\circ\psi,\qquad \Theta=\frac{1+\tau a}{\zeta^3},
\end{equation}
its chord constant is bounded below by $c_{\rm ch}\ell^+-C\tau$, and $x^+(0)=x(0)+\tau v(0)$, $y^+(0)=y(0)P(0)/|P(0)|$, where $|P(0)|=\ell\zeta(0)\ge\frac12\ell>0$ by \textup{(a)}.
\item At a common $\xi$, we have
\begin{equation}\label{eq:semi-expand}
\frac{\kc^+-\kt}{\tau}=-2a\kt-\mathfrak P[a]\,\partial_s\kt+\rho,\qquad\|\rho\|_{L^\infty}\le C\tau .
\end{equation}
\item Let $k\ge0$ be an integer and $A\ge0$. If $\|\kt\|_{k+1}+\|\kc\|_{k+1}\le A$, $\|b\|_k\le A(1+\|\kt\|_{k+2})$, $\|b\|_{k+1}\le A(1+\|\kt\|_{k+3})$ and $\tau(1+\|\kt\|_{k+3})\le1$, then, for $0\le\theta\le1$,
\begin{equation}\label{eq:semi-G3}
\|\kc^+-\kt\|_{k+1+\theta}\le C_{k,A}\,\tau\,(1+\|\kt\|_{k+2})^{1-\theta}(1+\|\kt\|_{k+3})^{\theta}.
\end{equation}
\item Let $r\ge3$ be an integer and $A,\varpi\ge0$. If $\|\kt\|_{r-1}+\|\kc\|_{r-1}\le A$, $\|b\|_{r-2}\le A(1+\|\kt\|_r)$, $\|b\|_{r-1}\le A(1+\|\kt\|_{r+1})+\varpi$ and $\tau(1+\|\kt\|_r)\le1$, then
\begin{equation}\label{eq:semi-metric}
\|\kc^+\|_r\le(1+C_{r,A}\tau)\|\kt\|_r+C_{r,A}\tau\|\kt\|_{r+1}+C_{r,A}\tau(1+\varpi),
\end{equation}
where $C_{r,A}$ and the threshold for $\tau$ do not depend on $\varpi$.
\item Let $n\ge0$ be an integer and $A\ge0$. If $\|v_n\|_{W^{n+1,\infty}}+\|v_t\|_{W^{n+1,\infty}}+\|\kc\|_{W^{n,\infty}}\le A$, then $\|\Phi-\mathrm{id}\|_{W^{n+1,\infty}}\le C_{n,A}\tau$.
\end{enumerate}
\end{lemma}

Parts (a)--(c) and (f) follow, after the change to normalized arclength, from the component equations in the proof of \cite[Lemma~3.2]{Rao2026}, the identities \cite[(3.6)--(3.11), (3.28)--(3.29)]{Rao2026} and the coordinate bound in the proof of \cite[Lemma~4.3]{Rao2026}. In the scheme of \cite{Rao2026}, the right-hand side of the BGN equation \cite[(3.5)]{Rao2026} is the normal velocity, $\partial_s^2x^{m+1}=v_n\bn$, while here it is the new curvature, $\partial_s^2w=-\kt\bn$. The identities therefore hold with $v_n$ replaced by $-\kt$ where it stands for this right-hand side, and with $v_n$ unchanged where it is the normal component of the velocity.

The next lemma is the counterpart of \cite[Lemmas~3.2 and 3.4]{Rao2026} for the operator of order three.

\begin{lemma}\label{lem:semi-step}
Let $\Gamma\in\mathcal U_{\rm ad}$ have length $\ell$ and curvature $\kc$. There are $\tau_R>0$ and $C$ such that, for $0<\tau\le\tau_R$, the resolvent equation
\begin{equation}\label{eq:semi-res-eq}
\kt+\tau\bigl(\mathcal P_\ell+\mathcal E_{\kc,\Gamma}\bigr)\kt=\kc
\end{equation}
has a unique solution $\kt\in H^\sigma$, and the following hold.
\begin{enumerate}
\item For $0\le r\le\sigma$, we have $\|\kt-\kc\|_r\le C\tau^{(\sigma-r)/3}$. In particular, $\|\kt\|_\sigma\le C$.
\item For $0\le\theta\le\frac12$, we have $\|\kt\|_{\sigma+3\theta}\le C\tau^{-\theta}$.
\item For $0\le s\le\sigma$, we have $\|\kt-\kc\|_s\le C\tau\,(1+\|\kc\|_{s+3/2})\|\kt\|_{s+3}$.
\item Let $r\ge3$ be real and $K\ge\|\kc\|_{r-1}$. There are $c>0$, $\tau_K>0$ and $C_K$, depending only on $r$, $K$ and the constants of $\mathcal U_{\rm ad}$, such that, for $\tau\le\tau_K$,
\begin{equation}\label{eq:semi-res}
\|\kt\|_r^2+c\,\tau\|\kt\|_{r+3/2}^2\le(1+C_K\tau)\|\kc\|_r^2+C_K\tau .
\end{equation}
\item The solution $\kt$ is real and smooth, and $\|\kt\|_{L^2}\ge\frac12\kc_{\min}$ for $\tau\le\tau_R$.
\end{enumerate}
\end{lemma}

The velocity of the step is formed from the solution of \cref{lem:semi-step} by $v_n=-\frac12\Lam_\Gamma\kt$ and $v_t=L_\kc^{-1}D_\kc v_n$, and $a$ and $b$ are as in \cref{lem:semi-geom}.

\begin{lemma}\label{lem:semi-vel}
In the setting of \cref{lem:semi-step}, the following hold.
\begin{enumerate}
\item We have $\|v_n\|_{\sigma-1}+\|v_t\|_\sigma+\|a\|_{\sigma-1}+\|b\|_{\sigma-2}\le C$.
\item For $s\ge1$ and $t\ge2$, we have
\[
\begin{gathered}
\|v_n\|_s\le C(1+\|\kt\|_{s+1}+\|\kc\|_{s-1/2}),\qquad\|v_t\|_t\le C(1+\|\kt\|_t+\|\kc\|_{t-1}),\\
\|b\|_s\le C(1+\|\kt\|_{s+2}+\|\kc\|_{s+1/2}).
\end{gathered}
\]
\end{enumerate}
\end{lemma}

Let $\Gamma^m\in\mathcal U_{\rm ad}$, let $\kt^m$ be the solution of \cref{lem:semi-step} on $\Gamma^m$, let $v^m$ be its velocity, and let $\Gamma^{m+1}$ be the curve $x^+$ of \cref{lem:semi-geom}. If $\Gamma^{m+1}\in\mathcal U_{\rm ad}$, let $\kt^{m+1}$ be the solution of \cref{lem:semi-step} on $\Gamma^{m+1}$. The defect $d^m$ is then defined by \eqref{eq:semi-drows}, and the length defect by
\[
d_\ell^m=\frac{\ell^{m+1}-\ell^m}\tau-\ell^m\int_\Tor\kt^mv_n^m\,\dd\xi .
\]

In the next lemma, $\Gamma^j\in\mathcal U_{\rm ad}$ for $j\le m_*$, $m_*\tau\le T$, the curve $\Gamma^{j+1}$ is the curve $x^+$ of \cref{lem:semi-geom} for the step from $\Gamma^j$, and $\kt^j$ is the solution of \cref{lem:semi-step} on $\Gamma^j$. We write
\[
\epsilon^j=\kt^j-\kc(t_j),\qquad\lambda_j=\ell^j-\ell(t_j),\qquad \mathcal E_E^j=\mathcal E_{\kc(t_j),\Gamma(t_j)},\qquad\mathcal P_E^j=\mathcal P_{\ell(t_j)} .
\]

\begin{lemma}\label{lem:semi-stab}
In this setting, let $D\ge0$ and $a_0\ge0$ satisfy
\[
\|d^j\|_{L^2}+|d_\ell^j|+\|d_E^j\|_\sigma\le D\quad(j<m_*),\qquad\|\kc^j-\kt^j\|_{\sigma-1}\le a_0\quad(j\le m_*).
\]
Then, there is $C_S$, depending only on the constants of $\mathcal U_{\rm ad}$, on (E) and on $T$, such that
\[
\|\epsilon^m\|_\sigma+|\lambda_m|\le C_S\bigl(\|\epsilon^0\|_\sigma+|\lambda_0|+a_0+D+\tau\bigr)\qquad(m\le m_*).
\]
\end{lemma}

The uniqueness of the step in \cref{prop:semi-A} and the linearization of the step in \cref{sec:app-p52} use the following variational form of \eqref{eq:semi}. Let $x$ be a smooth regular closed curve, and let $\Jm$ be the rotation of \cref{sec:scheme}, so that $\bn\,\dd s=\Jm x'\,\dd\xi$. The time-semidiscrete step seeks $F=\Ft(x)=x+\tau v$ and $c\in\R$ such that
$\int(F-x)\cdot\bn\,\dd s=0$ and
\begin{equation}\label{eq:p52-step}
\ax(F,\varphi)+\tau^{-1}\bx(F-x,\varphi)=c\int\varphi\cdot\bn\,\dd s
\qquad\text{for all }\varphi,
\end{equation}
with the forms $\ax$, $\bx$ of \cref{sec:scheme}. With
$\partial_s^2F=-\kt\bn$, the step is equivalent to $\kt=c-2\Lp v_n$, that is,
to $v_n=-\tfrac12\Lam_\Gamma\kt$ by \cref{lem:cont-gauge}.

\begin{lemma}[The single layer form of the step]\label{lem:p52-gauge}
Let $(F,\kt)$ be a smooth solution of \eqref{eq:p52-step}, and let $v_n=\tau^{-1}(F-x)\cdot\bn$. Then, $\int v_n\,\dd s=0$, and there is a constant $c_{\SL}$ with $\kt=c_{\SL}-2\SL v_n$.
\end{lemma}

\begin{proof}
The constraint gives $\int v_n\,\dd s=0$, so $\Lp v_n=\SL v_n-\avg{\SL v_n}$, and the claim holds with $c_{\SL}=c+2\avg{\SL v_n}$.
\end{proof}

We use the $\SL$-form, in which $\kt=c_{\SL}-2\SL v_n$. By \cref{lem:p52-saddle}, every smooth solution of \eqref{eq:p52-step} is the solution of the saddle point form below.

\begin{lemma}[Unique solvability of the saddle point form]\label{lem:p52-saddle}
Let $Y$ be a smooth curve with $|Y'|\le M$ and with a chord-arc constant $c>0$. Assume that there is $\mu>0$ with $\langle\mathcal L_Y\varsigma,\varsigma\rangle\ge\mu\|\varsigma\|_{H^{-1/2}}^2$ for every $\varsigma\in L^2$, and that $\mathcal L_Y$ maps $L^2$ onto $H^1$. These hypotheses hold for every curve of $\mathcal U_{\rm ad}$, where the second one follows from \cref{prop:ext-slfred}, since the first one makes $\mathcal L_Y$ injective. In the $\SL$-form of \cref{lem:p52-gauge}, the step from $Y$ is the saddle point problem for $(U,c_Y)\in H^1(\Tor;\R^2)\times\R$,
\[
a_Y(Y+U,\varphi)+2\tau^{-1}\big\langle\mathcal L_Y(U\cdot\Jm Y'),\varphi\cdot\Jm Y'\big\rangle=c_Y\int\varphi\cdot\Jm Y'\,\dd\xi,\qquad\int U\cdot\Jm Y'\,\dd\xi=0,
\]
for all $\varphi\in H^1(\Tor;\R^2)$, where $\langle\cdot,\cdot\rangle$ is the $L^2$ pairing in $\xi$. Its operator $A_Y$ is an isomorphism, so the problem has exactly one solution $X_Y$, and \eqref{eq:p52-step} with $x$ replaced by $Y$ has no more than one solution.
\end{lemma}

Let $\hat x^0=\hat X(\cdot,0)$. As long as $\hat x^m\in\mathcal U_{\rm ad}$, let $\kt^m$ be the solution of \cref{lem:semi-step} on $\hat x^m$, let $v^m$ be its velocity, and let $\hat x^{m+1}$ be the curve $x^+$ of \cref{lem:semi-geom}.

\begin{lemma}\label{prop:semi-A}
Assume (E). There are $\tau_A>0$, $\gamma>0$ and $C$ such that, for $0<\tau\le\tau_A$ and every $m$ with $m\tau\le T$, the curve $\hat x^m$ is defined and belongs to $\mathcal U_{\rm ad}$, and
\begin{equation}\label{eq:semi-A}
\|\kc^m-\kc(t_m)\|_\sigma+\|\kt^m-\kc(t_m)\|_\sigma+|\ell^m-\ell(t_m)|\le C\tau^\gamma .
\end{equation}
Moreover, if $(m+1)\tau\le T$, then the step \eqref{eq:semi} from $\hat x^m$ has a unique smooth, regular and embedded solution, its BGN curvature is $\kt^m$, and the normalized-arclength representative of the new curve is $\hat x^{m+1}$.
\end{lemma}

The proof is the induction of \cite[\S3.6]{Rao2026}, with two changes. First, at this stage we use only the low-regularity bound $O(\tau^{(3-\sigma)/6})$ of the initial error in $H^\sigma$, which one resolvent step from $\kc(0)$ gives. Indeed, the exact initial curve lies in the class $\mathcal U_{\rm ad}$ formed with $\sigma''=(\sigma+3)/2\in(\sigma,3)$ in place of $\sigma$, and \cref{lem:semi-step}(a) on this class, with $r=\sigma$, gives $\|\kt^0-\kc(0)\|_\sigma\le C\tau^{(\sigma''-\sigma)/3}=C\tau^{(3-\sigma)/6}$, while $\kc^0=\kc(0)$ and $\ell^0=\ell(0)$. The first-order bound of the initial error is used for \cref{thm:semi-B}, after the uniform higher-order bounds of \cref{lem:semi-Hr} have been established. Second, the class $\mathcal U_{\rm ad}$ of the new curve is obtained from its closeness to the exact curve, which is proved before the error of the new step is known.

\subsection{The linearized time-semidiscrete step}\label{sec:app-p52}

We fix $m$ with $(m+1)\tau\le T$ and the curve $x=x^m$, with $F=\Ft(x)=x^{m+1}$, $v=v^m=(F-x)/\tau$, $c=c^m$, $\kc=\kc^m$ and $\kt=\kt^m$.

\begin{lemma}[The adapted norm]\label{lem:p52-adnorm}
There are $\tau_X>0$, $c_X>0$ and $C_X$ such that, for $0<\tau\le\tau_X$, $0\le m\le T/\tau$ and $u,u_1\in H^1(\Tor;\R^2)$,
\[
c_X\|u\|_{H^1}\le\xnorm[x^m]{u}\le C_X\|u\|_{H^1},\qquad \xnorm[x^m]{u+u_1}\le\xnorm[x^m]{u}+\xnorm[x^m]{u_1}.
\]
\end{lemma}

The equivalence is the one stated after \cite[(5.8)]{Rao2026}.

\begin{lemma}[The single layer along smooth and piecewise linear directions]\label{lem:p52-dSL}
Let $x$ be a smooth closed curve with chord-arc constant $c_x>0$. The constants below depend only on $c_x$ and on bounds of the first six derivatives of $x$, unless stated otherwise.
\begin{enumerate}
\item\label{item:p52-dSL-a} If $e$ is piecewise linear and $\varsigma\in H^1$, or if $e$ is smooth and $\varsigma$ is a trigonometric polynomial, then $f_{\mathcal L}[e]\varsigma\in H^1$ and $\|f_{\mathcal L}[e]\varsigma\|_{H^1}\le C\|e\|_{H^1}\|\varsigma\|_{H^1}$. For a fixed smooth $\varsigma$, the map $e\mapsto f_{\mathcal L}[e]\varsigma$ extends from smooth directions to a bounded operator from $H^1(\Tor;\R^2)$ to $H^1$, whose norm depends also on bounds of finitely many derivatives of $\varsigma$, and the extension equals $f_{\mathcal L}[e]\varsigma$ for piecewise linear $e$.
\item\label{item:p52-dSL-b} Let $e$ be smooth or piecewise linear, and let $d_e=\|e'\|_{L^\infty}$. For $|\varepsilon|\,d_e\le c_x/8$ and $\varsigma\in L^2$, we have
\[
\big\|\mathcal L_{x+\varepsilon e}\varsigma-\mathcal L_x\varsigma-\varepsilon f_{\mathcal L}[e]\varsigma\big\|_{H^1}\le C\varepsilon^2d_e^2\|\varsigma\|_{L^2}.
\]
\end{enumerate}
\end{lemma}

\begin{lemma}[Differentiability of the step]\label{lem:p52-diff}
For every smooth direction $e$ there is $\varepsilon_e>0$ such that, for $|\varepsilon|<\varepsilon_e$, the step \eqref{eq:semi} from $x+\varepsilon e$ has a unique smooth, regular and embedded solution $\Ft(x+\varepsilon e)$, and $\varepsilon\mapsto\Ft(x+\varepsilon e)(\xi)$ is differentiable at $\varepsilon=0$ for every $\xi$. The derivative $u=D\Ft(x)[e]$ lies in $H^1$, and with $w=(u-e)/\tau$ and some $\delta c\in\R$, it satisfies, for all $\varphi\in H^1(\Tor;\R^2)$,
\begin{equation}\label{eq:p52-lin}
\ax(u,\varphi)+\ax'[e](F,\varphi)+\bx^{\SL}(w,\varphi)+\bx^{\SL\prime}[e](v,\varphi)
=\delta c\int\varphi\cdot\bn\,\dd s+c_{\SL}\int\varphi\cdot\Jm e'\,\dd\xi,
\end{equation}
and $\int w_n\,\dd s=-\int v\cdot\Jm e'\,\dd\xi$. Here, $\bx^{\SL}(w,\varphi)=2\int(\SL w_n)\varphi_n\,\dd s$, $\ax'[e](F,\varphi)=-\int(\bt\cdot\partial_se)\,\partial_sF\cdot\partial_s\varphi\,\dd s$, and
\[
\bx^{\SL\prime}[e](v,\varphi)=\underbrace{2\int\big(f_{\mathcal L}[e]\varsigma_v\big)\varphi_n\,\dd s}_{\text{kernel variation}}+\underbrace{2\int\mathcal L_x(v\cdot\Jm e')\,\varphi_n\,\dd s}_{\text{datum variation}}+\underbrace{2\int(\SL v_n)\,\varphi\cdot\Jm e'\,\dd\xi}_{\text{test variation}} .
\]
The map $e\mapsto D\Ft(x)[e]$ extends uniquely to a bounded linear operator on $H^1(\Tor;\R^2)$. For $e\in H^1(\Tor;\R^2)$, the extension $u$ and the number $\delta c$ are the solution, given by the isomorphism $A_x$ of \cref{lem:p52-saddle}, of \eqref{eq:p52-lin} and the linearized constraint, in which the terms that contain $e$ are bounded linear functionals of $e\in H^1$ and the kernel variation $f_{\mathcal L}[e]\varsigma_v$ is the extension of \cref{lem:p52-dSL}\ref{item:p52-dSL-a} for the smooth density $\varsigma_v$.
\end{lemma}

The differentiation is that of \cite[(5.11)--(5.12)]{Rao2026}.

\begin{lemma}[Reparametrization]\label{lem:p52-cov}
Let $\chi=\chi_x$ and $\hat x=x\circ\chi$.
Then, we have $\Ft(\hat x)=\Ft(x)\circ\chi$, and $D\Ft(\hat x)[e\circ\chi]=\big(D\Ft(x)[e]\big)\circ\chi$ for $e\in H^1(\Tor;\R^2)$. Moreover, $\xnorm[x]{u}=\xnorm[\hat x]{u\circ\chi}$.
\end{lemma}

In the normal row, which is \eqref{eq:p52-lin} tested with $\varphi=\psi\bn$, the terms of order one in $\tau$ are collected in
\[
\mathcal Je=\underbrace{(\partial_s^2v\cdot\bn)\,\partial_se_t}_{\text{curvature}}+\underbrace{\partial_s\big[(\partial_sv\cdot\bn)(\bt\cdot\partial_se)\big]-\kc(\partial_sv\cdot\bt)(\bt\cdot\partial_se)}_{\text{stiffness}}+\underbrace{\Dop\Lop^{-1}\mathcal G_xe}_{\text{tangential}} .
\]
This is the operator $-J_x$ of \cite[(5.15)]{Rao2026}, and it is bounded from $H^1$ to $H^{-1}$ because $v$ is smooth. Applying $\tfrac12\Lam_\Gamma$ to the normal row, with $w_n=(u_n-e_n)/\tau$ and $\mathcal H=\tfrac12\Lam_\Gamma\Sop$, we obtain
\begin{equation}\label{eq:p52-res}
(I+\tau\mathcal H)u_n=e_n+\tau f[e]-\tfrac{\tau^2}{2}\Lam_\Gamma\mathcal Je,
\end{equation}
where
\[
f[e]=\underbrace{\frac12\Lam_\Gamma\bigl[(\kc^2+\kc\kt+\Kop_\kc)e_n
-(\partial_s\kc)e_t\bigr]}_{\text{curvature terms}}\underbrace{-\Lam_\Gamma\bigl(f_{\mathcal L}[e]\varsigma_v\bigr)}_{\text{shape derivative}}\underbrace{-\,g_e}_{\text{frame}}
\]
does not contain $\tau^{-1}$.

\begin{lemma}\label{lem:p52-forcing}
The forcing terms of \eqref{eq:p52-res} satisfy $\|f[e]\|_{L^2}\le C\|e\|_{H^1}$ and
$\|\Lam_\Gamma\mathcal Je\|_{H^{-2}}\le C\|e\|_{H^1}$.
\end{lemma}

In normalized arclength, $\Lam_\Gamma=\Lam_0+Q_\Gamma$ and $\Sop=-\partial_s^2+\kc^2+\Kop_\kc$ give $\mathcal H=\mathcal H_0+\mathcal H_1$ with the Fourier multiplier $\mathcal H_0=\ell^{-3}|\partial_\xi|^3=\tfrac12\Lam_0(-\partial_s^2)$, whose symbol is $b|n|^3$ with $b=(2\pi/\ell)^3$, and $\mathcal H_1=\frac12Q_\Gamma(-\partial_s^2)+\frac12\Lam_\Gamma(\kc^2+\Kop_\kc)$, which is bounded from $H^1$ to $L^2$ because $Q_\Gamma$ is bounded from $H^{-1}$ to $L^2$ (\cref{thm:cont-N1}), $\Lam_\Gamma$ is bounded from $H^1$ to $L^2$, and $\Kop_\kc=D_\kc L_\kc^{-1}D_\kc$ is bounded on $H^1$ because $\kc$ is smooth and $L_\kc$ is invertible (\cref{lem:setting-L}). The next lemma replaces the second bound of \cite[(5.18)]{Rao2026}, and the counterpart of the first bound is the energy bound for $\vartheta$ in the proof of \cref{prop:stab}.

\begin{lemma}[The resolvent of $\mathcal H$]\label{lem:p52-resolvent}
Let $c_b=\max(1,b^{-1})$. If $0<\tau\le1$ and $(1+c_b)\tau^{2/3}\|\mathcal H_1\|_{H^1\to L^2}\le\tfrac12$, then $I+\tau\mathcal H$ is an isomorphism from $H^1$ onto $H^{-2}$, and $\|(I+\tau\mathcal H)^{-1}\|_{H^{-2}\to H^1}\le2^{5/2}c_b\tau^{-1}$.
\end{lemma}

\subsection{The fully discrete step}\label{map:fd}

\begin{lemma}[The discrete derivative]\label{lem:fd-diff-h}
Let $X=\Ih x$ with $h\le h_0$. Then, $F_{h,\tau}$ is Fr\'echet differentiable at $X$, and $DF_{h,\tau}(X)e_h=u_h$ is the solution of \eqref{eq:fd-lin-h}.
\end{lemma}

\bibliographystyle{abbrvnat}
\bibliography{paper8,paperI}

\end{document}